\documentclass[12pt]{amsart}

\usepackage[foot]{amsaddr}
\usepackage[utf8]{inputenc} 
\usepackage{amsmath, amssymb,bm, cases, mathtools, thmtools}
\usepackage{verbatim}
\usepackage{graphicx}\graphicspath{{figures/}}
\usepackage{multicol}
\usepackage{tabularx}
\usepackage[usenames,dvipsnames]{xcolor}
\usepackage{mathrsfs} 
\usepackage{url}
\usepackage{bbm}
\usepackage[%
    minnames=4,maxnames=99,maxcitenames=4,
    style=alphabetic,
    doi=false,url=false,
    firstinits=true,hyperref,natbib,backend=bibtex]{biblatex}
\renewbibmacro{in:}{%
  \ifentrytype{article}{}{\printtext{\bibstring{in}\intitlepunct}}}
\bibliography{time_series_bib} 

\usepackage[colorlinks,citecolor=blue,urlcolor=blue,linkcolor=RawSienna]{hyperref}
\usepackage{hypernat}
\usepackage{datetime}

\DeclareMathAlphabet\EuRoman{U}{eur}{m}{n}
\SetMathAlphabet\EuRoman{bold}{U}{eur}{b}{n}

\usepackage{euscript,microtype}

\usepackage[capitalize]{cleveref}

\crefname{lemma}{Lemma}{Lemmas}
\crefname{corollary}{Corollary}{Corollaries}
\crefname{theorem}{Theorem}{Theorems}

\makeatletter
\let\reftagform@=\tagform@
\def\tagform@#1{\maketag@@@{\ignorespaces\textcolor{gray}{(#1)}\unskip\@@italiccorr}}
\renewcommand{\eqref}[1]{\textup{\reftagform@{\ref{#1}}}}
\makeatother

\definecolor{WowColor}{rgb}{.75,0,.75}
\definecolor{SubtleColor}{rgb}{0,0,.50}

\newcounter{margincounter}

\declaretheorem[style=plain,numberwithin=section,name=Theorem]{theorem}
\declaretheorem[style=plain,sibling=theorem,name=Lemma]{lemma}

\declaretheorem[style=definition,sibling=theorem,name=Definition]{definition}
\declaretheorem[style=definition,qed=$\triangleleft$,sibling=theorem,name=Example]{example}
\declaretheorem[style=remark,qed=$\triangleleft$,sibling=theorem,name=Remark]{remark}

\numberwithin{theorem}{section}

\def\[#1\]{\begin{align}#1\end{align}}
\def\*[#1\]{\begin{align*}#1\end{align*}}

\newcommand{\Naturals}{\mathbb{N}}

\newcommand{\Reals}{\mathbb{R}}

\newcommand{\upto}{\!\uparrow\!}

\newcommand{\dee}{\mathrm{d}}

\DeclareMathOperator{\GMC}{GMC}
\DeclareMathOperator{\SLC}{SLC}

\DeclareMathOperator*{\newlim}{\mathrm{lim}\vphantom{\mathrm{infsup}}}
\DeclareMathOperator*{\newmin}{\mathrm{min}\vphantom{\mathrm{infsup}}}
\DeclareMathOperator*{\newmax}{\mathrm{max}\vphantom{\mathrm{infsup}}}
\DeclareMathOperator*{\newinf}{\mathrm{inf}\vphantom{\mathrm{infsup}}}
\DeclareMathOperator*{\newsup}{\mathrm{sup}\vphantom{\mathrm{infsup}}}
\renewcommand{\lim}{\newlim}
\renewcommand{\min}{\newmin}
\renewcommand{\max}{\newmax}
\renewcommand{\inf}{\newinf}
\renewcommand{\sup}{\newsup}

\renewcommand{\Pr}{\mathbb{P}}
\def\EE{\mathbb{E}}
\DeclareMathOperator*{\var}{var}

\newcommand{\floor}[1]{\lfloor #1 \rfloor}

\def\bone{\mathbf{1}}
\def\Ind{\bone}

\newcommand{\abs}[1]{\lvert #1 \rvert}

\newcommand{\bigO}{\mathcal{O}}
\newcommand{\Cov}{\textrm{Cov}}
\newcommand{\bW}{\boldsymbol{W}}

\newcommand{\bI}{\boldsymbol{I}}

\calclayout
\usepackage{setspace}
\title[Spectral Inference under Complex Temporal Dynamics]{Spectral Inference under Complex Temporal Dynamics}
\address{Department of Statistical Sciences, University of Toronto, Canada}
\author{Jun Yang}
\author{Zhou Zhou}
\email{jun@utstat.toronto.edu}
\email{zhou@utstat.toronto.edu}

\begin{document}
	
	\begin{abstract}
	We develop a unified theory and methodology for the inference of evolutionary Fourier power spectra for a general class of locally stationary and possibly nonlinear processes. In particular, simultaneous confidence regions (SCR) with asymptotically correct coverage rates are constructed for the evolutionary spectral densities on a nearly optimally dense grid of the joint time-frequency domain. A simulation based bootstrap method is proposed to implement the SCR. The SCR enables researchers and practitioners to visually evaluate the magnitude and pattern of the evolutionary power spectra with asymptotically accurate statistical guarantee. The SCR also serves as a unified tool for a wide range of statistical inference problems in time-frequency analysis ranging from tests for white noise, stationarity and time-frequency separability to the validation for non-stationary linear models.
	\end{abstract}

	\maketitle
	
	\begin{center}
		\begin{minipage}{1\linewidth}
			\setcounter{tocdepth}{1}
			\tableofcontents
		\end{minipage}
	\end{center}
\clearpage

\section{Introduction}\label{section_intro}
It is well known that the frequency content of many real-world stochastic processes evolves over time. Motivated by the limitations of the traditional spectral methods in analyzing non-stationary signals, time-frequency analysis has become one of the major research areas in applied mathematics and signal processing \cite{Cohen1995,Grochenig2001,Daubechies1990}. Based on various models or representations of the non-stationary signal and its time-varying spectra, time-frequency analysis aims at depicting temporal and spectral information simultaneously and jointly. Roughly speaking, there are three major classes of algorithms in time-frequency analysis: linear algorithms such as short time Fourier transforms (STFT) and wavelet transforms \cite{Allen1977,Meyer1992,Daubechies1992}; bilinear time-frequency representations such as the Wigner–Ville distribution and more generally the Cohen's class of bilinear time–frequency distributions \cite{Cohen1995,Hlawatsch1992} and nonlinear algorithms such as the empirical mode decomposition method \cite{Huang1998} and the synchrosqueezing transform \cite{Daubechies2011}. Though there exists a vast literature on defining and estimating the time-varying frequency content, statistical inference such as confidence region construction and hypothesis testing has been paid little attention to in time-frequency analysis.
 
{It is clear that the subject and the goals of time-frequency analysis and non-stationary time series analysis are highly overlapped.} Unfortunately it seems that the non-stationary spectral domain theory and methodology in the time series literature have been developed largely independently from time-frequency analysis. One major effort in non-stationary time series analysis lies in forming general classes of non-stationary time series models through their evolutionary spectral representation. Among others, \citet{Priestley1965} proposed the notion of evolutionary spectra in a seminar paper. In another seminal work, \citet{Dahlhaus1997} defined a general and theoretically tractable class of locally stationary time series models based on their time-varying spectral representation.  \citet{Nason2000} studied a class of locally stationary time series from an evolutionary wavelet spectrum perspective and investigated the estimation of the latter spectrum.  A second line of research in the non-stationary spectral domain literature involves adaptive estimation of the evolutionary spectra. See for instance \cite{Adak1998} for a binary segmentation based method,  \cite{Ombao2001} for an automatic estimation procedure based on the smooth localized complex exponential (SLEX) transform and  \cite{Fryzlewicz2006} for a Haar–Fisz technique for the estimation of the evolutionary wavelet spectra. On the statistical inference side,  there exists a small number of papers utilizing the notion of evolutionary spectra to test some properties, especially second order stationarity, of a time series. See for instance \cite{Paparoditis2010,Dette2011,Dwivedi2011,Jentsch2015} for tests of stationarity based on properties of the Fourier periodogram or spectral density. See also \cite{Nason2013} for a test of stationarity based on the evolutionary wavelet spectra. On the other hand, however, to date there have been no results on the joint and simultaneous inference of the evolutionary spectrum itself for general classes of non-stationary and possibly nonlinear time series to the best of our knowledge.

The purpose of the paper is to develop a unified theory and methodology for the joint and simultaneous inference of the evolutionary spectral densities for a general class of locally stationary and possibly nonlinear processes. From a time-frequency analysis perspective, the purpose of the paper is to provide a unified and asymptotically correct method for the simultaneous statistical inference of the STFT-based evolutionary power spectra, one of the most classic and fundamental algorithms in time-frequency analysis. Let $\{X^{(N)}_{i}\}_{i=1}^N$ be the observed time series or signal. One major contribution of the paper is that we establish a maximum deviation theory for the STFT-based spectral density estimates over a nearly optimally dense grid ${\mathcal G}_N$ in the joint time-frequency domain. Here the optimality of the grid refers to the best balance between computational burden and (asymptotic) correctness in depicting the overall time-frequency stochastic variation of the estimates. We refer the readers to \cref{subsection_nearlyoptimal} for a detailed definition and discussion of the optimality. The theory is established for a very general class of possibly nonlinear locally stationary processes which admit  a time-varying physical representation in the sense of \cite{Zhou2009} and serves as a foundation for the joint and simultaneous time-frequency inference of evolutionary spectral densities. Specifically, we are able to prove that the spectral density estimates on ${\mathcal G}_N$ are asymptotically independent quadratic forms of $\{X_{i}^{(N)}\}_{i=1}^N$.  Consequently, the maximum deviation of the spectral density estimates on  ${\mathcal G}_N$ behaves asymptotically like a Gumbel law.  The key technique used in the proofs is a joint time-frequency Gaussian approximation to a class of diverging dimensional quadratic forms of non-stationary time series,  which may have wider applicability in evolutionary power spectrum analysis. 

A second main contribution of the paper is that we propose a simulation based bootstrap method to implement simultaneous statistical inferences to a wide range of problems in time-frequency analysis.  The motivation of the bootstrap is to alleviate the slow convergence of the maximum deviation to its Gumbel limit. The bootstrap simply generates independent normally distributed pseudo samples of length $N$ and approximate the distribution of the target maximum deviation with that of the normalized empirical maximum deviations of the spectral density estimates from the pseudo samples. The similar idea was used in, for example \cite{Wu2007,Zhou2010}, {for different problems}. The bootstrap is proved to be asymptotically correct and performs reasonably well in the simulations.  One important application of the bootstrap is to construct simultaneous confidence regions (SCR) for the evolutionary spectral density, which enables researchers and practitioners to visually evaluate the magnitude and pattern of the evolutionary power spectra with asymptotically accurate statistical guarantee. In particular, the SCR helps one to visually identify which variations in time and/or frequency are genuine and which variations are likely to be produced by random fluctuations. See \cref{subsection_real_data} for two detailed applications in earthquake and explosion signal processing and finance. On the other hand, the SCR can be applied to a wide range of tests on the structure of the evolutionary spectra or the time series itself. Observe that typically under some specific structural assumptions, the time-varying spectra can be estimated with a faster convergence rate than those estimated by STFT without any prior information. Therefore a generic testing procedure is to estimate the evolutionary spectra under the null hypothesis and check whether the latter estimated spectra can be fully embedded into the SCR.  This is a very general procedure and it is asymptotically correct as long as the evolutionary spectra estimated under the null hypothesis converges faster than the SCR.  Furthermore, the test achieves asymptotically the power $1$ for local alternatives whose evolutionary spectra deviate from the null hypothesis with a rate larger than the order of the width of the SCR. Specific examples include tests for non-stationary white noise, weak stationarity and time-frequency separability as well as model validation for locally stationary ARMA models and so on. See \cref{subsection_applications} for a detailed discussion and \cref{subsection_real_data} for detailed implementations of the tests in real data. 

Finally, we would like to mention that, under the stationarity assumption, the inference of the spectral density is a classic topic in time series analysis.  There is a vast literature on the topic and we will only list a very small number of representative works. Early works on this topic include \cite{Parzen1957,Woodroofe1967,Brillinger1969,Anderson1971,Rosenblatt1984} among others where asymptotic properties of the spectral density estimates were established under various linearity, strong mixing and joint cumulant conditions.  For recent developments see  \cite{Liu2010,Paparoditis2012,Wu2018} among others.

The rest of the paper is organized as follows. We first formulate the problem in \cref{section_formulation}. In \cref{section_fourier}, we study the STFT and show that the STFTs are asymptotically independent Gaussian random variables under very mild conditions. In \cref{section_consistent_normal}, we study the asymptotic properties of the STFT-based spectral density estimates, including consistency and asymptotic normality. In \cref{section_max_dev}, we establish a maximum deviation theory for the STFT-based spectral density estimates over a nearly optimally dense grid in the joint time-frequency domain. In \cref{section_implement}, we discuss tuning parameter selection and propose a simulation-based bootstrap method to implement the simultaneous statistical inference. Simulations and real data analysis are given in \cref{section_simulation}. Proofs of the main results are deferred to \cref{section_proof} and many details of the proofs have been put in \cref{section_supplement}.
         
\section{Problem Formulation}\label{section_formulation}

We first define locally stationary time series and their instantaneous covariance and spectral density. Throughout the article, we assume the time series $\{X_i^{(N)}\}_{i=1}^N$ is centered, i.e. $\EE[X_i^{(N)}]=0$.  Furthermore, for a random variable $X$, define $\|X\|_q:=[\EE|X|^q]^{1/q}$ and use  $\|\cdot\|$ to denote $\|\cdot\|_2$ for simplicity.
\begin{definition}{(Locally stationary time series \cite{Zhou2009})}\label{def_local_stationary}
	We say $\{X_{i}^{(N)}\}_{i=1}^N$ is a locally stationary time series if there exists a nonlinear filter $G$
	such that
	\[
	X_{i}^{(N)}=G\left(i/N,\mathcal{F}_i\right), \quad i=1,\dots,N,
	\]
	where $\mathcal{F}_i=(\dots,\epsilon_0,\dots,\epsilon_{i-1},\epsilon_i)$ and $\epsilon_i$'s are i.i.d.\ random variables. {Furthermore, the nonlinear filter $G$ satisfies the stochastic Lipschitz continuity condition, $\SLC(q)$, for some $q>0$; that is, there exists $C>0$ such that for all $i$ and $u,s\in (0,1)$,  we have
		\[\label{def_SLC}
		\|G(u,\mathcal{F}_i)-G(s,\mathcal{F}_i)\|_q\le C|u-s|.
		\]}
\end{definition}

\begin{remark}\label{remark_rescale}
	For time series $X_1,X_2, \dots,X_N$, we rescale the time index as $t_i=i/N$, $i=1,\dots,N$. Then $\{t_i\}$ forms a dense grid in $[0,1]$. The rescaled time $u\in [0,1]$ is a natural extension of $\{t_i\}_{i=1}^N$ to be continuum. This rescaling provides an asymptotic device for studying locally stationary time series, which was first introduced by \citet{Dahlhaus1997}. In particular, the rescaling together with the stochastic Lipschitz continuity assumption ensure that for each $X_i$, there is a diverging number of data points in its neighborhood with similar distributional properties.
\end{remark}

\begin{example}(Locally stationary linear time series)
	Let $\epsilon_i$ be i.i.d.\ random variables and
	\[
	G(u,\mathcal{F}_i)=\sum_{j=0}^{\infty} a_j(u) \epsilon_{i-j},
	\]
	where $a_j(u)\in \mathcal{C}^1[0,1]$ for $j=0,1,\dots$. This model was considered in \cite{Dahlhaus1997}. Verification of the SLC assumption is discussed in \cite[Propositions 2 and 3]{Zhou2009}.
\end{example}

\begin{example}(Time varying threshold AR models)
	Let $\epsilon_i\in \mathcal{L}^q, q>0$ be i.i.d.\ random variables with distribution function $F_{\epsilon}$ and density $f_{\epsilon}$. Consider the model
	\[
	G(u,\mathcal{F}_i)=a(u)[G(u,\mathcal{F}_{i-1})]^++b(u)[-G(u,\mathcal{F}_{i-1})]^++\epsilon_i,\quad 0\le u\le 1,
	\]
	where $a(\cdot),b(\cdot)\in\mathcal{C}^1[0,1]$. Then if $\sup_u[|a(u)|+|b(u)|]<1$, the $\SLC(q)$ assumption holds.
	See also \cite[Section 4]{Zhou2009} for more discussions on checking the SLC assumption for locally stationary nonlinear time series.
\end{example}

For simplicity, we will use $X_i$ to denote $X_{i}^{(N)}$ in this paper. Without loss of generality, we assume $X_i=0$ for any $i>N$. We adopt the physical dependence measure \cite{Zhou2009} to describe the dependence structure of the time series.

\begin{definition}(Physical dependence measure)
	Let $\{\epsilon_i'\}$ be an i.i.d.\ copy of $\{\epsilon_i\}$. Consider the locally stationary time series $\{X_i\}_{i=1}^N$. Assume $\max_{1\le i\le N}\|X_i\|_p<\infty$. For $k\ge 0$, define the $k$-th physical dependence measure by
	\[
	\delta_p(k):=\sup_{0\le u\le 1}\|G(u,\mathcal{F}_k)-G(u,(\mathcal{F}_{-1},\epsilon_0',\epsilon_1,\dots,\epsilon_k))\|_p.
	\]
\end{definition}	
Next, we extend the geometric-moment contraction (GMC) condition \cite{Shao2007} to the non-stationary setting. 
\begin{definition}(Geometric-moment contraction)\label{def_GMC}
	We say that the locally stationary time series $\{X_i\}_{i=1}^N$ is $\GMC(p)$ if for any $k$ we have $\delta_p(k)=\bigO(\rho^k)$ for some $\rho\in(0,1)$.
\end{definition}

{Let
	$\mathcal{P}_k(X):=\EE(X\,|\,\mathcal{F}_k)-\EE(X\,|\,\mathcal{F}_{k-1})$
	and $\tilde{X}_{k}^{[\ell]}:=\EE(X_k\,|\,\epsilon_{k-\ell+1},\dots,\epsilon_k)$ be the $\ell$-dependent conditional expectations of $X_k$. 
	From the $\GMC(2)$ condition and $\sup_k \|X_k\|<\infty$, one can easily verify that 
	$\sup_k \sum_{j=-\infty}^{k} \|\mathcal{P}_j X_{k}\|<\infty$ and $\quad \lim_{\ell\to \infty}\sup_k\|X_{k}-\tilde{X}_{k}^{[\ell]}\|=0$.
	We refer to  \cref{remark_GMC} and \cite{Shao2007} for more discussions on the GMC condition.}

\begin{example}(Non-stationary nonlinear time series)\label{example_GMC}
	Many stationary nonlinear time series models are of the form 
	\[
	X_i=R(X_{i-1},\epsilon_i),
	\]
	where $\epsilon_i$ are i.i.d.\ and $R$ is a measurable function. A natural extension to a locally stationary setting is to incorporate the time index $u$ via
	\[
	X_i(u)=R(u,X_{i-1}(u),\epsilon_i),\quad 0\le u\le 1.
	\]
	\citet[Theorem 6]{Zhou2009} showed that one can have a non-stationary process $X_i=X_{i}^{(N)}=G(i/N,\mathcal{F}_i)$ and the $\GMC(\alpha)$ condition holds, if $\sup_u \|R(u,x_0,\epsilon_i)\|_{\alpha}<\infty$ for some $x_0$, and
	\[
	\sup_{u\in [0,1]}\sup_{x\neq y}\frac{\|R(u,x,\epsilon_0)-R(u,y,\epsilon_0)\|_{\alpha}}{|x-y|}<1.
	\]
	See \cite[Section 4.2]{Zhou2009} for more details.
\end{example}

\begin{definition}(Instantaneous covariance)\label{def_inst_covariance}
	Let $u\in [0,1]$. The instantaneous covariance at $u$ is defined by
	\[
	r(u,k):=\Cov\left(G(u,\mathcal{F}_{0}),G(u,\mathcal{F}_{k})\right).
	\]
\end{definition}
\begin{remark}\label{remark_SLC}
	The assumption of $\SLC(q)$ together with $\sup_i \EE|X_i|^p<\infty$, where $1/p+1/q=1$, implies the instantaneous covariance $r(u,k)$ is Lipschitz continuous.
	That is, for all $k$ and for all $u,s\in [0,1], u\neq s$, we have
	\[
	|r(u,k)-r(s,k)|/|u-s|\le C,
	\]
	for some finite constant $C$. The proof is given in \cref{proof_remark_SLC}.
	Therefore, uniformly on $u$, for any positive integer $n\le N$, we have
	\[\label{local_autocorrelation_condition}
	r(u+\delta_u,k)-r(u,k)=\bigO(n/N),\quad \forall - n/N\le \delta_u\le n/N.
	\]
	Particularly, if we choose $n=o(\sqrt{N})$ then $r(u+\delta_u,k)-r(u,k)=o(1/n), \forall - n/N\le \delta_u\le n/N$.
\end{remark}

Next, we define the evolutionary spectral density using the instantaneous covariance.
\begin{definition}(Instantaneous spectral density)\label{def_inst_spectral_density}
	Let $u\in[0,1]$. The spectral density at $u$ is defined by
	\[
	f(u,\theta):=\frac{1}{2\pi}\sum_{k\in\mathbb{Z}}r(u,k)\exp(\sqrt{-1}k\theta).
	\]
\end{definition}

{
\begin{remark}
	In the definition of instantaneous spectral density, $u\in [0,1]$ represents the rescaled time (see \cref{remark_rescale} for more discussions) and $\theta\in [0,2\pi)$ represents the frequency. Different from the usual spectral density for stationary process, the instantaneous spectral density is a two dimensional function of $u$ and $\theta$, which captures the spectral density variation in both time and frequency. The usual spectral density for stationary process is a one-dimensional function of $\theta$ and is static over time. The notion of instantaneous spectral density is useful for capturing the dynamics of the spectral evolution over time. 
\end{remark}
}

{
	\begin{remark}
		Note that, for any fixed time point $u$, $r(u,k)$ is a non-negative definite function on the integers. Hence Bochner's Theorem (or Herglotz Representation Theorem) implies that the covariance function $r(u,k)$ and the spectral density function $f(u,\theta)$ has a one-to-one correspondence at each rescaled time point $u$ under the GMC condition. Therefore, $r(u,k)$ defined in \cref{def_inst_covariance} has a one-one-one correspondence to the spectral density $f(u,\theta)$ defined in \cref{def_inst_spectral_density} for short range dependent locally stationary time series defined in our paper.
	\end{remark}
}

In this paper, we always assume $f_*:=\inf_{u,\theta} f(u,\theta)>0$, which is a natural assumption in the time series literature (see e.g. \cite{Shao2007,Liu2010}). Finally, we define the STFT, the local periodogram, and the STFT-based spectral density estimates.
\begin{definition} (Short-time Fourier transform)
	Let $\tau(\cdot)\le \tau_*<\infty$ be a kernel with support $[-1/2,1/2]$ such that $\tau\in \mathcal{C}^1([-1/2, 1/2])$ and $\int \tau^2(x)\dee x=1$. Let $n$ be the number of data in a local window and $\theta\in [0,2\pi)$. Then the STFT is defined  by 
	\[
	J_n(u,\theta):=\sum_{i=1}^N \tau\left(\frac{i-\lfloor uN\rfloor}{n}\right)X_i \exp(\sqrt{-1}\theta i).
	\]	
\end{definition}
\begin{definition}(Local periodogram)
	\[
	I_n(u,\theta):=\frac{1}{2\pi n}|J_n(u,\theta)|^2.
	\]
\end{definition}

\begin{remark}
	Note that defining
	\[
	\hat{r}(u,k):=\frac{1}{n}\sum_{i=1}^N \tau\left(\frac{i-\lfloor uN \rfloor}{n}\right)\tau\left(\frac{i+k-\lfloor uN \rfloor}{n}\right)X_i X_{i+k},
	\]
	then we can write $I_n(u,\theta)$ as 
	\[
	I_n(u,\theta)=\frac{1}{2\pi}\sum_{k=-n}^n \hat{r}(u,k)\exp(\sqrt{-1}\theta k).
	\]
\end{remark}

It is well known that $I_n(u,\theta)$ is an inconsistent estimator of $f(u,\theta)$ due to the fact that $\hat{r}(u,k)$ are inconsistent when $k$ is large. A natural and classic way to overcome this difficulty is to restrict the above summation to relatively small $k$'s only. This leads to the following.

\begin{definition}{(STFT-based spectral density estimator)}
Let $a(\cdot)$ be an even, Lipschitz continuous kernel function with support $[-1,1]$ and $a(0)=1$; let $B_n$ be a sequence of positive integers with $B_n\to \infty$ and $B_n/n\to 0$. Then the STFT-based spectral density estimator is defined by
\begin{equation}\label{equ_density_estimate}
\hat{f}_n(u,\theta):=\frac{1}{2\pi}\sum_{k=-B_n}^{B_n} \hat{r}(u,k) a(k/B_n)\exp(\sqrt{-1}k\theta).
\end{equation}
\end{definition}

{
	\begin{remark}
		 The modified $\hat{f}_n(u,\theta)$ in \cref{equ_density_estimate} is not always non-negative as it depends on the property of the kernel function $a(\cdot)$. According to \cite[pp.822]{Andrews1991}, if the kernel function further satisfies $\frac{1}{2\pi}\int_{-\infty}^{\infty}a(x)\exp(-\sqrt{-1}\theta x)\dee x\ge 0$ for any $\theta\in [0,2\pi)$, then the modified $\hat{f}_n(u,\theta)$ in \cref{equ_density_estimate} is always non-negative. For example, the Bartlett kernel, $a(x)=(1-|x|)\Ind_{\{|x|\le 1\}}$, and the Parzen kernel, $a(x)=(1-6x^2+6|x|^3)\Ind_{\{0\le |x|\le 1/2\}}+2(1-|x|)^3\Ind_{\{1/2<|x|\le 1\}}$.
	\end{remark}
}

\section{Fourier Transforms}\label{section_fourier}
In this section, we study the STFT and show that the STFTs are asymptotically independent and normally distributed under mild conditions. More specifically, when we consider frequencies $\{2\pi j/n: j=1,\dots,n\}$, we show that uniformly over a grid of $u$ and $j$, $\{J_n(u,2\pi j/n)\}$ are asymptotically independent and normally distributed random variables.
 
Denote the real and imaginary parts of $\{J_n(u,2\pi j/n)/\sqrt{\pi n f(u,2\pi j/n)}\}$ by
\[
\begin{split}
Z_{u,j}^{(n)}&=\frac{\sum_{k=1}^N \tau\left(\frac{k-\lfloor uN\rfloor}{n}\right)X_k\cos(k2\pi j/n)}{\sqrt{\pi n f(u,2\pi j/n)}},\\\quad Z_{u,j+m}^{(n)}&=\frac{\sum_{k=1}^N \tau\left(\frac{k-\lfloor uN\rfloor}{n}\right)X_k\sin(k2\pi j/n)}{\sqrt{\pi n f(u,2\pi j/n)}},\quad j=1,\dots,m,
\end{split}
\]
where $m:=\floor{(n-1)/2}$. Then, we have the following result.

\begin{theorem}\label{thm_fourier} {Assume $\GMC(2)$, $\SLC(2)$, and $\sup_k \EE(X_k^2)<\infty$.} Let $\Omega_{p,q}=\{c\in \mathbb{R}^{pq}: |c|=1\}$, where $|\cdot|$ denotes Euclidean norm, and $$Z_{U,J}=(Z_{u_1,j_1}^{(n)},\dots,Z_{u_1,j_p}^{(n)},\dots,Z_{u_q,j_1}^{(n)},\dots,Z_{u_q,j_p}^{(n)})^T$$ for $J=(j_1,\dots,j_p)$ satisfies $1\le j_1,\dots,j_p\le 2m$  and $U=(u_1,\dots,u_q)$ satisfies $0< u_1<\dots<u_q<1$.
Then for any fixed $p,q\in\mathbb{N}$, as $n\to \infty$, we have that
	\[
	\sup_{J}\sup_{c\in\Omega_{p,q}}\sup_x |P(c^TZ_{U,J}\le x)-\Phi(x)|=o(1),
	\]
	{where $\Phi(x)$ is the cumulative distribution function of the standard normal distribution.}
\end{theorem}	
\begin{proof}
	See \cref{proof_thm_fourier}.
\end{proof}
The above theorem shows that if we select any $p$ elements from the canonical frequencies $\{2\pi j/n, j=1,\dots,n\}$ and $q$ well-separated points from the re-scaled time, the STFTs are asymptotically independent on the latter time-frequency grid. Moreover, the vector formed by these STFTs is asymptotically jointly normally distributed. 

\section{Consistency and Asymptotic Normality}\label{section_consistent_normal}

In this section, we study the asymptotic properties of the smoothed periodogram estimator $\hat{f}_n(u,\theta)$.

\subsection{Consistency}\label{subsection_consistency}
The consistency result for the local spectral density estimate $\hat{f}_n(u,\theta)$ is as follows.
\begin{theorem}\label{thm_consistency}
	{Assume $\GMC(2)$, $\SLC(2)$, and there exists $\delta\in (0,4]$ such that $\sup_i \EE(|X_i|^{4+\delta})<\infty$.} Let $B_n\to \infty$, $B_n=\bigO(n^{\eta})$, $0<\eta<\delta/(4+\delta)$. Then
	\[\label{eq_thm_consistency}
	\sup_{u}\max_{\theta\in [0,\pi]}\sqrt{n/B_n} |\hat{f}_{n}(u,\theta)-\EE(\hat{f}_{n}(u,\theta))|=\bigO_{\Pr}(\sqrt{\log n}).
	\]
\end{theorem}
\begin{proof}
	See \cref{proof_thm_consistency}.
\end{proof}	
Later we will see from \cref{thm_max_dev} that the order $\bigO_{\Pr}(\sqrt{\log n})$ on the right hand side of \cref{eq_thm_consistency} is indeed optimal. 

\begin{remark}\label{remark_expectation_consistency} 
Assume $\sup_i \EE|X_i|^p<\infty$ with $p>4$ and $\SLC(q)$ with $1/p+1/q=1$. If we further assume the kernel $\tau(\cdot)$ is an even function and $r(u,k)$ is twice  continuously differentiable with respect to $u$, then {under $\GMC(2)$}, whenever $n=o(N^{2/3})$, $B_n=o(\min\{n,N^{1/3}\})$, and $\sup_u \sum_{k\in\mathbb{Z}} k^2|r(u,k)|<\infty$, if $a(\cdot)$ is locally quadratic at $0$, i.e.
	\[
	\lim_{u\to 0} u^{-2}[1-a(u)]=C,
	\]
	where $C$ is a nonzero constant, then
	we have
	\[\label{temp_for_consistency}
	\sup_u\sup_{\theta} \left[
	\EE\hat{f}_n(u,\theta)-f(u,\theta)-\frac{C}{B_n^2} f''(u,\theta)\right]=o(1/B_n^2),
	\] 
	where $f''(u,\theta):=-\frac{1}{2\pi}\sum_{k\in\mathbb{Z}}k^2 r(u,k)\exp(\sqrt{-1}k\theta)$.
	The proof is given in \cref{proof_remark_expectation_consistency}. {Therefore, the consistency of $\hat{f}_n(u,\theta)$ is implied by combining \cref{thm_consistency} and \cref{temp_for_consistency}.}
\end{remark}

\subsection{Asymptotic Normality}
Developing an asymptotic distribution for the local spectral density estimate is an important problem in spectral analysis of non-stationary time series. This allows one to perform statistical inference such as constructing point-wise confidence intervals and performing point-wise hypothesis testing. In the following, we derive a central limit theorem for $\hat{f}_n(u,\theta)$. 
\begin{theorem}\label{thm_normality}
	{Assume $\GMC(2)$, $\SLC(2)$, and $\sup_i \EE(|X_i|^{4+\delta})<\infty$ for some $\delta>0$}, $B_n\to \infty$ and $B_n=o(n/(\log n)^{2+8/\delta})$.  Then
	\[
	\sqrt{n/B_n}\{\hat{f}_{n}(u,\theta)-\EE(\hat{f}_{n}(u,\theta))\}\Rightarrow \mathcal{N}(0,\sigma^2_u(\theta)),
	\]
		where $\Rightarrow$ denotes weak convergence,
		$\sigma^2_{u}(\theta)=[1+\eta(2\theta)]f^2(u,\theta)\int_{-1}^{1}a^2(t)\dee t$
		and $\eta(\theta)=1$ if $\theta=2k\pi$ for some integer $k$ and $\eta(\theta)=0$ otherwise.
\end{theorem}
\begin{proof}
	See \cref{proof_thm_normality}.
\end{proof}


\section{Maximum Deviations
}	\label{section_max_dev}
The asymptotic normality for $\hat{f}_n(u,\theta)$ derived in the last section cannot be used to construct simultaneous confidence regions (SCR) over $u$ and $\theta$. For simultaneous spectral inference under complex temporal dynamics, one needs to know the asymptotic behavior of the maximum deviation of $\hat{f}_n(u,\theta)$ from $f(u,\theta)$ on the joint time-frequency domain, which is an extremely difficult problem. In this section, we establish a maximum deviation theory for the STFT-based spectral density estimates over a dense grid in the joint time-frequency domain. Such results serve as a theoretical foundation for the joint time-frequency inference of the evolutionary spectral densities.

\begin{itemize}
	\item{Condition (a):} Define  $\mathcal{U}=\{u_1,\dots,u_{C_n}\}$ where $C_n=|\mathcal{U}|$ and $\frac{n}{2N}< u_i< 1-\frac{n}{2N}, i=1,\dots,C_n$. For any $u_{i_1},u_{i_2}\in \mathcal{U}$ with $i_1\neq i_2$, we {assume that} $|u_{i_1}-u_{i_2}|\ge \frac{n}{N}(1-1/(\log B_n)^2)$.
	\item{Condition (b):} Assume $\sup_k\EE |X_k|^p<\infty$ where $p>4$, and $\SLC(q)$ where $1/p+1/q=1$. Let $\alpha$ be a constant such that $\frac{3}{4(p-1)}<\alpha<\frac{1}{4}$. Then assume $C_n=o[\min\{(nB_n)^{2\alpha(p-1)-1},B_n^{1+2\alpha(p-2)}n^{-2-2\gamma}\}]$ for some $\gamma>0$.
	\item {Condition (c):} Assume that $a(\cdot)$ is an even and bounded function with bounded support $[-1, 1]$, $\lim_{x\to 0}a(x)=a(0)=1$, $\int_{-1}^1 a^2(x)\dee x<\infty$, and $\sum_{j\in\mathbb{Z}}\sup_{|s-j|\le 1}| a(jx)-a(sx)|=\bigO(1)$ as $x\to 0$.
	\item {Condition (d):} There exists $0<\delta_1<\delta_2<1$ and $c_1,c_2>0$ such that for all large $n$, $c_1n^{\delta_1}\le B_n\le c_2 n^{\delta_2}$.
\end{itemize}
Note that Conditions (c) and (d) are very mild. Condition (a) implies that the time interval between any two time points on the grid $\mathcal{U}$ cannot be too close. Condition (b) implies that the total number of the selected time points is not too large.

\begin{remark}
	Condition (a) implies that $C_n\le \frac{N}{n}(1-\frac{n}{N})(1-\frac{1}{(\log B_n)^2})=\bigO(N/n)$. Although we do not assume $\{u_i\}$ to be equally spaced, we suggest in practice choosing $\{u_i\}$ equally spaced and $C_n=\frac{N}{n}(1-\frac{n}{N})(1-\frac{1}{(\log B_n)^2})$ to avoid the tricky problem on how to choose the $u_i$'s and the $C_n$.
\end{remark}

\begin{definition}(Dense Grid $\mathcal{G}_N$)
		Let $\mathcal{G}_N$ be a collection of time-frequency pairs such that $(u,\theta)\in \mathcal{G}_N$ if $u\in \mathcal{U}$ and $\theta\in\{\frac{i\pi}{B_n}, i=0,\dots,B_n\}$.
\end{definition}
The following theorem states that the maximum deviation of the spectral density estimates behaves asymptotically like a Gumbel distribution.
\begin{theorem}\label{thm_max_dev}
	Under $\GMC(2)$ and Conditions (a)--(d), we have that, for any $x\in \Reals$,
	\[\label{eq_max_dev_step1}
	\begin{split}
	&\Pr\left[\max_{(u,\theta)\in\mathcal{G}_N} \frac{n}{B_n}\frac{|\hat{f}_n(u,\theta)-\EE(\hat{f}_n(u,\theta))|^2}{f^2(u,\theta)\int_{-1}^1 a^2(t)\dee t}\right.\\
	&\qquad\left.-2\log B_n-2\log C_n+\log(\pi \log B_n+\pi\log C_n)\le x\right]\to e^{-e^{-x/2}}.
	\end{split}
	\]
\end{theorem}

\begin{proof}
	See \cref{proof_thm_max_dev}.
\end{proof}

\cref{thm_max_dev} states that the spectral density estimates $\hat{f}_n(u,\theta)$ on a dense grid $\mathcal{G}_N$ consisting of $C_n\times B_n$ total number of pairs of $(u,\theta)$ are asymptotically independent quadratic forms of $\{X_i\}_{i=1}^N$. Furthermore, the maximum deviation of the spectral density estimates on $\mathcal{G}_N$ converges to a Gumbel law. This result can be used to construct SCR for the evolutionary spectral densities. Note that \cref{thm_max_dev} is established for a very general class of possibly nonlinear locally stationary processes for the joint and simultaneous time-frequency inference of the evolutionary spectral densities.

\subsection{Near optimality of the grid selection}\label{subsection_nearlyoptimal}
Note that there is a trade-off on how dense the grid should be chosen. On the one hand, we hope the grid is dense enough to asymptotically correctly depict the whole time-frequency stochastic variation of the estimates. On the other hand, making the grid too dense is a waste of computational resources since it does not reveal any extra useful information on the overall variability of the estimates. In the following, we define the notion of asymptotically uniform variation matching of a sequence of dense grids. The purpose of the latter notion is to mathematically determine how dense a sequence of grids should be such that it will adequately capture the overall stochastic variation of the spectral density estimates on the joint time-frequency domain.

\begin{definition}(Asymptotically uniform variation matching of grids) Consider a given sequence of bandwidths $(n,B_n)$,
	and let $\{\tilde{\mathcal{G}}_N\}$ be a sequence of grids of time-frequency pairs $\{(u_i,\theta_j)\}$ with time and frequencies equally spaced i.e.\ $|u_{i+1}-u_{i}|=\delta_{\theta,n}$ and $|\theta_{j+1}-\theta_{j}|=\delta_{u,n}$, respectively. Then the sequence $\{\tilde{\mathcal{G}}_N\}$ is said to be asymptotically uniform variation matching if
	\begin{equation}
	\begin{split}
	\max_{\{u_i,\theta_j\}\in \tilde{\mathcal{G}}_N}\sup_{\{u: |u-u_i|\le \delta_{u,n}, \theta: |\theta-\theta_j|\le \delta_{\theta,n}\}}&\sqrt{n/B_n}\left|\left[\hat{f}_n(u,\theta)-\EE(\hat{f}_n(u,\theta))\right]-\left[\hat{f}_n(u_i,\theta_j)-\EE(\hat{f}_n(u_i,\theta_j))\right]\right|\\
	&=o_{\Pr}(\sqrt{\log n}).
	\end{split}
	\end{equation}
\end{definition}
Note that we have previously shown in \cref{thm_consistency} that the uniform stochastic variation of $\sqrt{n/B_n}\hat{f}_n(u,\theta)$ on $(u,\theta)\in (0,1)\times [0,\pi)$ has the order $\bigO_{\Pr}(\sqrt{\log n})$. In combination with \cref{thm_max_dev}, we can see the order $\bigO_{\Pr}(\sqrt{\log n})$ cannot be improved. Therefore, by a simple chaining argument, we can show if a sequence of grids $\{\tilde{\mathcal{G}}_N\}$ is an asymptotically uniform variation matching, then 
\begin{equation}
\begin{split}
&\sqrt{n/B_n}\left|\sup_{(u,\theta)\in (0,1)\times [0,\pi)}\left|\hat{f}_n(u,\theta)-\EE(\hat{f}_n(u,\theta))\right|-\max_{\{u_i,\theta_j\}\in \tilde{\mathcal{G}}_N}\left|\hat{f}_n(u_i,\theta_j)-\EE(\hat{f}_n(u_i,\theta_j))\right|\right|
=o_{\Pr}(\sqrt{\log n}).
\end{split}
\end{equation}
Hence, the uniform stochastic variation of $\hat{f}_n(u,\theta)$ on $(u,\theta)\in \tilde{\mathcal{G}}_N$ is asymptotically equal to the uniform stochastic variation of $\hat{f}_n(u,\theta)$ on $(u,\theta)\in (0,1)\times [0,\pi)$. In other words, $\max_{\{u_i,\theta_j\}\in \tilde{\mathcal{G}}_N}\left|\hat{f}_n(u_i,\theta_j)-\EE(\hat{f}_n(u_i,\theta_j))\right|$ and  $\sup_{(u,\theta)\in (0,1)\times [0,\pi)}\left|\hat{f}_n(u,\theta)-\EE(\hat{f}_n(u,\theta))\right|$ have the same limiting distribution.

However, a grid that is asymptotically uniform variation matching may be unnecessarily dense which causes a waste of computational resources without depicting any additional useful information. The optimal grid should balance between computational burden and asymptotic correctness in depicting the overall time-frequency stochastic variation of the estimates. {Furthermore, if the grid is too dense, the limiting distribution is different from our main result and is unknown to the best of our knowledge.} Therefore, we hope to choose a sequence of grids as sparse as possible provided it is (nearly) asymptotically uniform variation matching. 

Next, we show the sequence of grids used in \cref{thm_max_dev} is indeed nearly optimal in this sense. Recall that in \cref{thm_max_dev}, the interval between adjacent frequencies is of order $\delta_{\theta,n}=\Omega(1/B_n)$ and the averaged interval between two adjacent time indices is of order $\delta_{u,n}=\Omega(n/N)$, where we define $a_n=\Omega(b_n)$ if $1/a_n=\bigO(1/b_n)$. 
In the following, we show that if we choose a sequence of slightly denser grids with $\delta_{\theta,n}=\bigO\left(\frac{1}{B_n(\log n)^{\alpha}}\right)$ and $\delta_{u,n}=\bigO\left(\frac{n}{N(\log n)^{\alpha}}\right)$ where $\alpha$ is any fixed positive constant, then the latter sequence of grids is asymptotically uniform variation matching. Since $\alpha$ can be chosen arbitrarily close to zero, the dense grids in \cref{thm_max_dev} are nearly optimal.
 
\begin{theorem}\label{thm_near_optimality} Under the assumptions of \cref{thm_max_dev}, a sequence of grids with equally spaced time and frequency intervals $\delta_{u,n}$ and $\delta_{\theta,n}$ is asymptotically uniform variation matching if
	$\delta_{u,n}=\bigO\left(\frac{n}{N(\log n)^{\alpha}}\right)$ and $\delta_{\theta,n}=\bigO\left(\frac{1}{B_n(\log n)^{\alpha}}\right)$ for some $\alpha>0$.
\end{theorem}
\begin{proof}
	See \cref{proof_thm_near_optimality}.
\end{proof}

\subsection{Applications of the Simultaneous Confidence Regions}\label{subsection_applications}
In this subsection, we illustrate several applications of the proposed SCR for joint time-frequency inference. These examples include testing time-varying white noise (\cref{example_test_white}), testing stationarity (\cref{example_test_stationary}), testing time-frequency separability or correlation stationarity (\cref{example_test_correlation}), and validating time-varying ARMA models (\cref{example_test_tvARMA}). 

These examples demonstrate that our maximum deviation theory can serve as a foundation for the joint and simultaneous time-frequency inference. In particular, as far as we know, there is no existing methodology in the literature for testing time-frequency separability of locally stationary time series, nor model validation for time-varying ARMA models, although they are certainly very important problems. On the other hand, our proposed SCR serves as an asymptotically valid and visually friendly tool for the above purposes (see \cref{example_test_correlation,example_test_tvARMA}). 

In order to implement the tests, observe that typically under some specific structural assumptions, the time-varying spectra can be estimated with a faster convergence rate than those estimated by the STFT. Therefore, to test the structure of the evolutionary spectra under the null hypothesis, a generic procedure is to check whether the estimated spectra under the null hypothesis can be fully embedded into the SCR. Note that this very general procedure is asymptotically correct as long as the evolutionary spectra estimated under the null hypothesis converges faster than the SCR. The test achieves asymptotic power $1$ for local alternatives whose evolutionary spectra deviate from the null hypothesis with a rate larger than the order of the width of the SCR.

\begin{example}(Testing time-varying white noise)\label{example_test_white}
	White noise is a collection of uncorrelated random variables with mean $0$ and time-varying variance $\sigma^2(u)$. It can be verified that testing time-varying white noise is equivalent to testing the following null hypothesis:
	\[\label{eq_test_WN}
	H_0:\quad\forall \theta,\quad f(u,\theta)=g(u),\quad u\in [0,1]
	\]  
	for some time-varying function $g(\cdot)$.
	{Consider the following optimization problem:
	\[
	g_0(u):=\arg\min_{\tilde{g}} \frac{1}{\pi}\int_0^{\pi} |f(u,\theta)-\tilde{g}(u)|^2\dee \theta.
	\]
	That is, we would like to find a function of $u$ which is closest to $f(u,\theta)$ in $L_2$ distance. Direct calculations show that $g_0(u)=\frac{1}{\pi}\int_{0}^{\pi}f(u,\theta)\dee\theta$.
	Therefore, under the null hypothesis we can estimate the function $g$ in \cref{eq_test_WN} by
	\[
	\hat{g}(u):=\frac{1}{\pi}\int_{0}^{\pi} \hat{f}_n(u,\theta)\dee \theta\approx \frac{1}{\pi}\int_{0}^{\pi} f(u,\theta)\dee \theta = g_0(u).
	\]
	}
	It can be shown that under the null hypothesis the convergence rate of $\hat{g}(u)$ uniformly over $u$ is $\bigO_{\Pr}(\sqrt{\log n}/\sqrt{n})$, which is faster than the rate of SCR which is $\bigO_{\Pr}(\sqrt{\log n}/\sqrt{n/B_n})$. Therefore, we can apply the proposed SCR to test time-varying white noise.
\end{example}

\begin{example}(Testing stationarity)\label{example_test_stationary} Under the null hypothesis that the time series is stationary, it is equivalent to testing
	\[\label{eq_test_stationary}
	H_0:\quad\forall u,\quad  f(u,\theta)=h(\theta),\quad  \theta\in[0,\pi]
	\]
	for some function $h(\cdot)$.
	{Consider the following optimization problem:
	\[
	h_0(\theta):=\arg\min_{\tilde{h}}	\int_{0}^{1}|f(u,\theta)-\tilde{h}(\theta)|^2\dee u.
	\]
	That is, we would like to find a function of $\theta$ which is closest to $f(u,\theta)$ in $L_2$ distance. Direct calculations show that $h_0(\theta)=\int_0^1 f(u,\theta)\dee u$.
	Therefore, under the null hypothesis, we can estimate the function $h$ in \cref{eq_test_stationary} by 
	\[
	\hat{h}(\theta):=\int_0^1\hat{f}_n(u,\theta)\dee u \approx \int_0^1 f(u,\theta)\dee u = h_0(\theta).
	\]
	}
	It can be shown that the convergence rate of $\hat{h}(\theta)$ uniformly over $\theta$ is $\bigO_{\Pr}(\sqrt{\log n}/\sqrt{N/B_n})$, which is faster than the rate $\bigO_{\Pr}(\sqrt{\log n}/\sqrt{n/B_n})$ of the SCR. Therefore, we can apply the proposed SCR to test stationarity.
\end{example}

\begin{example}(Testing time-frequency separability or correlation stationarity)\label{example_test_correlation}
	We call a non-stationary time series time-frequency separable if $f(u,\theta)=g(u)h(\theta)$ for some functions $g(\cdot)$ and $h(\cdot)$. {If a non-stationary time series is time-frequency separable, the frequency curves across different times are parallel to each other. Similarly, the time curves across different frequencies are parallel to each other as well. Therefore, the property of time-frequency separability enables one to model the temporal and spectral behaviors of the time-frequency function separately. Furthermore,} it can be verified that testing time-frequency separability is equivalent to testing correlation stationarity for locally stationary time series, i.e. $\textrm{corr}(X_i,X_{i+k})=l(k)$, for some function $l(\cdot)$.  Without loss of generality, we can formulate the null hypothesis as
	\[
	H_0:\quad f(u,\theta)=C_0 g(u)h(\theta),
	\]
	for some constant $C_0$ and $\int_0^1 g(u)\dee u=1$ and $\int_0^{\pi}h(\theta)=1$.
%
%
	Under the null hypothesis, we can estimate $C_0$, $g(u)$ and $h(\theta)$ by
	\[
	\hat{C}_0&:=\int_0^{\pi}\int_0^1 \hat{f}_n(u,\theta)\dee u\dee \theta\approx \int_0^{\pi}\int_0^1 f(u,\theta)\dee u\dee \theta=C_0,\\ \hat{g}(u)&:=\frac{1}{\hat{C}_0}\int_0^{\pi}\hat{f}_n(u,\theta)\dee \theta\approx \frac{1}{C_0}\int_0^{\pi} f(u,\theta)\dee\theta=g(u),\\
	\hat{h}(\theta)&:=\frac{1}{\hat{C}_0}\int_0^1 \hat{f}_n(u,\theta)\dee u\approx \frac{1}{C_0}\int_0^1 f(u,\theta)\dee u=h(\theta),
	\] 
	and we can estimate $f(u,\theta)$ by $\hat{C}_0\hat{g}(u)\hat{h}(\theta)$. It can be shown that the convergence rates of $\hat{C}_0$, $\hat{g}(u)$, and $\hat{h}(\theta)$ are $\bigO_{\Pr}(1/\sqrt{N})$, $\bigO_{\Pr}(\sqrt{\log n}/\sqrt{n})$, and $\bigO_{\Pr}(\sqrt{\log n}/\sqrt{N/B_n})$, respectively. All of them are faster than the convergence rate of the SCR which is $\bigO_{\Pr}(\sqrt{\log n}/\sqrt{n/B_n})$.
	Therefore, we can apply the proposed SCR to test the null hypothesis.
\end{example}

\begin{example}(Validating time-varying ARMA models)\label{example_test_tvARMA}  Consider 
the null hypothesis that the time series follows the following time-varying ARMA model
	\[
	H_0: \sum_{i=0}^p a_i(t/N)X_{t-i}=\sum_{j=0}^q b_j(t/N)\epsilon_{t-j}
	\]
	where $a_0(u)=1$, $a_i(\cdot),b_i(\cdot)\in \mathcal{C}^1[0,1]$, and $\epsilon_i$ are uncorrelated random variables with mean $0$ and variance $1$.
	Under the null hypothesis, $\{X_i\}$ is a locally stationary time series with spectral density
	\[
	f(u,\theta)=\frac{1}{2\pi}\frac{\left|\sum_{j=0}^q b_j(u)\exp(\sqrt{-1}2\pi \theta j) \right|^2}{\left|\sum_{i=0}^p a_i(u)\exp(\sqrt{-1}2\pi\theta i)\right|^2}.
	\]
	The spectral density  can be fitted using the generalized Whittle's method \cite{Dahlhaus1997}, where $a_i(t/N)$ and $b_i(t/N)$ are estimated by minimizing a generalized Whittle function and $p$ and $q$ are selected, for example, by AIC. Note that under the null hypothesis, the spectral density estimated using Whittle's method has a convergence rate $\bigO_{\Pr}(\sqrt{\log n}/\sqrt{n})$ which is faster than the rate $\bigO_{\Pr}(\sqrt{\log n}/\sqrt{n/B_n})$ by the STFT-based methods without prior information. Therefore, 
	to test the fitted non-parametric time-varying ARMA model, we can plot the  non-parametric spectral density using the estimated time-varying parameters $a_i(\cdot)$ and $b_i(\cdot)$.  
	Under the null hypothesis, the non-parametric spectral density should fall within our SCR with the prescribed probability asymptotically. 
\end{example}

{
	The benefits of spectral domain approach to various hypothesis testing problems depend on the specific application. For example, for tests of stationarity, the test based on evolutionary spectral density is technically easier than the corresponding tests in the time domain. The main reason is that the time domain test needs to consider time-invariance of $r(u,k)$ for a diverging number of $k$ and hence is a high-dimensional problem. On the other hand, the spectral domain test of stationarity only needs to check that $f(u,\theta)$ does not depend on $u$. Similar arguments apply to the test of white noise. For another example,  we proposed a frequency domain method for the problem of model validation of non-stationary linear models. However, technically it is difficult to approach this problem from the time domain. Furthermore, for many time series signals in engineering applications, the most important information is embedded in the frequency domain. Therefore, in engineering and signal processing applications, frequency domain methods are typically more favourable and are widely used. Therefore frequency-domain-based tests are preferable in many such applications.
}


\section{Bootstrap and Tuning Parameter Selection}\label{section_implement}
In \cref{subsection_bootstrap}, we propose a simulation based bootstrap method to implement simultaneous statistical inferences. The motivation of the bootstrap procedure is to alleviate the slow convergence of the maximum deviation to its Gumbel limit in \cref{thm_max_dev}. We discuss methods for tuning parameter selection in \cref{subsection_selection}.

\subsection{The Bootstrap Procedure}\label{subsection_bootstrap}
Although \cref{thm_max_dev} shows that SCR can be constructed using the Gumbel distribution, the convergence rate in \cref{thm_max_dev} is too slow to be useful in moderate samples. We propose a bootstrap procedure to alleviate the slow convergence of the maximum deviations. One important application of the bootstrap is to construct SCR in moderate sample cases.

Let $\{\epsilon_1,\dots,\epsilon_N\}$ be i.i.d.\ $\mathcal{N}(0,1)$ random variables. Defining 
\[
\hat{r}^{\epsilon}(u,k):=\frac{1}{n}\sum_{i=1}^N \tau\left(\frac{i-\lfloor uN \rfloor}{n}\right)\tau\left(\frac{i+k-\lfloor uN \rfloor}{n}\right)\epsilon_i \epsilon_{i+k}
\]
and 
\[
\hat{f}^{\epsilon}_n(u,\theta):=\frac{1}{2\pi}\sum_{k=-B_n}^{B_n} \hat{r}^{\epsilon}(u,k) a(k/B_n)\exp(\sqrt{-1}k\theta),
\]
it can be easily verified that the following analogy of \cref{thm_max_dev} holds.
\[\label{eq_max_dev_bootstrap}
\begin{split}
&\Pr\left[\max_{(u,\theta)\in\mathcal{G}_N} \frac{n}{B_n}\frac{|\hat{f}^{\epsilon}_n(u,\theta)-\EE(\hat{f}^{\epsilon}_n(u,\theta))|^2}{[f^{\epsilon}(u,\theta)]^2\int_{-1}^1 a^2(t)\dee t}\right.\\
&\qquad\left.-2\log B_n-2\log C_n+\log(\pi \log B_n+\pi\log C_n)\le x\right]\to e^{-e^{-x/2}}.
\end{split}
\]
Therefore, we propose to construct the SCR for $\{\hat{f}_n(u,\theta)\}$ using the empirical distribution of $\hat{f}^{\epsilon}_n(u,\theta)$. More specifically, we generate $\{\epsilon_i\}_{i=1}^N$ independently for $N_{\textrm{MC}}$ times. Let $\bar{f}^{\epsilon}_n(u,\theta)$ be the sample mean of $\{\hat{f}^{\epsilon}_{n,m}(u,\theta), m=1,\dots, N_{\textrm{MC}}\}$ from the $N_{\textrm{MC}}$ Monte Carlo experiments. Then we compute the empirical distribution of
\[
\max_{(u,\theta)\in\mathcal{G}_N}\frac{|\hat{f}^{\epsilon}_{n,m}(u,\theta)-\bar{f}^{\epsilon}_n(u,\theta)|^2}{[\bar{f}^{\epsilon}_n(u,\theta)]^2},\quad m=1,\dots,N_{\textrm{MC}}
\]
to approximate the distribution of
\[
\max_{(u,\theta)\in \mathcal{G}_N}\frac{|f(u,\theta)-\hat{f}_n(u,\theta)|^2}{[\hat{f}_n(u,\theta)]^2},
\]
which can be employed to construct the SCR. For example, for a given $\alpha\in (0,1)$, we estimate the $(1-\alpha)$-th quantile $\gamma^2_{1-\alpha}$ from the bootstrapped distribution using $\hat{f}^{\epsilon}_n(u,\theta)$, which also approximately satisfies
\[
\Pr\left(\max_{(u,\theta)\in\mathcal{G}_N}\frac{|f(u,\theta)-\hat{f}_n(u,\theta)|^2}{[\hat{f}_n(u,\theta)]^2}\le \gamma^2_{1-\alpha}\right)=1-\alpha.
\]
Therefore,  the constructed SCR is
\[\label{tmp_SCR_width}
\max\{0, (1-\gamma_{1-\alpha})\hat{f}(u,\theta)\}\le f(u,\theta)\le(1+\gamma_{1-\alpha})\hat{f}(u,\theta),\quad \forall (u,\theta)\in\mathcal{G}_N.
\]
Note that in small sample cases, the  lower bound for the confidence region can be $0$ if the estimated $\gamma_{1-\alpha}$ is larger than $1$. This happens when $N$ is not large enough and large $B_n$ and $C_n$ are selected. 
For large sample sizes, the estimated $\gamma_{1-\alpha}^2$ is typically much smaller than $1$. In that case, we can further use the following approximation
\[
\frac{|f(u,\theta)-\hat{f}_n(u,\theta)|^2}{[\hat{f}_n(u,\theta)]^2}
\approx [\log (f(u,\theta)/\hat{f}_n(u,\theta)]^2.
\]
Then the SCR can be constructed as
\[
\exp(-\gamma_{1-\alpha})\hat{f}_n(u,\theta)\le f(u,\theta)\le \exp(+\gamma_{1-\alpha})\hat{f}_n(u,\theta),\quad \forall (u,\theta)\in\mathcal{G}_N.
\]

Overall, the practical implementation is given as follows
\begin{enumerate}
	\item Select $B_n$ and $n$ using the tuning parameter selection method described in \cref{subsection_selection};
	\item Compute the critical value using bootstrap described in \cref{subsection_bootstrap};
	\item Compute the spectral density estimates by \cref{equ_density_estimate};
	\item Compute the SCR defined in \cref{subsection_bootstrap} using the spectral density estimates and the critical value obtained by the bootstrap.
\end{enumerate}

{Note that the validity of the proposed bootstrap procedure is asymptotically justified by \cref{thm_max_dev}. On the other hand, theoretical justification for the superiority of the bootstrap procedure for moderate samples is extremely difficult, as it requires deriving higher order asymptotics of the maximum deviation of the time-varying spectral densities. We will investigate this problem in some future work. }

\subsection{Tuning parameter selection}\label{subsection_selection}

Choosing $B_n$ and $n$ in practice is a non-trivial problem. In our Monte Carlo experiments and real data analysis, we find that the minimum volatility (MV) method \cite{Politis1999,Zhou2013} performs reasonably well. Specifically, the MV method uses the fact that the estimator $\hat{f}_n(u,\theta)$ becomes stable when the block size $n$ and the bandwidth $B_n$ are in an appropriate range. More specifically, we first set a proper interval for $n$ as $[n_l, n_r]$. In our simulations and data analysis, {we choose $n_l=2N^{\eta}$ and $n_r=3N^{\eta}$ if $N\le 1000$, $n_l=2.5N^{\eta}$ and $n_r=4N^{\eta}$ if $1000< N\le 2000$ and $n_l=3N^{\eta}$ and $n_r=5N^{\eta}$ if $N>2000$, where $\eta=0.48$. Although the rule for setting $n_l$ and $n_r$ is ad-hoc, it works well in our simulations and data analysis.} In practice, one can also either choose $n_l$ and $n_r$ based on prior knowledge of the data, or select them by visually evaluating the fitted evolutionary spectral densities. A reasonable value of $n$ should not produce too rough or too smooth estimates of the spectral density. {We remark that $n_l$ and $n_u$ are only upper and lower bounds of the candidate bandwidths.} {Simulations show that the simulated coverage probabilities are typically not sensitive to the choices of $n_l$ and $n_r$. 
} 
In order to use the MV method, we first form a two-dimensional grid of all candidate pairs of $(n,B_n)$ such that $n\in [n_l,n_r]$ and $B_n<n/\log(n)$. Then, for each candidate pair $(n,B_n)$, we estimate $\hat{f}_n(u,\theta)$ using the candidate pair for a fixed time-frequency grid of $(u,\theta)$. Next, we compute the average variance of the spectral density estimates $\hat{f}_n(u,\theta)$ over the neighborhood of each candidate pair on the two-dimensional grid of all candidate pairs of $(n,B_n)$. Finally, we choose the pair of $(n,B_n)$ which gives the lowest average variance. We refer to \cite{Politis1999,Zhou2013} for more detailed discussions of the MV method. 

{Note that cross validation is another popular method for choosing bandwidths \cite{Dahlhaus2019}. However, it is a difficult task to implement cross validation in the context of time-varying spectral density estimation. 
Finally, it is well-known that choosing theoretically optimal bandwidths is an extremely difficult problem. We hope to investigate this problem in some future work.
}

\section{Simulations and Data Analysis}\label{section_simulation}
In this section, we study the performance of the proposed SCR via simulations and real data analysis. In \cref{subsection_accuracy_bootstrap}, the accuracy of the proposed bootstrap procedure is studied; The accuracy of tuning parameter selection is considered in \cref{subsection_window};
The accuracy and power for hypothesis testing is studied in \cref{subsection_accuracy_hypo_testing}; Finally, we perform real data analysis in \cref{subsection_real_data}. Throughout this section, the kernel $\tau(\cdot)$ is chosen to be a re-scaled Epanechnikov kernel such that $\int\tau^2(x)\dee x=1$, and the kernel $a(\cdot)$ is a re-scaled tri-cube kernel such that $a(0)=1$. The two kernel functions are defined as follows.
\[
\tau(x):=
\begin{cases}
\frac{\sqrt{30}}{4}(1-4x^2),& \text{if } \abs{x}< 1/2,\\
0,& \text{otherwise},
\end{cases}
\qquad 
a(x):=
\begin{cases}
(1-\abs{x}^3)^3,& \text{if } \abs{x}<1,\\
0,& \text{otherwise}.
\end{cases}
\]
In all the simulations, we ran Monte Carlo experiments for $N_{\textrm{MC}}=10000$. The results for SCR and hypothesis testing are obtained by averaging over $1000$ independent datasets.
\subsection{Accuracy of Bootstrap}\label{subsection_accuracy_bootstrap}
In this subsection, we study the accuracy of the proposed bootstrap procedure for moderate finite samples (e.g. $N=400$ or $N=800$). We consider different examples of locally stationary time series models described in the following  \cref{example_AR,example_WN,example_MS,example_TH,example_BL}. 

\begin{example}\label{example_AR}(Time-varying AR model) 
	 We have 
	\[
	X_i=a(i/N)X_{i-1}+\epsilon_i,
	\]
	where $\{\epsilon_i\}$ are i.i.d.\ $\mathcal{N}(0,1)$. In this example, we choose $a(u)=0.3\cos(2\pi u)$. Then the model is locally stationary in the sense that the AR$(1)$ coefficient $a(u)=0.3\cos(2\pi u)$ changes smoothly on the interval $[0,1]$. The simulated uncoverage probabilities of the SCR are shown in \cref{table_bootstrap_AR}.
	\begin{table}
		\caption{Simulated Uncoverage Probabilities for \cref{example_AR}}\label{table_bootstrap_AR}
		\begin{tabular}{cc|cc|cc}	
	\hline
	&&\multicolumn{2}{c}{$N=400$}  &\multicolumn{2}{c}{$N=800$} \\
	\cline{3-6}
	$n$ & $B_n$ & $\alpha=0.05$ & $\alpha=0.1$& $\alpha=0.05$ & $\alpha=0.1$\\
	\hline
	$72$ & $36$ & $0.03$& $0.06$ & $0.03$ & $0.06$ \\
	$72$ & $32$ & $0.04$ & $0.07$ & $0.04$ & $0.08$ \\	
	$72$ & $28$ & $0.04$ & $0.08$ & $0.05$ & $0.10$ \\	
	\hline
	$54$ & $36$ & $0.03$ & $0.06$ & $0.03$ & $0.06$ \\
	$54$ & $32$ & $0.04$ & $0.08$ & $0.04$ & $0.08$ \\
	$54$ & $28$ & $0.04$ & $0.09$ & $0.05$ & $0.09$ \\
	\hline
	$36$ & $32$ & $0.05$ & $0.11$ & $0.05$ & $0.11$ \\
	$36$ & $28$ & $0.07$ & $0.13$ & $0.07$ & $0.14$ \\
	\hline
\end{tabular}
	\end{table}
\end{example}

\begin{example}(Time-varying ARCH model)\label{example_WN}
	Consider the following time-varying ARCH$(1)$ model: 
	\[
	X_{i}=\epsilon_i\sqrt{a_0(i/N)+a_1(i/N)X_{i-1}^2},
	\] where $\{\epsilon_i\}$ are i.i.d.\ standard normally distributed random variables, $a_0(u)>0, a_1(u)>0$ and $a_0(u)+a_1(u)<1$. Note that $\{X_i\}$ is a white noise sequence. In this example, {we choose $a_0(u)=0.6$ and $a_1(u)=0.3\sin(\pi u)$.}
	The simulated uncoverage probabilities of the SCR are shown in \cref{table_bootstrap_WN}. 
		\begin{table}
		\caption{Simulated Uncoverage Probabilities for \cref{example_WN}}\label{table_bootstrap_WN}
		\begin{tabular}{cc|cc|cc}
	\hline
	&&\multicolumn{2}{c}{$N=400$}  &\multicolumn{2}{c}{$N=800$} \\
	\cline{3-6}
	$n$ & $B_n$ & $\alpha=0.05$  & $\alpha=0.1$ & $\alpha=0.05$ & $\alpha=0.1$ \\
	\hline
	$72$ & $36$ & $0.03$ & $0.06$  & $0.02$ & $0.05$\\
	$72$ & $32$ & $0.04$ & $0.08$ & $0.04$ & $0.08$ \\	
	$72$ & $28$ & $0.05$ & $0.10$ & $0.05$ & $0.11$ \\
	\hline
	$54$ & $36$ & $0.03$ & $0.06$ & $0.03$ & $0.07$ \\
	$54$ & $32$ & $0.05$ & $0.09$ & $0.04$ & $0.09$ \\
	$54$ & $28$ & $0.05$ & $0.10$ & $0.04$ & $0.10$ \\
	\hline
	$36$ & $32$ & $0.06$ & $0.12$ & $0.06$ & $0.12$ \\
	$36$ & $28$ & $0.08$ & $0.15$ & $0.08$ & $0.14$ \\				
	\hline
\end{tabular}
	\end{table}
\end{example}
 
\begin{example}(Time-varying Markov switching model)\label{example_MS}	Suppose $\{S_i\}$ is a Markov chain on state space $\{0,1\}$ with transition matrix $P$. Consider the following time-varying Markov switching model 
	\[
	X_{i}=
	\begin{cases}
	b_1(i/n)X_{i-1}+\epsilon_i ,& \textrm{if}\ S_i=0,\\
	b_2(i/n)X_{i-1}+\epsilon_i,& \textrm{if}\ S_i=1.
	\end{cases}
	\]
	where $\{\epsilon_i\}$ are i.i.d.\ standard normally distributed random variables,  $|b_1|<1$, and $|b_2|<1$. In this example, we choose $P=\begin{bmatrix} 0.9 & 0.1\\ 0.5 & 0.5\end{bmatrix}$, $ b_1(u)=0.4\cos(2\pi u)$, and $b_2(u)=0.1\sin(2\pi u)$. The simulated uncoverage probabilities of the SCR are shown in \cref{table_bootstrap_MS}.
			\begin{table}
		\caption{Simulated Uncoverage Probabilities for \cref{example_MS}}\label{table_bootstrap_MS}
		\begin{tabular}{cc|cc|cc}
	\hline
	&&\multicolumn{2}{c}{$N=400$}  &\multicolumn{2}{c}{$N=800$} \\
	\cline{3-6}
	$n$ & $B_n$ & $\alpha=0.05$ & $\alpha=0.1$ & $\alpha=0.05$ & $\alpha=0.1$ \\
	\hline
	$72$ & $36$ & $0.04$ & $0.09$ & $0.04$ & $0.09$ \\			
	$72$ & $32$ & $0.05$ & $0.11$ & $0.05$ & $0.12$ \\
	$72$ & $28$ & $0.06$ & $0.12$ & $0.06$ & $0.13$ \\
	\hline
	$54$ & $36$ & $0.06$ & $0.10$ & $0.05$ & $0.10$ \\
	$54$ & $32$ & $0.06$ & $0.11$ & $0.06$ & $0.11$ \\
	$54$ & $28$ & $0.06$ & $0.12$ & $0.07$ & $0.13$ \\
	\hline
	$36$ & $32$ & $0.07$ & $0.14$ & $0.08$ & $0.14$ \\
	$36$ & $28$ & $0.07$ & $0.14$ & $0.08$ & $0.15$ \\
	\hline				
\end{tabular}
	\end{table}
\end{example}

\begin{example}(Time-varying threshold AR model)\label{example_TH}
	Suppose $\{\epsilon_i\}$ are i.i.d.\ standard normally distributed random variables and consider the following threshold AR model
	\[
	X_i=a(i/N)\max(0,X_{i-1})+b(i/N)\max(0,-X_{i-1})+\epsilon_i,
	\]
	where $\sup_{u\in [0,1]} [|a(u)|+|b(u)|]<1$.
	In this example, we choose $a(u)=0.3\cos(2\pi u)$ and $b(u)=0.3\sin(2\pi u)$. The simulated uncoverage probabilities of the SCR are shown in \cref{table_bootstrap_TH}. 
			\begin{table}
		\caption{Simulated Uncoverage Probabilities for \cref{example_TH}}\label{table_bootstrap_TH}
		\begin{tabular}{cc|cc|cc}
	\hline
	&&\multicolumn{2}{c}{$N=400$}  &\multicolumn{2}{c}{$N=800$} \\
	\cline{3-6}
	$n$ & $B_n$ & $\alpha=0.05$ & $\alpha=0.1$ & $\alpha=0.05$ & $\alpha=0.1$\\
	\hline
	$72$ & $36$ & $0.05$ & $0.11$ & $0.06$ & $0.12$ \\			
	$72$ & $32$ & $0.06$ & $0.12$ & $0.06$ & $0.13$ \\	
	$72$ & $28$ & $0.07$ & $0.14$ & $0.08$ & $0.15$ \\
	\hline
	$54$ & $36$ & $0.05$ & $0.10$ & $0.06$ & $0.12$ \\
	$54$ & $32$ & $0.05$ & $0.11$ & $0.06$ & $0.12$ \\
	$54$ & $28$ & $0.06$ & $0.12$ & $0.07$ & $0.13$ \\
	\hline
	$36$ & $32$ & $0.08$ & $0.14$ & $0.08$ & $0.14$ \\
	$36$ & $28$ & $0.08$ & $0.14$ & $0.09$ & $0.17$\\					
	\hline
\end{tabular}
	\end{table}
\end{example}

\begin{example}(Time-varying bilinear process)\label{example_BL}
	Let $\{\epsilon_i\}$ be i.i.d.\ standard normally distributed random variables and consider the following model
	\[
	X_i=b(i/N)X_{i-1}+\epsilon_i+c(i/N)X_{i-1}\epsilon_{i-1},
	\]
	where $b^2(u)+c^2(u)<1$.
	In this example, we choose $b(u)=0.3\cos(2\pi u)$ and $c(u)=0.1\sin(2\pi u)$. The simulated uncoverage probabilities of the SCR are shown in \cref{table_bootstrap_BL}.
		\begin{table}
		\caption{Simulated Uncoverage Probabilities for \cref{example_BL}}\label{table_bootstrap_BL}
		\begin{tabular}{cc|cc|cc}
	\hline
	&&\multicolumn{2}{c}{$N=400$}  &\multicolumn{2}{c}{$N=800$} \\
	\cline{3-6}
	$n$ & $B_n$ & $\alpha=0.05$ & $\alpha=0.1$ & $\alpha=0.05$ & $\alpha=0.1$ \\
	\hline
	$72$ & $36$ & $0.04$ & $0.08$ & $0.05$ & $0.08$ \\			
	$72$ & $32$ & $0.05$ & $0.11$ & $0.05$ & $0.10$ \\	
	$72$ & $28$ & $0.06$ & $0.13$ & $0.05$ & $0.11$ \\
	\hline
	$54$ & $36$ & $0.03$ & $0.07$ & $0.04$ & $0.09$ \\
	$54$ & $32$ & $0.05$ & $0.10$ & $0.06$ & $0.11$ \\
	$54$ & $28$ & $0.06$ & $0.13$ & $0.07$ & $0.13$ \\
	\hline
	$36$ & $32$ & $0.06$ & $0.11$ & $0.07$ & $0.14$ \\
	$36$ & $28$ & $0.09$ & $0.16$ & $0.09$ & $0.17$ \\				
	\hline
\end{tabular}
	\end{table}
\end{example}

According to the results in \cref{table_bootstrap_AR,table_bootstrap_WN,table_bootstrap_MS,table_bootstrap_TH,table_bootstrap_BL}, one can see that the proposed bootstrap works well when $B_n$ and $n$ are chosen in a relatively wide range. In the next subsection, we discuss the MV method for selecting $B_n$ and $n$ in practice.

\subsection{Accuracy of Tuning Parameter Selection}\label{subsection_window}
	\begin{table}
	\caption{Simulated Uncoverage Probabilities with Tuning Parameters Selected by the MV Method ({Numbers in the Parentheses Represent the Average} Width of the SCR)}\label{table_window}
	\begin{tabular}{cc|c|c|c|c|c|c|c|c}
		\hline
		\multicolumn{2}{c}{}\vline &\multicolumn{4}{c}{$N=400$} \vline &\multicolumn{4}{c}{$N=800$} \\
		\cline{3-10}
		\multicolumn{2}{c}{} \vline & $n$ & $B_n$ &  $\alpha=0.05$ & $\alpha=0.1$ & $n$ & $B_n$ & $\alpha=0.05$ & $\alpha=0.1$\\
		\cline{1-10}
		\multicolumn{2}{c}{\cref{example_AR} (\cref{table_bootstrap_AR})} \vline & $54$ & $30$ &$0.06$ ($0.51$)& $0.10$ ($0.45$)& $72$ & $30$ &$0.05$ ($0.47$)& $0.09$ ($0.43$)\\
		\multicolumn{2}{c}{\cref{example_WN} (\cref{table_bootstrap_WN})} \vline & $52$ & $32$ &$0.05$ ($0.37$)& $0.09$ ($0.33$)& $69$ & $30$ &$0.04$ ($0.35$)& $0.09$ ($0.30$)\\
		\multicolumn{2}{c}{\cref{example_MS} (\cref{table_bootstrap_MS})} \vline & $54$ & $32$ &$0.06$ ($0.52$)& $0.11$ ($0.47$)& $72$ & $32$ &$0.05$ ($0.49$)& $0.12$ ($0.43$)\\
		\multicolumn{2}{c}{\cref{example_TH} (\cref{table_bootstrap_TH})} \vline & $50$ & $32$ &$0.06$ ($0.53$)& $0.12$ ($0.45$)& $70$ & $32$ &$0.06$ ($0.48$)& $0.12$ ($0.43$)\\
		\multicolumn{2}{c}{\cref{example_BL} (\cref{table_bootstrap_BL})} \vline & $52$ & $32$ &$0.04$ ($0.52$)& $0.09$ ($0.46$)& $69$ & $31$ &$0.06$ ($0.44$)& $0.11$ ($0.40$)\\
		\hline
	\end{tabular}
\end{table}

We apply the MV method described in \cref{subsection_selection} to select the tuning parameters for \cref{example_AR,example_WN,example_MS,example_TH,example_BL}. For all examples, $N=400$ and $N=800$ are considered. The bootstrap accuracy is shown in  \cref{table_window}. 
{Furthermore, according to \cref{tmp_SCR_width}, we also included the average width of the SCR over $(u,\theta)\in\mathcal{G}_N$  in \cref{table_window}.}
{From \cref{table_window}, we can see that the coverage probabilities of the SCR with bandwidths selected by the MV method are accurate. Furthermore, the average width of the SCR decreases as $N$ increases.}

\subsection{Accuracy and Power of Hypothesis Testing}\label{subsection_accuracy_hypo_testing}
In this subsection, we study the accuracy and power of hypothesis testing using the proposed SCR. We consider \cref{example_power_ARCH} for testing stationarity and \cref{example_power_MA} for testing time-varying white noise. Furthermore, we also consider another example of non-parametric ARMA model validation, which is given in \cref{example_parametric_AR}.

\begin{figure}
		\centering
		\begin{minipage}{0.5\textwidth}
			\centering
	\caption{Simulated Powers for Testing Stationarity for \cref{example_power_ARCH}}
	\includegraphics[width=\textwidth]{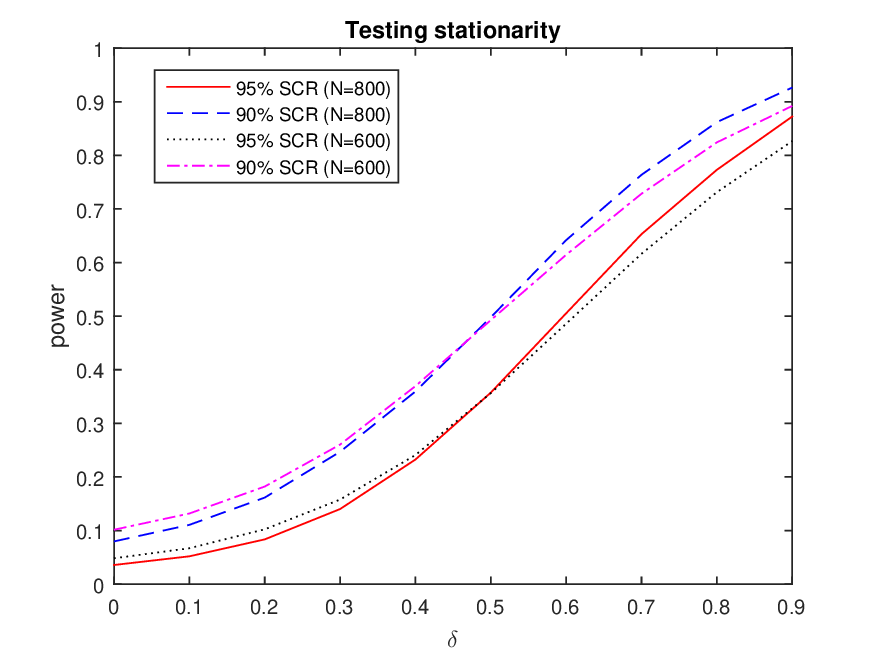}\label{figure_power_ARCH}
		\end{minipage}\hfill
	\begin{minipage}{0.5\textwidth}
		\centering
	\caption{Simulated Powers for Testing TV White Noise for \cref{example_power_MA}}
	\includegraphics[width=\textwidth]{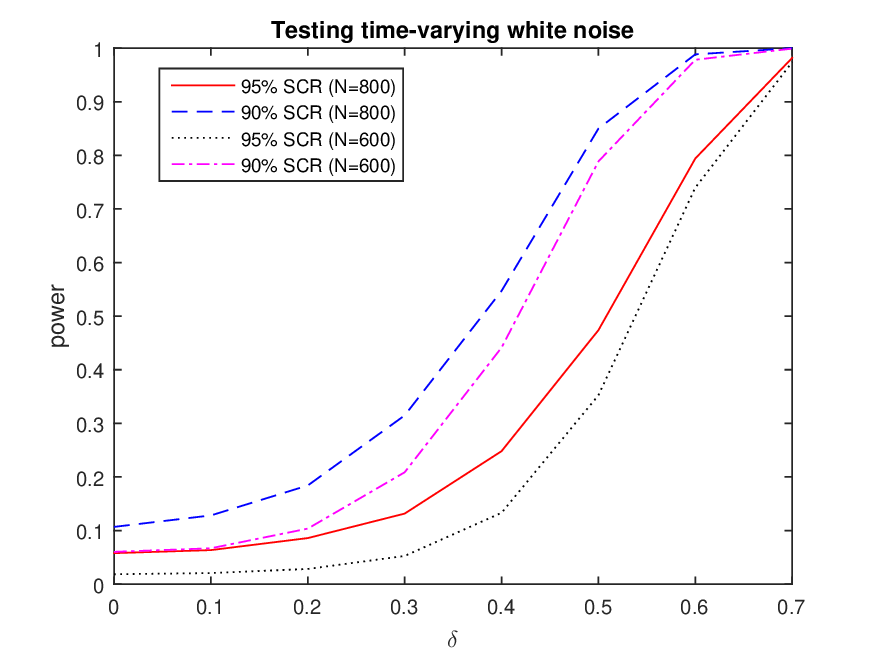}\label{figure_power_MA}
	\end{minipage}
\end{figure}

\begin{example}(Time-varying ARCH model)\label{example_power_ARCH}
	Consider the following model 
	\[
	X_i=\sigma_i\epsilon_i,\quad \sigma^2_i=a_0(i/N)+a_1(i/N) X_{i-1}^2,
	\]
	where $a_0(u)=0.3$ and $a_1(u)=0.2+\delta u$. Observe that when $\delta=0$, the model is stationary. When $\delta=0$, the accuracy of the hypothesis testing for stationarity is studied for two cases, one with $N=400$ and the other with $N=800$, where $n$ and $B_n$ are selected by the MV method. We have shown the simulated Type I error rates of the SCR in \cref{table_hypo}. 
	Next, we study the power of the hypothesis testing for stationarity using the proposed SCR by
	increasing $\delta$. We study both $0.05$ and $0.1$ level tests. The simulated powers of the SCR for $N=800$ {and $N=600$} are shown in \cref{figure_power_ARCH}. {One can see that, for both $N=800$ and $N=600$, the simulated Type I error rates of the SCR (when $\delta=0$) are accurate. Furthermore, the simulated power of the SCR increases with $N$.}

\end{example}

\begin{example}(Time-varying MA model)\label{example_power_MA}
 Consider the following model:
		\[
		X_i=a_0(i/N)\epsilon_i + a_1(i/N)\epsilon_{i-1}
		\] 
	where we let $a_0(u)=0.7+0.9\cos(2\pi u)$ and   $a_1(u)=\delta a_0(u)$. Clearly, when $\delta=0$, the model generates a time-varying white noise. When $\delta=0$, we study the accuracy of the hypothesis testing for time-varying white noise using the proposed SCR. The accuracy by the SCR is shown in \cref{table_hypo}, one with $N=800$ and the other with $N=1200$. The tuning parameters $n$ and $B_n$ are selected by the MV method. 
	We then test time-varying white noise using our proposed SCR by
	increasing $\delta$ for $N=800$ {and $N=600$. The simulated powers of the SCR are shown in \cref{figure_power_MA}.} {According to \cref{figure_power_MA}, one can see that the simulated coverage probabilities for $N=600$ are slightly below the nominal level. This is because the structure of the time series in this example is complicated. A sample size with $N=600$ is not large enough for the local stationarity of the time series to be fully captured statistically.
	On the other hand, for sample size $N=800$, the simulated Type I error rates are accurate and the powers are significantly higher than the case of $N=600$.} 

\end{example}

\begin{table}
	\caption{Simulated Accuracy of Hypothesis Testing}\label{table_hypo}
	\begin{tabular}{cc|c|c|c|c|c|c}
		\hline
		\multicolumn{2}{c}{Nominal Level} \vline &  & $\alpha=0.05$ &  $\alpha=0.1$ &  & $\alpha=0.05$ &  $\alpha=0.1$ \\
		\cline{1-8}
		\multicolumn{2}{c}{\cref{example_power_ARCH}} \vline & $N=400$ & $0.06$ & $0.12$ & $N=800$ & $0.04$ & $0.09$\\
		\multicolumn{2}{c}{\cref{example_power_MA}} \vline & $N=800$ & $0.05$ & $0.10$ & $N=1200$ & $0.04$ & $0.09$\\
		\multicolumn{2}{c}{\cref{example_parametric_AR}} \vline & $N=400$ & $0.04$ & $0.09$ & $N=800$ & $0.06$ & $0.12$\\
		\hline
	\end{tabular}
\end{table}

\begin{example}(Validating time-varying AR model)\label{example_parametric_AR}
	Consider the following time-varying AR model
	\[
	\sum_{j=0}^p a_j(i/N)X_{i-j}=\sigma(i/N)\epsilon_i,
	\]
	where $a_0(u)=1$, $a_j(\cdot)$ and $\sigma(\cdot)$ are smooth functions, $\epsilon_i$ are i.i.d.\ with mean $0$ and variance $1$. Then $\{X_i\}$ is a locally stationary time series with spectral density
	\[
	f(u,\theta)=\frac{\sigma^2(u)}{2\pi}\left|\sum_{j=0}^p a_j(u)\exp(\sqrt{-1}2\pi\theta j)\right|^{-2}.
	\]
	In this example, we generate time series with $p=1$, $a_1(u)=0.3+0.2u$, $\sigma(u)=1+0.3u+0.2u^2$, and length $N=400$ or $N=800$. 
	
	For each generated time series, we fit a time-varying AR model with $p=1$ by minimizing the local Whittle likelihood \cite{Dahlhaus1997}. 
We can then test if the spectral density of the fitted non-parametric time-varying AR model falls into the proposed SCR. The simulated coverage probabilities of the SCR are shown in \cref{table_hypo}, where $n$ and $B_n$ are selected by the MV method. We can see that, under the null hypothesis, the non-parametric time-varying AR model is validated since the simulated Type I error rates match quite well with the nominal levels of the proposed SCR test.
\end{example}


\subsection{Real Data Analysis }\label{subsection_real_data}

In this subsection, we present some real data analysis . We study an earthquake and explosion data set from seismology in \cref{example_earthquake_explosion} and then daily SP500 return from finance in \cref{example_SP500}. Observe that all time series are relatively long with $N>2000$. For tuning parameter selection, we use the MV method to search $(n,B_n)$ within the region $B_n<n/\log(n)$. Hypothesis tests are performed, including testing stationarity, time-varying white noise, and time-frequency separability on all the data sets. 


	\begin{table}
	\caption{Real Data: p-values for Testing (a) Stationarity, (b) Time-Varying White Noise, (c) Time-Frequency Separability (Correlation Stationarity).  }\label{table_real_data_test}
	\begin{tabular}{c|c|c|c}
		\hline
		$H_0$ &{Stationarity}  &{TV White Noise} &{Separability}\\
		\cline{1-4}
		\hline
		Earthquake & $0.0011^{**}$ &  $0.012^{*}$ & $0.064^{+}$ \\	
		Explosion  & $0.0005^{***}$ & $0.033^{*}$ & $0.61$ \\	
		SP500 & $0.0001^{***}$ & $0.99$ & $0.99$ \\	
		SP500 (Abs) & $0.0004^{***}$ & $0.037^{*}$ & $0.048^*$ \\	
		\hline
	\end{tabular}

	Signif. codes: $(***)<0.001\le (**)<0.01\le (*)<0.05\le (+)<0.1$. 
\end{table}

\begin{example}{(Earthquakes and explosions \cite{Shumway2017})}\label{example_earthquake_explosion}
In this example, we study an earthquake signal and an explosion signal from a seismic recording station \cite{Shumway2017}. The recording instruments in Scandinavia are observing earthquakes and mining explosions with one of each shown in \cref{fig_earth} and \cref{fig_explose}, respectively. The two time series (see \cref{fig_earth} and \cref{fig_explose}) each has length $N=2048$ representing two phases or arrivals along the surface, denote by phase $P$: $\{X_i: i=1,\dots,1024\}$ and phase $S$: $\{X_i: i=1025,\dots,2048\}$. The general problem of interest is in distinguishing or discriminating between waveforms generated by earthquakes and those generated by explosions. {The original data came from the technical report by \citet{Blandford1993}. According to \cite{Blandford1993}, the original earthquake and explosion signals have been filtered with a $3$-pole highpass Butterworth filter with the corner frequency at $1$ Hz to improve the signal-to-noise ratio. Then the amplitudes of the waveforms have been rescaled so the maximum amplitude for each signal is equal. According to \cite[Figure 2a and 2b]{Blandford1993}, the unit for time is $0.02$ second and the values of the earthquake and explosion data are rescaled to be no more than $1$.
}

From the time domain (see \cref{fig_earth,fig_explose}), one can observe that rough amplitude ratios of the first phase $P$ to the second phase $S$ are different for the two data sets, which tend to be smaller for earthquakes than for explosions.  From the spectral density estimates and their confidence regions, the $S$ component for the earthquake (see \cref{fig_earth}) shows power at the low frequencies only, and the power remains strong for a long time. In contrast, the explosion (see \cref{fig_explose}) shows power at higher frequencies than the earthquake, and the power of the $P$ and $S$ waves does not last as long as in the case of the earthquake.

Moreover, we notice from the confidence region at selected times and frequencies that the spectral density of explosion has the similar shape at different times, as well as at different frequencies (see \cref{fig_explose_time,fig_explose_freq}). However, the spectral density of earthquakes does not seem to have this property (see \cref{fig_earth_time,fig_earth_freq}). This may suggest that the explosion data are  correlation stationary or time-frequency separable. 
We further perform hypothesis tests on both data sets to confirm our observation (see \cref{table_real_data_test}). The p-values for testing stationarity and time-varying white noise for both earthquake and explosion are quite small, which implies that earthquake and explosion time series are not stationary and not time-varying white noise. However, the p-values for the hypothesis of time-frequency separability (i.e., correlation stationary) is $0.61$ for explosion, but $0.064$ for earthquake. This interesting result discovers a potential important difference between earthquake and explosion:  at least from the analyzed data, explosion tends to be time-frequency separable (correlation stationary) but earthquake does not.

{There are two main benefits from knowing that explosion time series are time-frequency separable but earthquake time series are not. First, this reveals an important structural property of the time-frequency behavior for explosion signals. Since time-frequency separability implies  the time curves for different frequencies are parallel and the frequency curves for different times are parallel as well, this directly suggests a parsimonious model for explosion time series using two one-dimensional models. Second, for the classification of earthquake and explosion signals,  time-frequency separability provides a non-linear feature of the explosion that could potentially serve the purpose. Since most commonly used features for classification are linear features, time-frequency separability is potentially important for feature extraction in order to improve the accuracy in classification tasks. However, since we only have analyzed one pair of earthquake and explosion signals, further studies with a large database of earthquake and explosion signals are needed to confirm this property for explosions which we leave to a future work.}


\begin{figure}
	\centering
	\begin{minipage}{0.5\textwidth}
		\centering
		\caption{Analysis of Earthquake Data}\includegraphics[width=\textwidth]{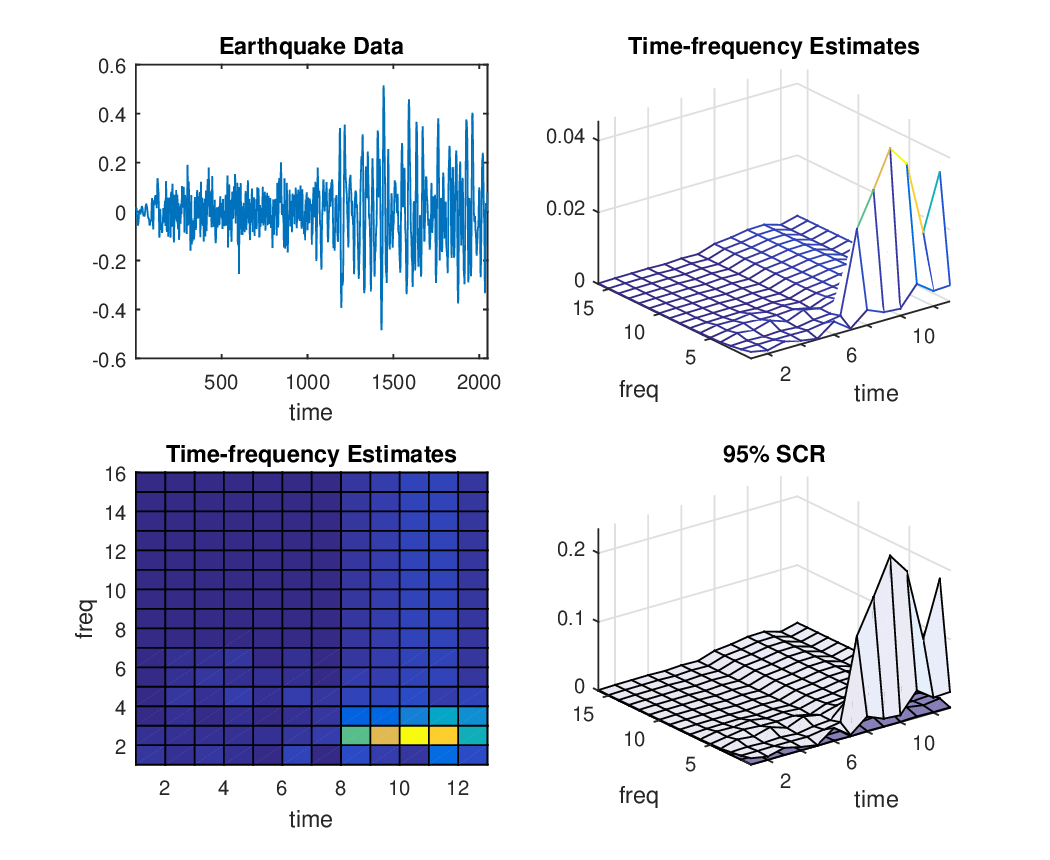}\label{fig_earth}
	\end{minipage}\hfill
	\begin{minipage}{0.5\textwidth}
		\centering
		\caption{Analysis of Explosion Data}\includegraphics[width=\textwidth]{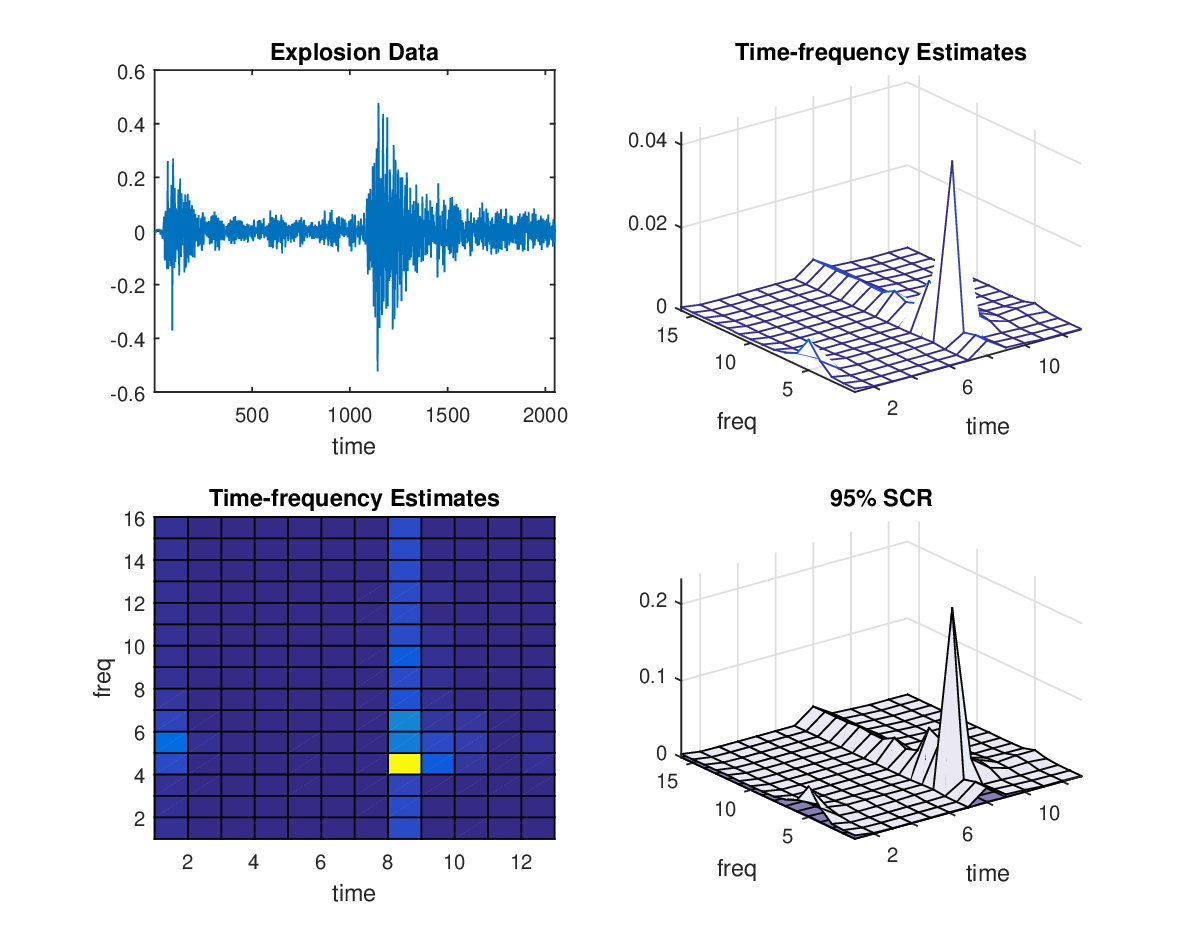}\label{fig_explose}
	\end{minipage}
\end{figure}

\begin{figure}
	\centering
	\begin{minipage}{0.5\textwidth}
		\centering
		\caption{Earthquake Data: Selected Times}\includegraphics[width=\textwidth]{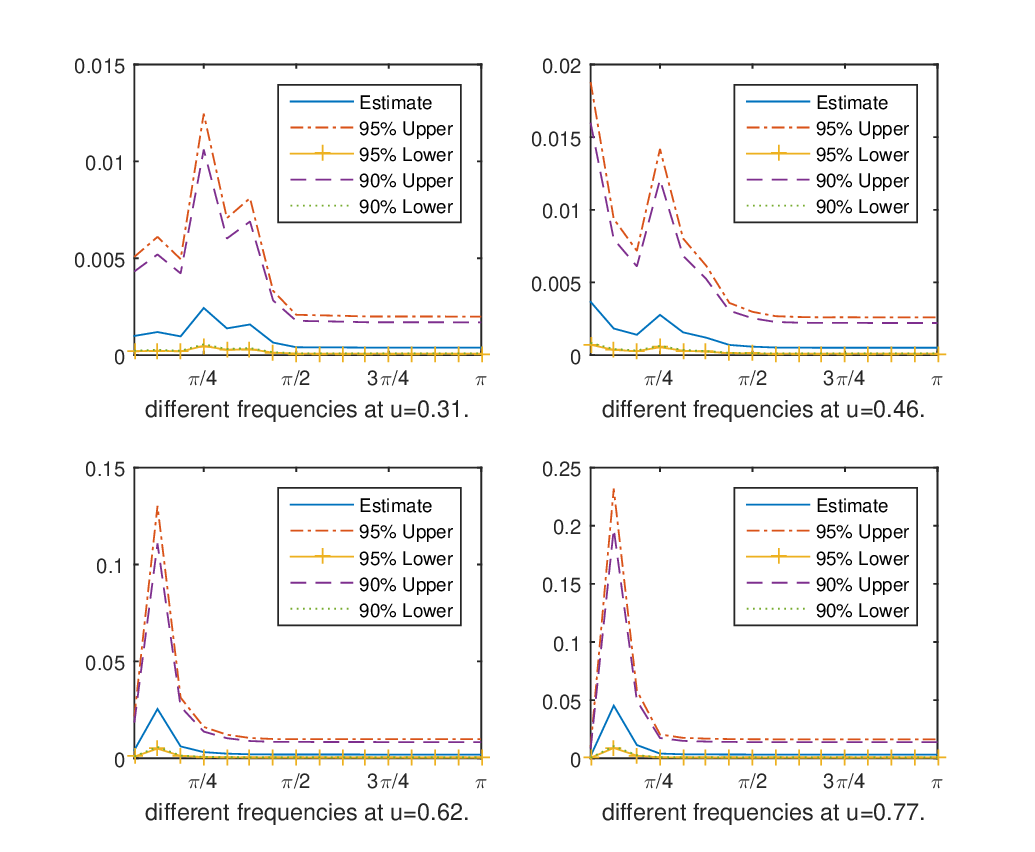}\label{fig_earth_time}
	\end{minipage}\hfill
	\begin{minipage}{0.5\textwidth}
		\centering
			\caption{Earthquake Data: Selected Frequencies}\includegraphics[width=\textwidth]{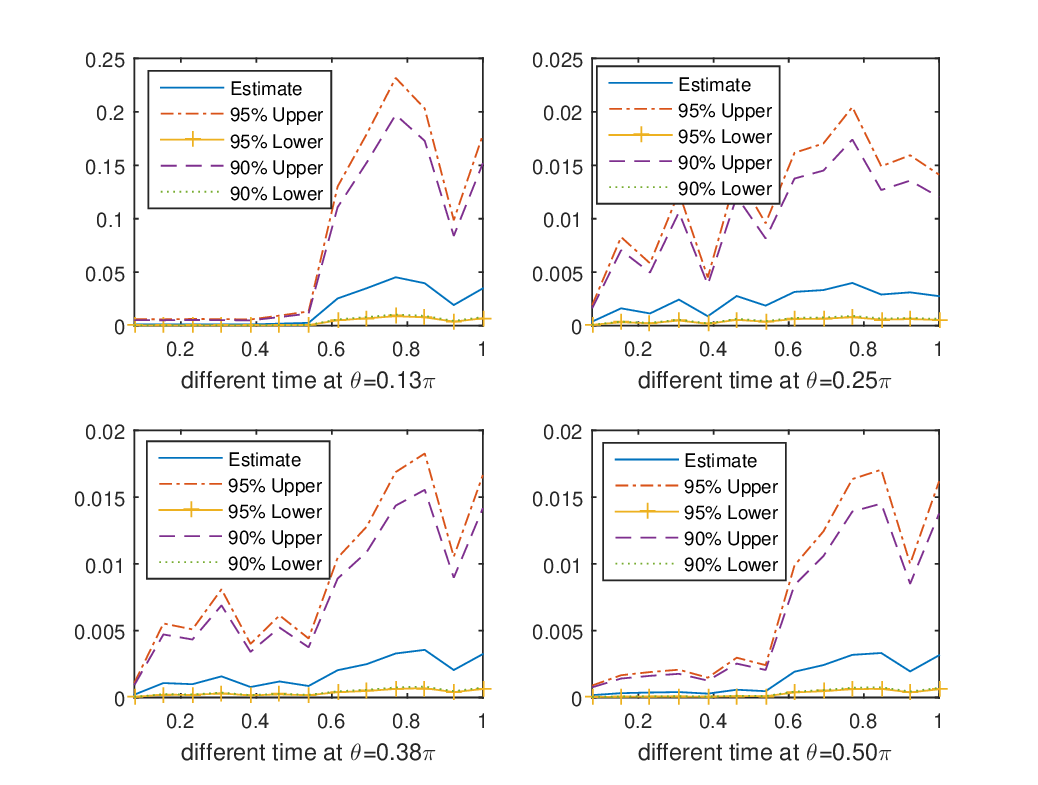}\label{fig_earth_freq}
	\end{minipage}
\end{figure}

\begin{figure}
	\centering
	\begin{minipage}{0.5\textwidth}
		\centering
			\caption{Explosion Data: Selected Times}\includegraphics[width=\textwidth]{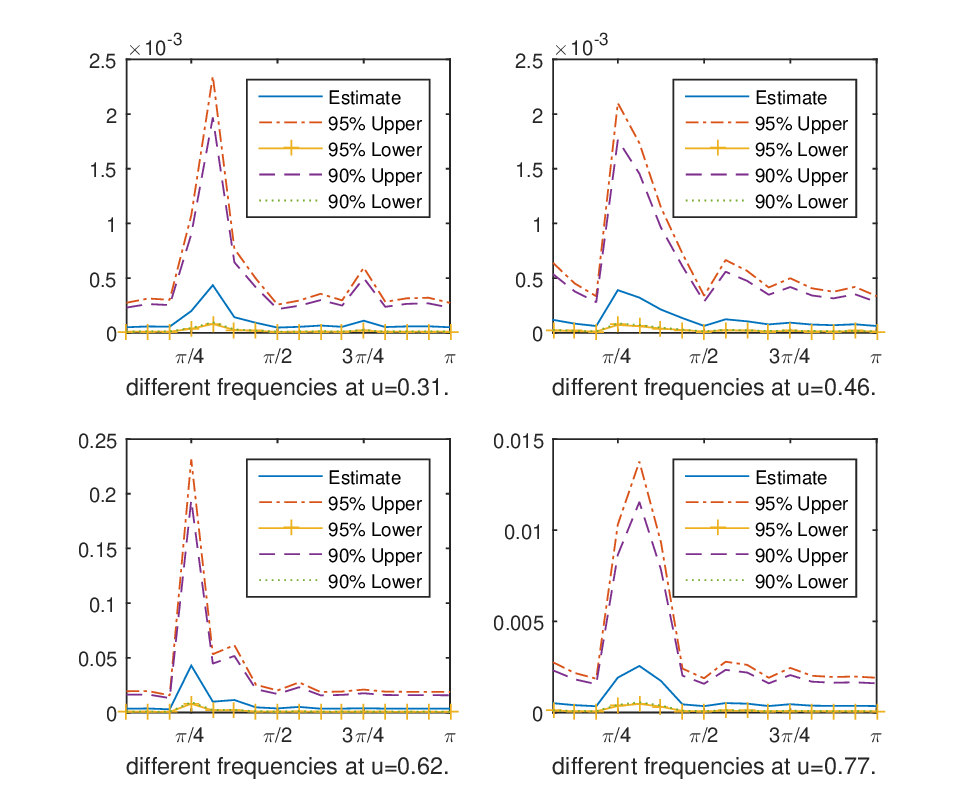}\label{fig_explose_time}
	\end{minipage}\hfill
	\begin{minipage}{0.5\textwidth}
		\centering
			\caption{Explosion Data: Selected Frequencies}\includegraphics[width=\textwidth]{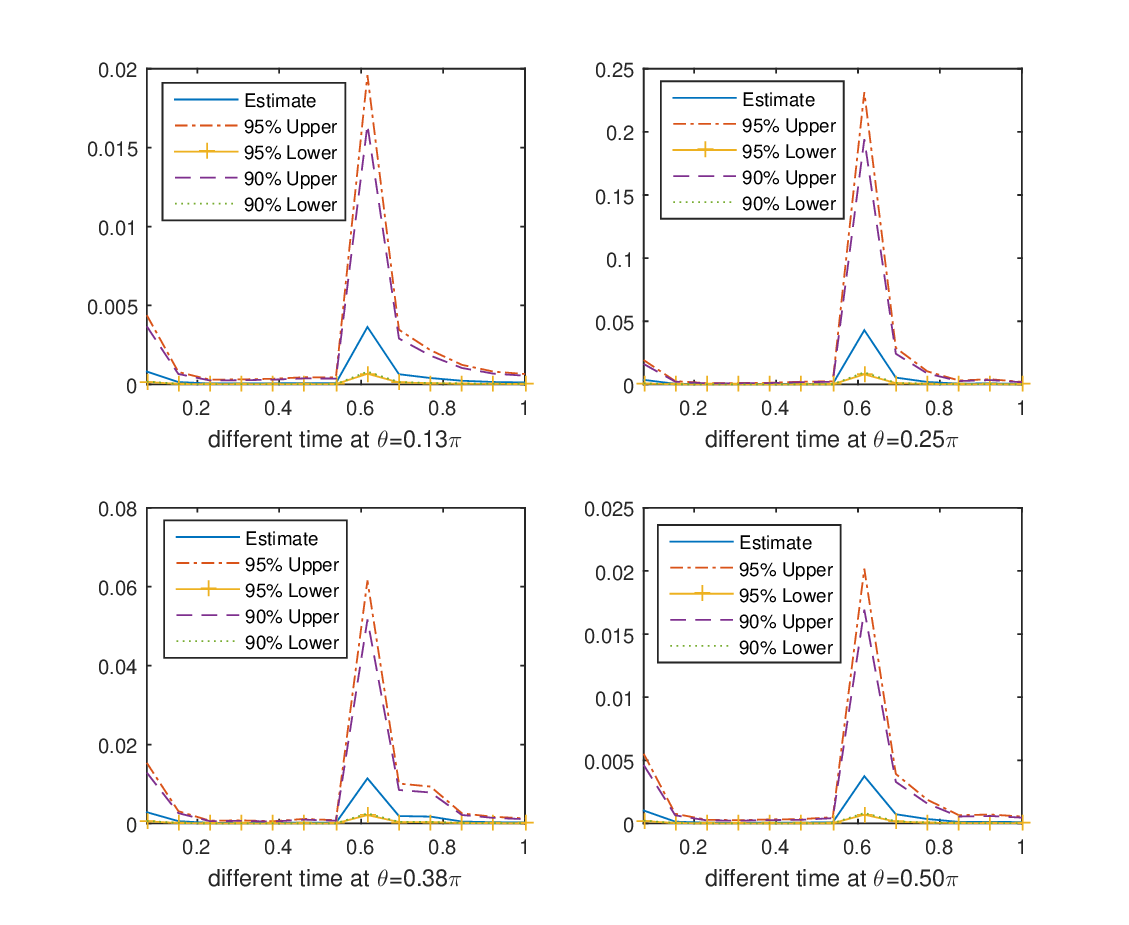}\label{fig_explose_freq}
	\end{minipage}
\end{figure}

\end{example}

\begin{example}{(SP500 daily returns)}\label{example_SP500}
	
In this example, we analyze daily returns of SP500 from September 23rd, 1991 to August 17th, 2018. We plot the original time series, the spectral density estimates and their confidence regions in \cref{fig_SP500}. Observing that the SCR in \cref{fig_SP500} appears to be quite flat over frequencies, it is reasonable to ask if the time series may be modeled as time-varying white noise. Actually, in the finance literature, it is commonly believed that stock daily returns behave like time-varying white noise. We further confirm this observation by performing hypothesis tests. The results (see \cref{table_real_data_test}) show that the SP500 time series is not stationary but it is likely to be a time-varying white noise since the p-value for testing time-varying white noise is $0.99$. Furthermore, the p-value for testing time-frequency separability is also quite large which is $0.99$.

Next, we turn our focus to the absolute value of SP500 daily returns. Volatility forecasting, i.e. forecasting future absolute values or squared values of the return, is a key problem in finance. The celebrated ARCH/GARCH models are equivalent to exponential smoothings of the absolute or squared returns. The optimal weights in the smoothing are determined fully by the evolutionary spectral density. Hence, to optimally forecast the evolutionary volatility, one way is to fit the absolute returns by an appropriate non-stationary linear model, then apply the fitted model to forecast the future volatility. To date, to our knowledge, there exists no methodology for validating non-stationary linear models.  In the following, we demonstrate that the proposed SCR is a useful tool for validating non-stationary linear models for absolute SP500 daily returns.

We first remove the local mean of the original SP500 time series by kernel smoothing.  The spectral density estimates and the SCRs are shown in \cref{fig_SP500abs,fig_SP500abs_time,fig_SP500abs_freq}. We observe from the plots that the spectral density of the absolute SP500 returns behaves quite differently from the original SP500 time series. For example, unlike the case for the original SP500 time series, the SCR for the absolute SP500 in \cref{fig_SP500abs_time} is not flat over frequencies anymore. We perform the same hypothesis tests again to the absolute SP500 time series. The results (see \cref{table_real_data_test}) show that the p-value for testing time-varying white noise is $0.037$, which is much smaller than that of the original SP500 time series. Furthermore, the p-value for testing time-frequency separability is $0.048$ which is also much smaller than the one for the original SP500 data.

\begin{figure}
	\centering
	\begin{minipage}{0.5\textwidth}
		\centering
			\caption{Analysis of Daily Returns of SP500 }\includegraphics[width=\textwidth]{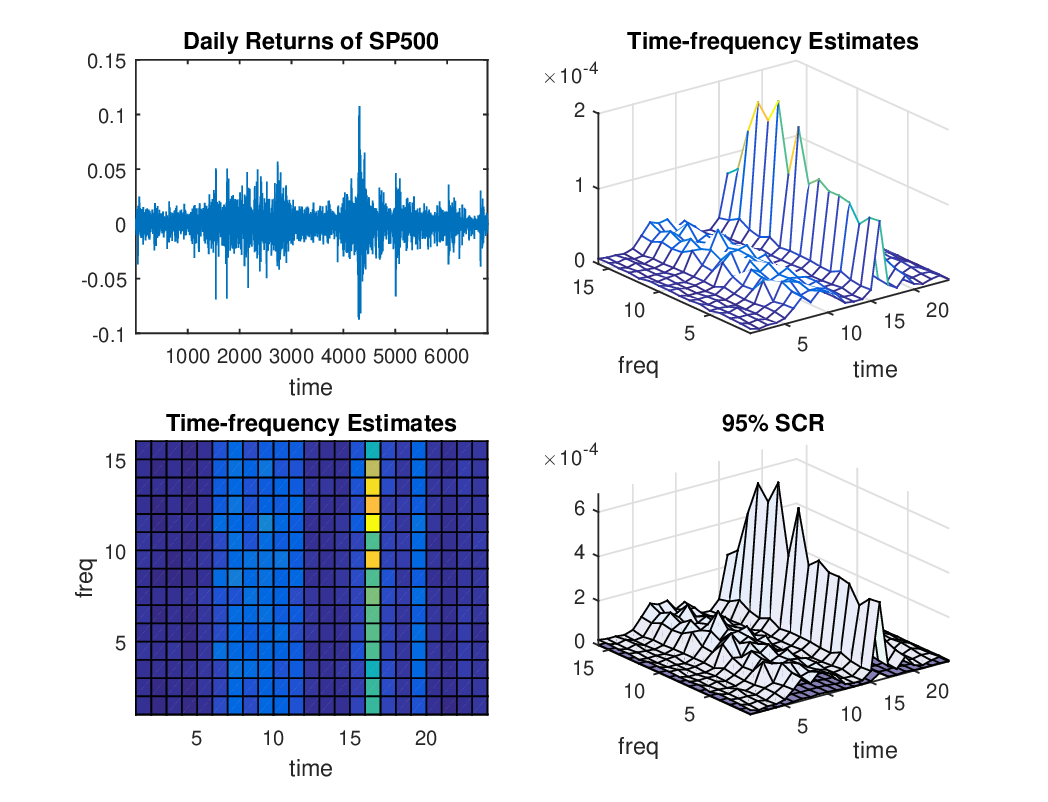}\label{fig_SP500}
	\end{minipage}\hfill
	\begin{minipage}{0.5\textwidth}
		\centering
			\caption{Analysis of Absolute SP500 Returns  }\includegraphics[width=\textwidth]{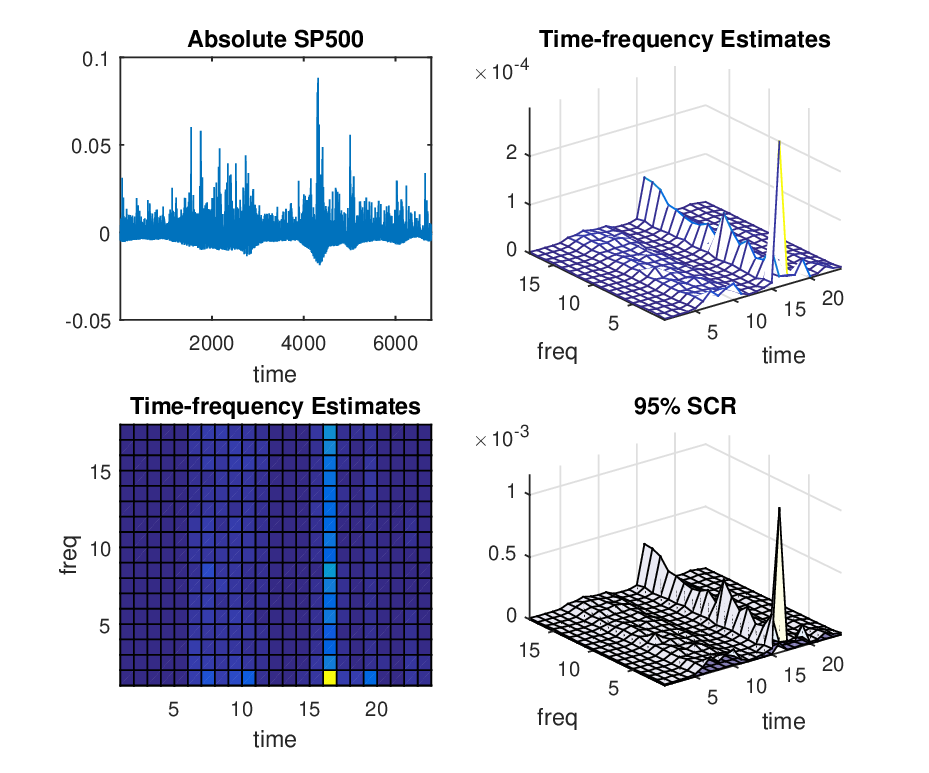}\label{fig_SP500abs}
	\end{minipage}
\end{figure}

\begin{figure}
	\centering
	\begin{minipage}{0.5\textwidth}
		\centering
			\caption{Absolute SP500 Return: Selected Times}\includegraphics[width=\textwidth]{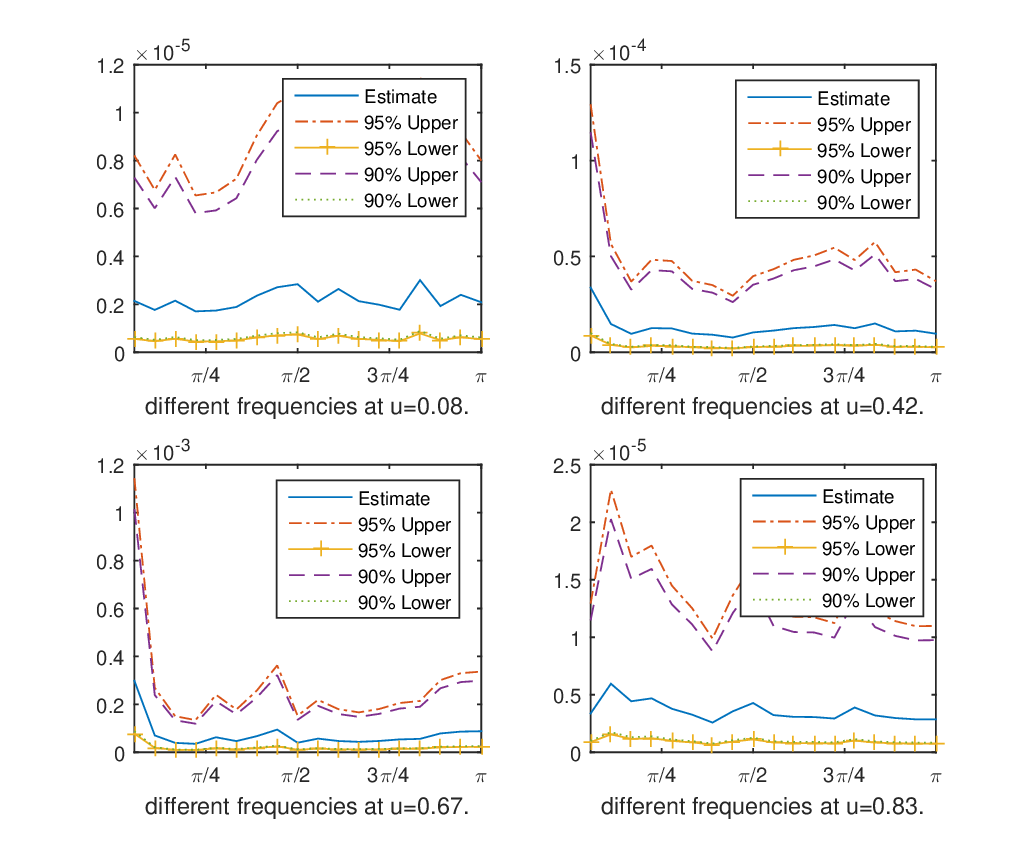}\label{fig_SP500abs_time}
	\end{minipage}\hfill
	\begin{minipage}{0.5\textwidth}
		\centering
				\caption{Absolute SP500 Return: Selected Frequencies}\includegraphics[width=\textwidth]{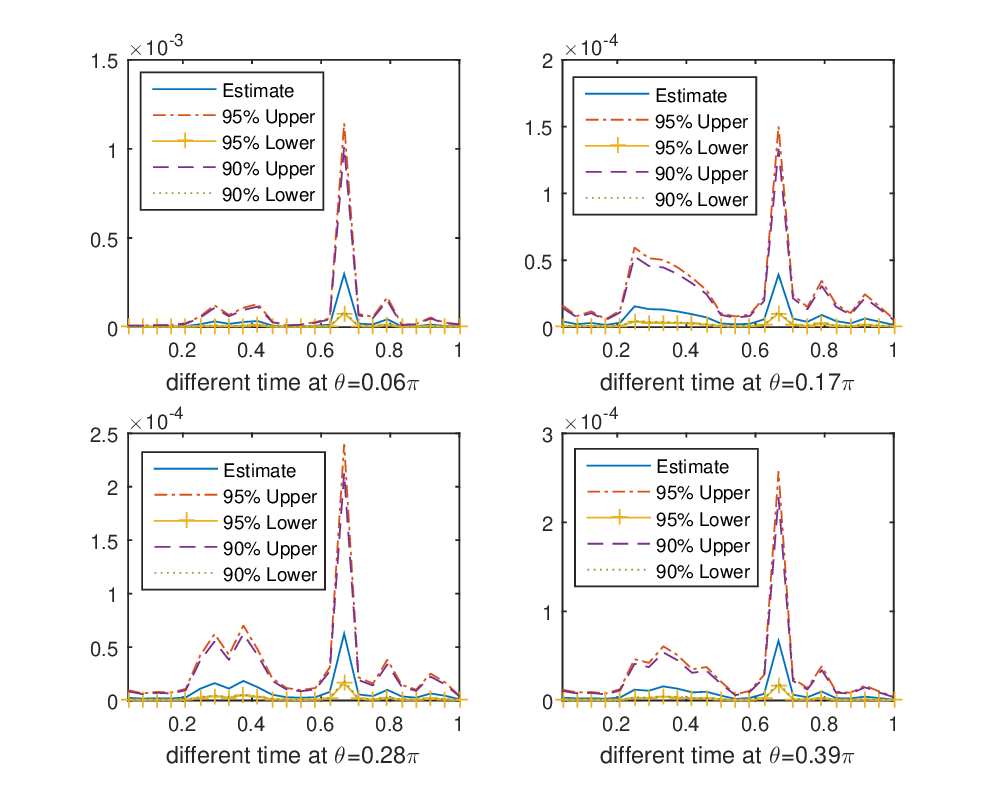}\label{fig_SP500abs_freq}
	\end{minipage}
\end{figure}

Finally, we fit time-varying non-stationary linear models for the absolute SP500 daily returns with mean removed by kernel smoothing. We first fit various time-varying AR or ARMA models 
\[
\sum_{i=0}^p a_i(t/N)X_{t-i}=\sum_{j=0}^q b_j(t/N)\epsilon_{t-j}
\]
to the absolute returns by minimizing the local Whittle likelihood \cite{Dahlhaus1997}. We then validate if the fitted spectral densities from the time-varying AR or ARMA models fall into the proposed SCR. The p-values for validating time-varying AR/ARMA models are shown in \cref{table_real_data_AR}. One can see that, the p-values for the tv-AR models are quite small, which implies that no tv-AR models up to order $5$ is appropriate for fitting absolute SP500 daily returns. For tv-ARMA models, the p-value for the tv-ARMA$(1,1)$ model equals $0.019$. This suggests that this tv-ARMA model is not appropriate for fitting the absolute SP500 daily returns either. In contrast, the corresponding p-value for validating the tv-ARMA$(2,1)$ model is $0.79$. This interesting observation suggests that the tv-ARMA$(2,1)$ model may be appropriate to fit the absolute returns. We further plot the spectral densities of the fitted time-varying AR$(1)$, AR$(4)$, AR$(5)$, ARMA$(1,1)$, ARMA$(2,1)$, and ARMA$(3,1)$ models in \cref{fig_parametric_SP500abs}.  From \cref{fig_parametric_SP500abs}, one can see that the fitted spectral densities by the tv-AR models are quite different from the STFT-based spectral density estimates. For tv-ARMA models, the spectral density estimates by the tv-ARMA$(1,1)$ model are not close to the STFT-based spectral density estimates either. Therefore, based on the proposed SCR, we conclude that the tv-ARMA$(2,1)$ model is an appropriate candidate for the analyzed data and can be used for short-term future volatility forecasting.

	\begin{table}
	\caption{p-values for Validating Time-varying ARMA Models to Absolute SP500}\label{table_real_data_AR}
	\begin{tabular}{c|c|c|c}
		\hline
		Model & p-value & Model & p-value\\
		\cline{1-4}
		\hline
		tv-AR$(1)$  &   $0.0066^{**}$ & tv-ARMA$(1,1)$ & $0.019^{*}$  \\	
		tv-AR$(2)$ &  $0.0015^{**}$ & tv-ARMA$(2,1)$ & $0.79$ \\	
		tv-AR$(3)$ &  $0.0015^{**}$ & tv-ARMA$(3,1)$ & $0.77$ \\	
		tv-AR$(4)$ &  $0.0012^{**}$ & tv-ARMA$(4,1)$ & $0.78$ \\	
		tv-AR$(5)$ &  $0.0012^{**}$ & tv-ARMA$(5,1)$  & $0.84$ \\ 
		\hline
	\end{tabular}

	Signif. codes: $(***)<0.001\le (**)<0.01\le (*)<0.05\le (+)<0.1$. 
\end{table}


\begin{figure}
	\caption{Fitting Absolute SP500 Daily Returns to Time-varying ARMA Models}\includegraphics[width=0.68\textwidth]{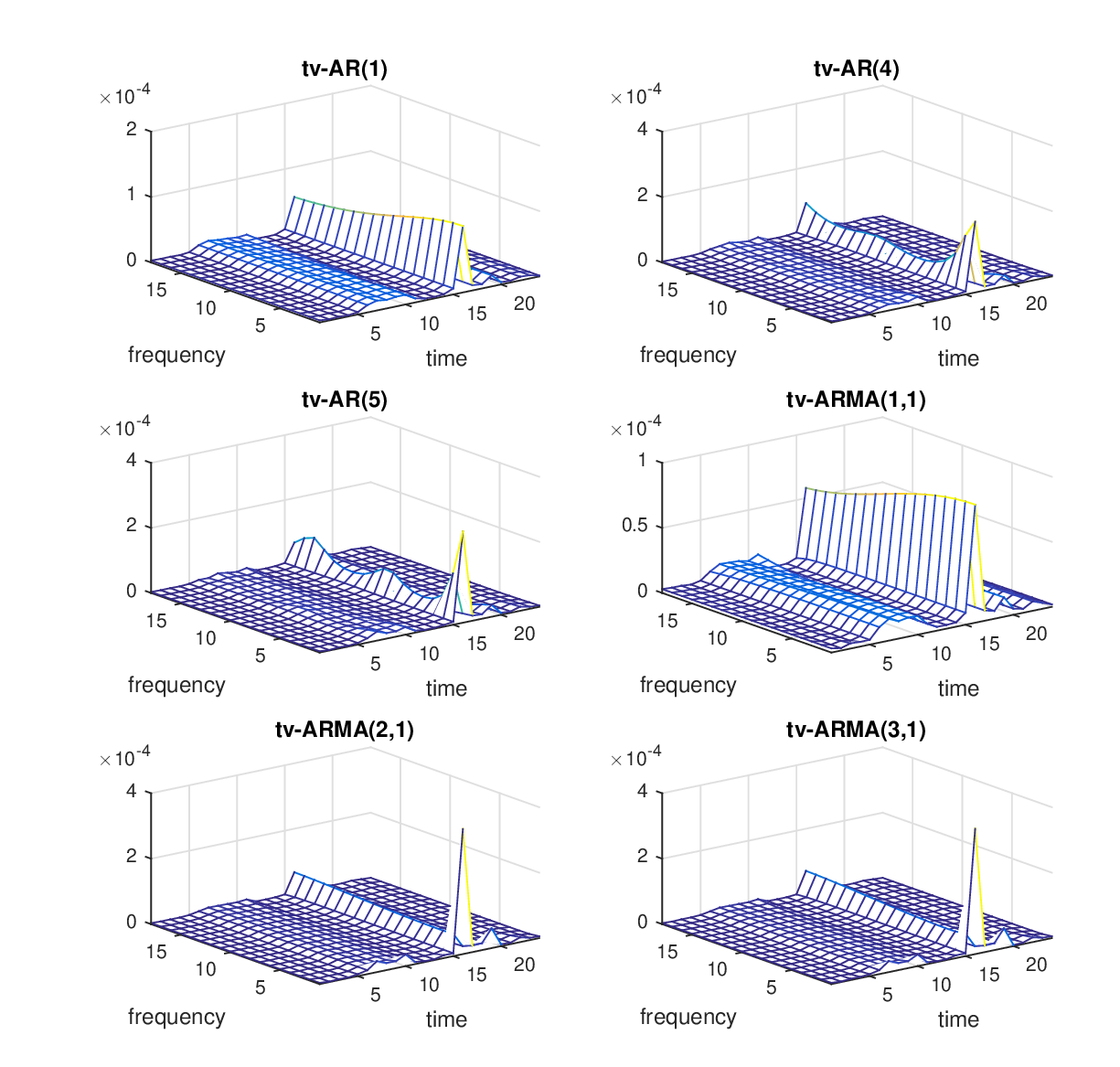}\label{fig_parametric_SP500abs}
\end{figure}


\end{example}

\section{Proofs of Main Results}\label{section_proof}

\subsection{Proof of \cref{thm_fourier}}\label{proof_thm_fourier}
	We prove \cref{thm_fourier} in two steps. In the first step, we show in \cref{subsubsection_proof_fourier_step1} that \cref{thm_fourier} is true for $q=1$. In this case, we let $\Omega_{p}=\{c\in \mathbb{R}^p: |c|=1\}$, $Z_{u,J}=(Z_{u,j_1}^{(n)},\dots,Z_{u,j_p}^{(n)})^T$ for $J=(j_1,\dots,j_p)$ satisfies $1\le j_1,\dots,j_p\le 2m$ (recall that  $m=\floor{(n-1)/2}$).
	We prove for any fixed $p\in\mathbb{N}$, as $n\to \infty$, we have that
	\[\label{proof_fourier_step1}
	\sup_{u}\sup_{J}\sup_{c\in\Omega_p}\sup_x |P(c^TZ_{u,J}\le x)-\Phi(x)|=o(1).
	\]
	In the second step of the proof, we show in \cref{subsubsection_proof_fourier_step2}  that for fixed $q\in\mathbb{N}$, for any given $0<u_1<\dots< u_q<1$, we have
	$\{(c^{(i)})^TZ_{u_i,J}, i=1,\dots,q\}$ are asymptotically independent uniformly over $\{c^{(i)}\in \mathbb{R}^p: |c^{(i)}|=1\}$ for $i=1,\dots,q$. Finally,  \cref{thm_fourier} is proved by combining the two parts.
	\subsubsection{Proof of \cref{proof_fourier_step1}}\label{subsubsection_proof_fourier_step1}
	
	We denote $2\pi j/n$ by $\theta_j$ in this proof. 
	With out loss of generality, we restrict $J=\{j_1,\dots,j_p\}\in \{1,\dots,m\}$. Let $c=(c_1,\dots,c_p)$, define $\mu_{u,k}:=\sum_{\ell=1}^p\frac{c_{\ell}\cos(k\theta_{j_{\ell}})}{\sqrt{\pi f(u,\theta_{j_l})}}$. Then
	$\mu_{u,k}\le  \sum_{\ell=1}^p\frac{|c_{\ell}|}{\sqrt{\pi f_*}}\le \frac{p}{\sqrt{\pi f_*}} =:\mu_*,\quad \forall c\in \Omega_p,\quad \forall J$.
	Furthermore, defining
	\[\label{proof_fourier_def_Tun}
	T_{u,n}:=\sum_{k=1}^N \mu_{u,k}\tau\left(\frac{k-\lfloor uN\rfloor}{n}\right)X_k,\quad
	{\tilde{T}_{u,n}:=\sum_{k=1}^N \mu_{u,k}\tau\left(\frac{k-\lfloor uN\rfloor}{n}\right)\tilde{X}_k^{[\ell]}},
	\] 
	and $\eta:=\left(\frac{\|T_{u,n}-\tilde{T}_{u,n}\|}{\sqrt{n}}\right)^{1/2}$,
	we have the following key lemmas.
	\begin{lemma}\label{lemma1a} Under the assumptions of \cref{thm_fourier}, we have
		\[
		\lim_{n\to \infty}\sup_J\sup_c\sup_u \left|\frac{\|T_{u,n}\|^2}{n}-1\right|^2=0.
		\]
	\end{lemma}
	\begin{proof}
		See \cref{proof_lemma1a}.
	\end{proof}
	\begin{lemma}\label{lemma1b}
		Under the assumptions of \cref{thm_fourier}, we have
		\[
		\lim_{\ell\to \infty}\sup_J\sup_c\sup_u\frac{\|T_{u,n}-\tilde{T}_{u,n}\|}{\sqrt{n}}=0.
		\]
	\end{lemma}
	\begin{proof}
		See \cref{proof_lemma1b}.
	\end{proof}

	\begin{lemma}\label{lemma1c} Under the assumptions of \cref{thm_fourier}, we have
		\[
		\begin{split}
		\sup_x& \left|\Pr\left(\frac{T_{u,n}}{\sqrt{n}}\le x\right)-\Phi\left(\frac{x}{\|T_{u,n}\|/\sqrt{n}}\right)\right|=\bigO\left(\Pr\left(\left|\frac{T_{u,n}-\tilde{T}_{u,n}}{\sqrt{n}}\right|\ge \eta\right)+\delta_n+\eta^2\right),
		\end{split}
		\]
		where $\delta_n\to 0$ as $n\to \infty$ uniformly over $J$, $c$ and $u$.
	\end{lemma}
	\begin{proof}
		See \cref{proof_lemma1c}.
	\end{proof}
	Using the above results, we can then prove \cref{proof_fourier_step1} as follows. First,
	{by Chebyshev inequality and $\eta=\left(\frac{\|T_{u,n}-\tilde{T}_{u,n}\|}{\sqrt{n}}\right)^{1/2}$, we have}
	\[
	\Pr\left(\left|\frac{T_{u,n}-\tilde{T}_{u,n}}{\sqrt{n}}\right|\ge \eta\right)\le\frac{\EE(T_{u,n}-\tilde{T}_{u,n})^2/n}{\eta^2}=\eta^2.
	\]
	Next, according to \cref{lemma1a}, uniformly over $J$, $c$ and $u$, for any fixed $\ell$, as $n\to \infty$, we have that
	\[
	\begin{split}
	\sup_x \left|\Pr\left(\frac{T_{u,n}}{\sqrt{n}}\le x\right)-\Phi\left(\frac{x}{\|T_{u,n}\|/n}\right)\right|&\to \sup_x \left|\Pr\left(\frac{T_{u,n}}{\sqrt{n}}\le x\right)-\Phi\left(x\right)\right|.
	\end{split}
	\]
	By \cref{lemma1c}, we have that
	\[
	\sup_x \left|\Pr\left(\frac{T_{u,n}}{\sqrt{n}}\le x\right)-\Phi\left(\frac{x}{\|T_{u,n}\|/n}\right)\right|=\bigO(2\eta^2+\delta_n).
	\]
	Note that $\delta_n\to 0$ as $n\to \infty$. Also, {by \cref{lemma1b}, uniformly over $J,c,u,n$, we have $\eta\to 0$ as $\ell\to \infty$}. Finally, letting $n\to\infty$ then $\ell\to \infty$, we have that
	$\sup_x \left|\Pr\left(\frac{T_{u,n}}{\sqrt{n}}\le x\right)-\Phi\left(x\right)\right|\to 0$,
	uniformly over $J$, $c$, and $u$.
	\subsubsection{Proof of asymptotically independence of 	$\{(c^{(i)})^TZ_{u_i,J}, i=1,\dots,q\}$ }\label{subsubsection_proof_fourier_step2}
	We can write $T_{u_i,n}$ and $\tilde{T}_{u_i,n}$ defined in \cref{proof_fourier_def_Tun} as $T_{u_i,n,c^{(i)}}$ and $\tilde{T}_{u_i,n,c^{(i)}}$. Then by \cref{lemma1b}, it suffices to show that $\{\tilde{T}_{u_i,n,c^{(i)}}, i=1,\dots,q\}$ are asymptotically independent uniformly over $\{c^{(i)}\in \mathbb{R}^p: |c^{(i)}|=1\}$. Note that in the definition of $\tilde{T}_{u_i,n,c^{(i)}}$, $\tilde{X}_k$ is $\ell$-dependent, therefore, $\tilde{T}_{u_1,n,c^{(i)}}$ and $\tilde{T}_{u_2,n,c^{(i)}}$ with $u_2>u_1$ are independent if $\lfloor(u_2-u_1)N\rfloor>\ell+2n$. Since $0<u_1<\dots<u_q<1$ are fixed, $\min_{i\neq j} |u_i-u_j|>0$ is bounded away from zero. Therefore, $\{\tilde{T}_{u_i,n,c^{(i)}}, i=1,\dots,q\}$ are independent if $\ell<\lfloor(\min_{i\neq j} |u_i-u_j|)N\rfloor-2n$. Choosing $\ell=o(n)$ and $n=o(N)$, we have $\{\tilde{T}_{u_i,n,c^{(i)}}, i=1,\dots,q\}$ are asymptotically independent.
\subsection{Proof of \cref{thm_consistency}}\label{proof_thm_consistency}
	Throughout the proof, we use  $\|\cdot\|$ to denote $\|\cdot\|_2$ for simplicity. We define $X_{u,i,n}:=\tau\left(\frac{i-\lfloor n/2\rfloor}{n}\right)X_{\floor{uN}+i-\floor{n/2}}$. For simplicity we will omit the index $n$ and use $X_{u,i}$ for $X_{u,i,n}$. Define
	$Y_{u,i}:=Y_{u,i}(\theta)=\frac{1}{2\pi} \sum_{k=-B_n}^{B_n} X_{u,i}X_{u,i+k} a(k/B_n)\cos(k\theta)$,
	$g_{n}(u,\theta):=\sum_{i=1}^n Y_{u,i}(\theta)$, and {$ h_{n}(u,\theta):=\frac{1}{\sqrt{nB_n}}g_n(u,\theta)-\sqrt{n/B_n}\hat{f}_n(u,\theta)$},
	we have that
	\[
	\begin{split}
	&\sqrt{n/B_n}\{\hat{f}_{n}(u,\theta)-\EE(\hat{f}_{n}(u,\theta))\}=\frac{g_n(u,\theta)-\EE(g_n(u,\theta))}{\sqrt{nB_n}}-h_n(u,\theta)+\EE(h_n(u,\theta)).
	\end{split}
	\]
	Next, denote {$\tilde{X}_k^{[\ell]}$} as the $\ell$-dependent conditional expectation of $X_k$, {$\tilde{X}_{u,i}^{[\ell]}$} as the $\ell$-dependent conditional expectation of $X_{u,i}$, and $\tilde{Y}_{u,i}$ as the correspondence of sum using {$\tilde{X}_{u,i}^{[\ell]}$} instead of $X_{u,i}$, and $\tilde{g}_n$ as the correspondence of $g_n$ using $\tilde{Y}_{u,i}$ instead of $Y_{u,i}$. Note that under $\GMC(2)$ and $\sup_i\EE|X_i|^{4+\delta}<\infty$, we know $\GMC(4)$ holds. Then we have the following results.
	
	\begin{lemma}\label{lemma3c}
		Under the assumptions of \cref{thm_consistency}, $\GMC(4)$ holds with $0<\rho<1$, then
		{\[\label{tmp1}
		\sup_{\theta}\sup_u\|h_n(u,\theta)\|=(nB_n)^{-1/2} \bigO(B_n),
		\]
		\[
		\sup_{\theta}\sup_u\sup_i \|Y_{u,i}-\tilde{Y}_{u,i}\|=\bigO(B_n\rho^{\ell/4}),
		\]
		\[\label{tmp2}
		\sup_{\theta}\sup_u\|g_n(u,\theta)-\tilde{g}_n(u,\theta)\|=o(1).
		\]
		}
	\end{lemma}	
	\begin{proof}
		See \cref{proof_lemma3c}.
	\end{proof}
		
	Next, we apply the block method to $\{\tilde{Y}_{u,i}(\theta)\}$. Define
	\[\label{eq_block_method}
	U_{u,r}(\theta):=\sum_{i=(r-1)(p_n+q_n)+1}^{(r-1)(p_n+q_n)+p_n} \tilde{Y}_{u,i}(\theta),\quad 
	V_{u,r}(\theta):=\sum_{i=(r-1)(p_n+q_n)+p_n+1}^{r(p_n+q_n)} \tilde{Y}_{u,i}(\theta), \quad 1,\dots,k_n,
	\]
	where $k_n:=\lfloor n/(p_n+q_n)\rfloor$. 
	Let $p_n=q_n=\floor{n^{1-4\eta/\delta}(\log n)^{-8/\delta-4}}$ (i.e. same block length) and $\ell=\ell_n=\floor{-9\log n/\log \rho}$ (Note $B_n=o(p_n)$ since $\eta<\delta/(4+\delta)$). Then $U_{u,r}(\theta), r=1,\dots,k_n$ are independent (not identically distributed) block sums with block length $p_n$, and {$V_{u,r}(\theta), r=1,\dots,k_n-1$ are independent block sums with block length $q_n$.}
	Define $U_{u,r}'(\theta):=U_{u,r}(\theta)\Ind(|U_{u,r}(\theta)|\le d_n)$ where $d_n=\floor{\sqrt{nB_n}(\log n)^{-1/2}}$. Then we have the following results.
	\begin{lemma}\label{lemma3d}
		Under the assumptions of \cref{thm_consistency},  we have that
		\[\label{tmp3}
		\sup_u\EE(\max_{\theta}|V_{u,k_n}(\theta)|)=\bigO(\sqrt{p_n\ell_n}B_n),
		\]
		\[\label{tmp4}
		\sup_u \EE(\max_{\theta}|h_n(u,\theta)|)=o(1),
		\]
		\[\label{tmp5}
		\sup_u \max_r\max_{\theta} \var(U_{u,r}(\theta))=\bigO(p_nB_n).
		\]
		{Furthermore, we have that
		\[\label{tmp6}
		\var(U_{u,r}'(\theta))=\var(U_{u,r}(\theta))[1+o(1)],
		\]
		where the $o(1)$ term holds uniformly over $\theta$, $r$ and $u$.}
	\end{lemma}
	\begin{proof}
		See \cref{proof_lemma3d}.
	\end{proof}		
	\begin{lemma}\label{lemma3e}
		Let $U_{u,i}(\theta)$ be one of the block sums with block length $p_n$. Then we have that
		\[
		\sup_u\sup_i\sup_{\theta}\|U_{u,i}(\theta)\|_{2+\delta/2}=\bigO(\ell_n\sqrt{p_nB_n}).
		\]
	\end{lemma}
	\begin{proof}
	See \cref{proof_lemma3e}.
	\end{proof}
	Using the previous results \cref{tmp1,tmp2,tmp4}, we have that
	\[\label{tmp7}
	\begin{split}
	&\sup_{u}\max_{\theta}\sqrt{n/B_n} |\hat{f}_{n}(u,\theta)-\EE(\hat{f}_{n}(u,\theta))|\\
	&\le \frac{ \sup_u\max_{\theta} |\tilde{g}_n(u,\theta)-\EE(\tilde{g}_n(u,\theta))|+o(1)}{\sqrt{nB_n}}+\bigO_{\Pr}(\sqrt{B_n/n})+o_{\Pr}(1)\\
	&\le\frac{ \sup_u\max_{\theta} |\sum_{r=1}^{k_n} U_{u,r}(\theta)-\EE(\sum_{r=1}^{k_n} U_{u,r}(\theta))|}{\sqrt{nB_n}}\\
	&\qquad+\frac{\sup_u\max_{\theta}|\sum_{r=1}^{k_n-1} V_{u,r}(\theta)-\EE(\sum_{r=1}^{k_n-1} V_{u,r}(\theta))|}{\sqrt{nB_n}}\\
	&\qquad+ \frac{\sup_u\max_{\theta}|V_{u,k_n}(\theta))-\EE(V_{u,k_n}(\theta)))|}{\sqrt{nB_n}}+\bigO_{\Pr}(\sqrt{B_n/n})+o_{\Pr}(1).
	\end{split}
	\]
	{First, we can show that the third term of the right hand side of \cref{tmp7} is  $o_{\Pr}(\sqrt{\log n})$. This is because by \cref{tmp3}, it suffices to show $\frac{\sqrt{p_n\ell_n}B_n}{\sqrt{nB_n}}=o(\sqrt{\log n})$ and this can be easily verified using $p_n=n^{1-4\eta/\delta}(\log n)^{-8/\delta-4}$, $B_n=\bigO(n^{\eta})$ and $\delta\le 4$.}
	
	{Next, we show that the right hand side of the first two terms of \cref{tmp7} have a order of $\bigO_{\Pr}(\sqrt{\log n})$.} Let $H_{u,n}(\theta)=\sum_{r=1}^{k_n}[U_{u,r}(\theta)-\EE(U_{u,r}(\theta))]$ and $H_{u,n}'(\theta)=\sum_{r=1}^{k_n}[U_{u,r}'(\theta)-\EE(U_{u,r}'(\theta))]$. Let $\theta_j=\pi j/t_n, j=0,\dots, t_n$ where $t_n=\floor{B_n\log(B_n)}$. Then, since both $H_{u,n}$ and $H_{u,n}'$ have trigonometric polynomial forms, we can apply the following result from \cite[Corollary 2.1]{Woodroofe1967}.
	\begin{lemma}\label{lemma_woodroofe}
		Let $p(\lambda)=\sum_{v=-k}^k\alpha_v\exp(iv\lambda)$ be a trigonometric polynomial. Let $\lambda_i=\pi (i/rk), |i|\le rk$. Then
		$\max_{|\lambda|\le\pi} |p(\lambda)|\le \max_{|i|\le rk} |p(\lambda_i)/(1-3\pi r^{-1})|$.
	\end{lemma}
	\begin{proof}
		See \cite[Corollary 2.1]{Woodroofe1967}.
	\end{proof}		
	By setting $k=B_n$ and $r=\log(B_n)$ in \cref{lemma_woodroofe}, we get
	\[
	\max_{\theta} |H_{u,n}(\theta)|\le \frac{1}{1-3\pi/\log(B_n)}\max_{j\le t_n} |H_{u,n}(\theta_j)|.
	\]
	By \cref{tmp5,tmp6}, there exists a constant $C_1$ such that $$\sup_u\max_r\max_{\theta} \textrm{var}(U_{u,r}'(\theta))\le C_1p_nB_n.$$ Let $\alpha_n:=(C_1 n B_n \log n)^{1/2}$, by the union upper bound,
	\[
	\Pr(\max_{0\le j\le t_n} |H_{u,n}'(\theta_j)|\ge 4\alpha_n)\le  \sum_{j=0}^{t_n} \Pr(|H_{u,n}'(\theta_j)|\ge 4\alpha_n).
	\]
	Then we apply Bernstein's inequality (see  \cref{lemma_bernstein}) 
	to $\Pr(|H_{u,n}'(\theta_j)|\ge 4\alpha_n)$. This leads to, uniformly over $u$ and $\theta_j$,
	\[
	\begin{split}
	\Pr(|H_{u,n}'(\theta_j)|\ge 4\alpha_n)&\le \exp\left(\frac{-16\alpha_n^2}{2k_n C_1p_n B_n+\frac{8}{3} d_n\alpha_n}\right)\le C\exp\left(-\frac{nB_n \log n}{nB_n}\right).
	\end{split}
	\]
	Therefore, uniformly over $u$, we have that
	$\Pr(\max_{0\le j\le t_n} |H_{u,n}'(\theta_j)|\ge 4\alpha_n)= \bigO(t_n)\bigO(1/n)=o(1)$.
	Let $U_{u,n}^*(\theta)=U_{u,n}(\theta)-U_{u,n}'(\theta)$ and $H_{u,n}^*(\theta)=H_{u,n}(\theta)-H_{u,n}'(\theta)$. By the union upper bound and Chebyshev's inequality
	\[
	\begin{split}
	\Pr(\max_{0\le j\le t_n} |H_{u,n}^*(\theta_j)|\ge 4\alpha_n)&\le \sum_{j=0}^{t_n} \Pr(|H_{u,n}^*(\theta_j)|\ge 4\alpha_n)\le\sum_{j=0}^{t_n}\frac{\sum_{i=1}^{k_n}\textrm{var}(U_{u,i}^*(\theta_j))}{16\alpha_n^2}.
	\end{split}
	\]

	Using  \cref{lemma3e}, 
	$\sup_u\max_i\sup_{\theta}\|U_{u,i}(\theta)\|_{2+\delta/2 }=\bigO(\ell_n\sqrt{p_nB_n})$,
	and
	\[\label{tmp_inequalities}
	\textrm{var}(U_{u,i}^*\Ind_{|U_{u,i}^*|>d_n})=d_n^2\textrm{var}\left(\frac{U_{u,i}^*}{d_n}\Ind_{|U_{u,i}^*|>d_n}\right)\le d_n^2 \EE\left[{\left(\frac{U_{u,i}^*}{d_n}\right)}^{2+\delta/2}\right],
	\] 
	we have that
	\[
	\begin{split}
	&\sum_{j=0}^{t_n}\frac{\sum_{i=1}^{k_n}\textrm{var}(U_{u,i}^*(\theta_j))}{16\alpha_n^2}=\bigO\left(\frac{t_nk_n(\sqrt{p_nB_n}\ell_n)^{2+\delta/2}}{\alpha_n^2d_n^{\delta/2}}\right)\\
	&=\bigO\left(\frac{(B_n\log B_n)(n/p_n)(\sqrt{p_nB_n}\log n)^{2+\delta/2}}{(nB_n\log n)(nB_n)^{\delta/4}(\log n)^{-\delta/4}}\right)\\
	&=\bigO\left(\frac{(p_nB_n)^{1+\delta/4}(\log n)^{2+\delta/2}}{p_n(nB_n)^{\delta/4}(\log n)^{-\delta/4}}\right)=\bigO(p_n^{\delta/4}(B_n/n)^{\delta/4}(\log n)^{2+\delta/2+\delta/4}).
	\end{split}
	\]
	Using $p_n=n^{1-4\eta/\delta}(\log n)^{-8/\delta-4}$
	we have
	$p_n^{\delta/4}=(n^{\delta/4-\eta})(\log n)^{-2-\delta}$.
	Therefore,
	\[
	\begin{split}
	\sum_{j=0}^{t_n}\frac{\sum_{i=1}^{k_n}\textrm{var}(U_{u,i}^*(\theta_j))}{16\alpha_n^2}&=\bigO\left(\frac{t_nk_n(\sqrt{p_nB_n}\ell_n)^{2+\delta/2}}{\alpha_n^2d_n^{\delta/2}}\right)=\bigO(n^{-\eta}B_n^{\delta/4}(\log n)^{-\delta/4}).
	\end{split}
	\]
	Finally, $B_n=\bigO(n^{\eta}),\delta\le 4$ implies $B_n^{\delta/4}=\bigO(n^\eta)$, so we have that
	$\sum_{j=0}^{t_n}\frac{\sum_{i=1}^{k_n}\textrm{var}(U_{u,i}^*(\theta_j))}{16\alpha_n^2}=o(1)$.
	Therefore, uniformly over $u$, we have  $\max_{\theta}|H_{u,n}'(\theta)|=\bigO_{\Pr} (\alpha_n)$ and $\max_{\theta}|H_{u,n}^*(\theta)|=\bigO_{\Pr} (\alpha_n)$. Then 
	$\max_{\theta} |H_{u,n}(\theta)|=\max_{\theta} |H_{u,n}'(\theta)+H_{u,n}^*(\theta)|=\bigO_{\Pr} (\alpha_n)=\bigO_{\Pr}(\sqrt{nB_n\log n})$.
	So \cref{tmp7} has the order of $\bigO_{\Pr}(\sqrt{\log n})$.
\subsection{Proof of \cref{thm_normality}}\label{proof_thm_normality}
	Throughout the proof, we use  $\|\cdot\|$ to denote $\|\cdot\|_2$ for simplicity. We define $Y_{u,i}$, $g_n$, $h_n$, {$\tilde{X}_k^{[\ell]}$}, $\tilde{Y}_{u,i}$, $\tilde{g}_n$ 
	the same as in \cref{proof_thm_consistency}. Therefore, \cref{lemma3c} holds.
	Next, we apply the block method to $\{\tilde{Y}_{u,i}(\theta)\}$. 
	Define
	\[
	U_{u,r}(\theta):=\sum_{i=(r-1)(p_n+q_n)+1}^{(r-1)(p_n+q_n)+p_n} \tilde{Y}_{u,i}(\theta),\quad 
	V_{u,r}(\theta):=\sum_{i=(r-1)(p_n+q_n)+p_n+1}^{r(p_n+q_n)} \tilde{Y}_{u,i}(\theta), \quad 1,\dots,k_n,
	\]
	where $k_n:=\lfloor n/(p_n+q_n)\rfloor$.
	Let $\psi_n=n/(\log n)^{2+8/\delta}$, $p_n=\floor{\psi_n^{2/3}B_n^{1/3}}$, and $q_n=\floor{\psi_n^{1/3}B_n^{2/3}}$. Then we have $p_n,q_n\to\infty$ and $q_n=o(p_n)$. Since $\ell_n=\bigO(\log n)$, we have $2B_n+\ell_n=o(q_n)$ and $k_n=\floor{n/(p_n+q_n)}\to \infty$.	Note that $U_{u,r}(\theta), r=1,\dots,k_n$ are independent (not identically distributed) block sums with block length $p_n$, and $V_{u,r}(\theta), r=1,\dots,k_n$ are independent block sums with block length $q_n$. Now the proof of \cref{lemma3d} still follows.
	
	Defining $a_n/b_n\to 1$ by $a_n\sim b_n$, we have the following result.
	\begin{lemma}\label{lemma3b}
		Let the sequence $s_n\in \mathbb{N}$ satisfy $s_n\le n$, $s_n=o(n)$
		and $B_n=o(s_n)$. 
		Under $\textrm{GMC}(4)$ we have that
		\[
		\left\|\sum_{i=-s_n/2}^{s_n/2} \{Y_{u,i}(\theta)-\EE(Y_{u,i}(\theta))\}\right\|^2\sim s_nB_n\sigma^2_{u}(\theta),
		\]
		where
		$\sigma^2_{u}(\theta)=[1+\eta(2\theta)]f^2(u,\theta)\int_{-1}^{1}a^2(t)\dee t$
		and $\eta(\theta)=1$ if $\theta=2k\pi$ for some integer $k$ and $\eta(\theta)=0$ otherwise.
	\end{lemma}
	\begin{proof}
		See \cref{proof_lemma3b}.
	\end{proof}
		
	According to \cref{lemma3b,lemma3c}, for each block $U_{u,r}, r=1,\dots,k_n$, we have that
	\[
	\begin{split}
	\|U_{u,r}-\EE(U_{u,r})\|&=\left\|\sum_{j\in\mathcal{L}_r}\{\tilde{Y}_{u,j}-\EE(\tilde{Y}_{u,j})\}\right\|\\
	&=\left\|\sum_{j\in\mathcal{L}_r}\{Y_{u,j}-\EE(Y_{u,j})\}\right\|+\bigO\left(\sum_{j\in\mathcal{L}_r}\|Y_{u,j}-\tilde{Y}_{u,j}\|\right)\\
	&\sim (p_nB_n\sigma^2_{u})^{1/2}+\bigO(p_nB_n\rho^{\ell_n/4})\sim (p_nB_n\sigma^2_{u})^{1/2},
	\end{split}
	\]
	where $\mathcal{L}_r=\{j\in\mathbb{N}: (r-1)(p_n+q_n)+1\le j\le r(p_n+q_n)-q_n\}$. Similarly, we can also show that
	$\|V_{u,r}-\EE(V_{u,r})\|\sim (q_nB_n\sigma^2_u)^{1/2}+\bigO(q_nB_n\rho^{\ell_n/4})$.
	Then, since $q_n=o(p_n)$, we have that 
	\[
	\textrm{var}\left(\sum_{r=1}^{k_n-1}V_{u,r} +V_{u,k_n}\right)=(k_n-1)\bigO(q_nB_n\sigma^2_u)+\bigO((p_n+q_n)B_n)=o(nB_n)
	\] 
	which implies that 
	$\frac{\sum_r (V_{u,r}-\EE(V_{u,r}))}{\sqrt{nB_n}}\Rightarrow 0$.
	Also, by \cref{tmp1}, we have that
	$\textrm{var}(h_n(u,\theta))=\bigO(B_n/n)=\bigO((\log n)^{-2-8/\delta})$,
	which implies that
	$h_n(u,\theta)-\EE(h_n(u,\theta))\Rightarrow 0$.
	Therefore, by 
	\[
	\begin{split}
	&\sqrt{n/B_n}\{\hat{f}_{n}(u,\theta)-\EE(\hat{f}_{n}(u,\theta))\}=\frac{g_n(u,\theta)-\EE(g_n(u,\theta))}{\sqrt{nB_n}}-h_n(u,\theta)+\EE(h_n(u,\theta)),
	\end{split}
	\]
	we only need to show that
	$\frac{\sum_r (U_{u,r}-\EE(U_{u,r}))}{\sqrt{nB_n}}\Rightarrow \mathcal{N}(0,\sigma^2_u)$.
	We can check the conditions of \cref{lemma_berry-esseen}
	(the Berry--Esseen lemma) as follows.
	\[
	\EE\left(\frac{U_{u,r}-\EE(U_{u,r})}{\sqrt{nB_n}}\right)=0,\quad
	\sum_r \frac{\|U_{u,r}-\EE(U_{u,r})\|^2}{nB_n}\sim k_n \frac{p_nB_n\sigma^2_u}{nB_n}\sim \sigma^2_u.
	\]
	{By \cref{lemma3e}, we know $\|U_{u,r}\|_{2+\delta/2}=\bigO(\ell_n\sqrt{p_nB_n})$, which implies
	\[
	\sum_r \frac{\|U_{u,r}-\EE(U_{u,r})\|^{2+\delta/2}_{2+\delta/2}}{(nB_n)^{1+\delta/4}}=\bigO\left(k_n\frac{(\ell_n\sqrt{p_nB_n})^{2+\delta/2}}{(nB_n)^{1+\delta/4}}\right)=\bigO(\ell_n k_n^{-\delta/4}).
	\]
	}	
	Note that {$k_n=\floor{n/(p_n+q_n)}\sim n\psi^{-2/3}B_n^{-1/3}\sim n^{1/3}(\log n)^{(4/3+16/3\delta)}B_n^{-1/3}$, which implies $k_n^{-1}=\bigO((\log n)^{-4/3-16/3\delta})$.}
	{Then $\ell_n k_n^{-\delta/4}=\bigO((\log n)(\log n)^{(-\delta/3-4/3)})=\bigO((\log n)^{(-\delta/3-1/3)})\to 0$.} 
	Therefore, the result holds by \cref{lemma_berry-esseen}.
\subsection{Proof of \cref{thm_max_dev}}\label{proof_thm_max_dev}
	{Define $D_n=C_nB_n$, $\theta_i=\frac{i\pi}{B_n}, i=0,\dots,B_n$, and  $\alpha_{n,k}=a(k/B_n)\cos(k\theta)$. We use the previous definitions of $X_{u,k}$ and the $\ell$-dependent {$\tilde{X}_{u,k}^{[\ell]}$} as in \cref{proof_thm_consistency}. } Let $g_n(u,\theta):=[2\pi n \hat{f}_n(u,\theta)-\sum_{k=1}^n X_{u,k}^2]-\EE[2\pi n \hat{f}_n(u,\theta)-\sum_{k=1}^n X_{u,k}^2]$, where $\ell=\lfloor n^{\gamma}\rfloor$ for fixed $\gamma>0$ which is close to zero. Note that
	\[
	\begin{split}
	\hat{f}_n(u,\theta)-\EE(\hat{f}_n(u,\theta))
	&=\frac{1}{2\pi n}\sum_{1\le k,k'\le n}\alpha_{n,k-k'} [X_{u,k} X_{u,k'}-\EE (X_{u,k} X_{u,k'})]\\ &=\frac{1}{2\pi n}\left(g_n(u,\theta)+\sum_{k=1}^n(X_{u,k}^2-\EE X_{u,k}^2)\right).
	\end{split}
	\]
	Therefore, we have
	$g_n(u,\theta)=\sum_{1\le k,k'\le n, k\neq k'}\alpha_{n,k-k'} [X_{u,k} X_{u,k'}-\EE (X_{u,k} X_{u,k'})]$.
	Then let $\tilde{g}_n(u,\theta)$ be the corresponding version of $g_n(u,\theta)$ using $\ell$-dependent {$\{\tilde{X}_{u,k}^{[\ell]}\}$} instead of $\{X_{u,k}\}$. Define {$X_{u,k}'=\tilde{X}_{u,k}^{[\ell]}\Ind_{\left|\tilde{X}_{u,k}^{[\ell]}\right|\le (nB_n)^{\alpha}}$} where $\alpha<\frac{1}{4}$.
	Next, let $\bar{X}_{u,k}:=X_{u,k}'-\EE X_{u,k}'$ and define 
	\[
	\begin{split}
	\bar{g}_n&=
	2\sum_{1\le s<k\le n}\alpha_{n,k-s} [\bar{X}_{u,k} \bar{X}_{u,s}-\EE (\bar{X}_{u,k} \bar{X}_{u,s})]\\
	&=2\sum_{k=2}^n\bar{X}_{u,k}\sum_{s=1}^{k-1}\alpha_{n,k-s}\bar{X}_{u,s}-2\EE \sum_{k=2}^n\bar{X}_{u,k}\sum_{s=1}^{k-1}\alpha_{n,k-s}\bar{X}_{u,s}.
	\end{split}
	\]
	In the following, we show that $g_n(u,\theta)$ can be approximated by $\tilde{g}_n(u,\theta)$.
	\begin{lemma}\label{lemma_max_dev_1} Under the assumptions of \cref{thm_max_dev}, we have $\max_{u\in\mathcal{U}}\max_{0\le i\le B_n}\EE |g_n(u,\theta_i)-\tilde{g}_n(u,\theta_i)|=o(n^{1+\gamma}\rho^{\floor{n^{\gamma}}})$ and
		$\max_{u\in\mathcal{U}}\max_{0\le i\le B_n} \frac{|g_n(u,\theta_i)-\tilde{g}_n(u,\theta_i)|}{\sqrt{nB_n}}=o_{\Pr}(1)$.
	\end{lemma}
	\begin{proof}
		See  \cref{proof_lemma_max_dev_1}.
	\end{proof}
	
	Next, we show that $\tilde{g}_n(u,\theta)$ can be approximated by $\bar{g}_n(u,\theta)$.
	\begin{lemma}\label{lemma_max_dev_2}
		Under the assumptions of \cref{thm_max_dev}, we have that
		\[
		\EE\left(\max_{u\in\mathcal{U}}\max_{\theta}\frac{|\tilde{g}_n(u,\theta)-\bar{g}_n(u,\theta)|}{\sqrt{nB_n}}\right)=o(1).
		\]
	\end{lemma}
	\begin{proof}
		See  \cref{proof_lemma_max_dev_2}.
	\end{proof}	

	According to \cref{lemma_max_dev_1} and \cref{lemma_max_dev_2}, together with $\max_{i}|\tilde{g}_n(u,\theta_i)-\bar{g}_n(u,\theta_i)|\le \max_{\theta}|\tilde{g}_n(u,\theta)-\bar{g}_n(u,\theta)|$, we have that
	$\max_{u\in\mathcal{U}}\max_{0\le i\le B_n}\frac{|g_n(u,\theta)-\tilde{g}_n(u,\theta)|^2}{nB_n}=o_{\Pr}(1)$
	and
	\[
	\begin{split}
	&\Pr\left(\max_{u\in\mathcal{U}}\max_{0\le i\le B_n}\frac{|\tilde{g}_n(u,\theta_i)-\bar{g}_n(u,\theta_i)|^2}{nB_n}\ge y\right)\le \frac{\EE\left(\max_{u\in\mathcal{U}}\max_{\theta}\frac{|\tilde{g}_n(u,\theta)-\bar{g}_n(u,\theta)|^2}{nB_n}\right)}{y}=o(1).
	\end{split}
	\]
	Since $\max_u\max_i |\EE \tilde{g}_n(u,\theta_i)-\EE \bar{g}_n(u,\theta_i)|\le \EE(\max_u\max_i |\tilde{g}_n(u,\theta_i)-\bar{g}_n(u,\theta_i)|)$, it suffices to show that, for $D_n=B_nC_n$, we have that
	$$
	\Pr\left[\max_{0\le i\le B_n, u\in\mathcal{U}} \frac{|\bar{g}_n(u,\theta_i)-\EE(\bar{g}_n(u,\theta_i))|^2}{4\pi^2n B_n f_n^2(u,\theta_i)\int_{-1}^1 a(t)\dee t}-2\log D_n +\log(\pi \log D_n)\le x\right]\to e^{-e^{-x/2}}.
	$$

	Let $p_n=\floor{B_n^{1+\beta}}$, $q_n=B_n+\ell$, $\ell=\floor{n^{\gamma}}$ and $k_n=\floor{n/(p_n+q_n)}$, where $\gamma$ is small enough and $\beta>0$ is sufficiently close to zero. Split the interval $[1,n]$ into alternating big and small blocks $H_j$ and $I_j$ by
	\[
	\begin{split}
	H_j&=[(j-1)(p_n+q_n)+1,jp_n+(j-1)q_n],\quad 1\le j\le k_n,\\
	I_j&=[jp_n+(j-1)q_n+1,j(p_n+q_n)],\quad 1\le j\le k_n,\\
	I_{k_n+1}&=[k_n(p_n+q_n)+1,n].
	\end{split}
	\] 
	Define $\bar{Y}_{u,k}:=\bar{X}_{u,k}\sum_{s=1}^{k-1}\alpha_{n,k-s}\bar{X}_{u,s}$. Then
	$\bar{g}_n=\sum_{k=1}^n(\bar{Y}_{u,k}-\EE \bar{Y}_{u,k})$.
	For $1\le j\le k_n+1$, let
	\[
	U_j(u,\theta):=\sum_{k\in H_j}(\bar{Y}_{u,k}-\EE\bar{Y}_{u,k}),\quad
	V_j(u,\theta):=\sum_{k\in I_j}(\bar{Y}_{u,k}-\EE\bar{Y}_{u,k}).
	\]
	Then 
	$\bar{g}_n=\sum_{j=1}^{k_n} U_j +\sum_{j=1}^{k_n+1} V_j$.
	Next, define a truncated and normalized version of $U_j$ as
	\[
	\bar{U}_j(u,\theta):=U_j(u,\theta)\Ind\left(\frac{|U_j(u,\theta)|}{\sqrt{nB_n}}\le\frac{1}{(\log B_n)^4}\right)-\EE U_j(u,\theta)\Ind\left(\frac{|U_j(u,\theta)|}{\sqrt{nB_n}}\le\frac{1}{(\log B_n)^4}\right).
	\]
	In the following, we show that $\bar{g}_n(u,\theta_i)-\EE(\bar{g}_n(u,\theta_i))$ can be approximated by $\sum_{j=1}^{k_n}\bar{U}_j(u,\theta_i)$.
	\begin{lemma}\label{lemma_max_dev_3}
		Under the assumptions of \cref{thm_max_dev}, we have that
		\[
		\max_{u\in\mathcal{U}}\max_{0\le i\le B_n} \frac{\left|\bar{g}_n(u,\theta_i)-\EE(\bar{g}_n(u,\theta_i))-\sum_{j=1}^{k_n}\bar{U}_j(u,\theta_i)\right|}{\sqrt{nB_n}}=o_{\Pr}(1).
		\]
	\end{lemma}	
	\begin{proof}
		See  \cref{proof_lemma_max_dev_3}.
	\end{proof}
	Furthermore, we show in the following that $\sum_{j=1}^{k_n}\bar{U}_j(u,\theta_i)$ can be ignored if $i\notin [(\log B_n)^2, B_n-(\log B_n)^2]$.
	\begin{lemma}\label{lemma_max_dev_4}
		Under the assumptions of \cref{thm_max_dev}, we have that
		\[
		\Pr\left(\max_{u\in\mathcal{U}}\max_{i\notin [(\log B_n)^2, B_n-(\log B_n)^2]}\frac{\left|\sum_{j=1}^{k_n}\bar{U}_j(u,\theta_i)\right|}{\sqrt{nB_n}}\ge x\sqrt{\log (B_nC_n)}\right)=o(1).
		\]
	\end{lemma}	
	\begin{proof}
		See  \cref{proof_lemma_max_dev_4}.
	\end{proof}
	Finally, we complete the proof of \cref{eq_max_dev_step1} by the following result.
	\begin{lemma}\label{lemma_max_dev_5}
		Under the assumptions of \cref{thm_max_dev}, we have that
		\[
		\begin{split}
		&\Pr\left[\max_{u\in\mathcal{U}}\max_{(\log B_n)^2\le i\le B_n-(\log B_n)^2} \frac{\left|\sum_{j=1}^{k_n}\bar{U}_j(u,\theta_i)\right|^2}{4\pi^2n B_n f_n^2(u,\theta_i)\int_{-1}^1 a(t)\dee t}\right.\\
		&\qquad\left.-2\log D_n+\log(\pi \log D_n)\le x\right]\to e^{-e^{-x/2}}.
		\end{split}
		\]
	\end{lemma}	
	\begin{proof}
		See  \cref{proof_lemma_max_dev_5}.
	\end{proof}


\section*{Acknowledgement}
{The authors are grateful to the anonymous referees for their many helpful comments and suggestions which significantly improved the quality of the paper. }

\AtNextBibliography{\footnotesize}
\printbibliography

\appendix
\section{Supplemental Material}\label{section_supplement}
\renewcommand{\refname}{Supplemental Material}



\begin{remark}\label{remark_GMC}
	Denote $X_{u,i}:=G(i/N,\mathcal{F}_{u,i})$ where $\mathcal{F}_{u,i}=(\dots,\epsilon_{\floor{uN}},\epsilon_{\floor{uN}+1},\dots,\epsilon_{\floor{uN}+i})$. Let $\epsilon_k'$ be an i.i.d.\ copy of $\epsilon_k$ and $X_{u,i}':=G(i/N,\mathcal{F}_{u,i}')$
	where $\mathcal{F}_{u,i}'=(\dots,\epsilon_0',\dots,\epsilon_{\floor{uN}}',\epsilon_{\floor{uN}+1},\dots,\epsilon_{\floor{uN}+i})$ is a coupled version of $\mathcal{F}_{u,i}$. Then under $\GMC(p)$, $p>0$, there exist $C>0$ and $0<\rho=\rho(p)<1$ that do not depend on $u$, such that for any $u$ and $i$, we have
	\[\label{new_GMC}
	\sup_u\EE(|X_{u,i}'-X_{u,i}|^{p})\le C \rho^i.
	\]
	This is because, when $\GMC(p)$ holds, we have $\sup_u \EE(|X_{u,i}'-X_{u,i}|^p)\le \sum_{k=i}^{\infty} \delta_p(k)\le\bigO(\sum_{k=i}^{\infty}\rho^k)=\bigO(\rho^i)$.
	
	Furthermore, it can be easily shown that if $\textrm{GMC}(2)$ holds, then $\sup_u |r(u,k)|=\bigO(\rho^k)$ for some $\rho\in(0,1)$. Also, if $\sup_i \|X_i\|_p<\infty$ and $\GMC(\alpha)$ holds with any given $\alpha>0$, then $X_i$ is $\GMC(\alpha)$ with any $\alpha\in(0,p)$. In particular, if $\GMC(\alpha)$ holds with some $\alpha\ge 2$, then we must have $\sup_u\sum_{k=-\infty}^{\infty} |r(u,k)|<\infty$ since $\sup_u |r(u,k)|=\bigO(\rho^k)=o(k^{-2})$. Also, if $\GMC(2)$ holds as well as $\sup_i \EE( |X_i|^{4+\delta})<\infty$ for some $\delta>0$, then $\GMC(4)$ holds.
\end{remark}

\begin{lemma}{  (Berry-Esseen)}\label{lemma_berry-esseen}
	If $\{X_i, i\ge 1\}$ are independent random variables with $\EE(X_i)=0$,  $s_n^2=\sum_{i=1}^n \EE(X_i^2)>0$, $\sum_{i=1}^n \EE|X_i|^{2+\delta}<\infty$, for some $\delta\in (0,1]$ and $S_n=\sum_{i=1}^n X_i$, there exists a universal constant $C_{\delta}$ such that
	\[
	\sup_{-\infty<x<\infty} \left|\Pr(S_n<x s_n)-\Phi(x)\right|\le C_{\delta} \left(\frac{\sum_{i=1}^n \EE|X_i|^{2+\delta}}{s_n^{2+\delta}}\right).
	\]
\end{lemma}
\begin{proof}
	See \cite[pp. 304]{Chow1988}.
\end{proof}

\subsection{Proof of \cref{lemma1a}}\label{proof_lemma1a}
	Define $d_{u,n}(h)=\frac{1}{n}\sum_{k=1+h}^n \mu_{u,k}\mu_{u,k-h}$ for $0\le h\le n-1$ and $d_{u,n}(h)=0$ if $h\ge n$. Since
	\[
	\sum_{k=1}^n \cos(k\theta_{j_{\ell}})\cos((k+h)\theta_{j_{\ell'}})=\frac{n}{2}\cos(h\theta_{j_{\ell}})\Ind_{\{j_{\ell}=j_{\ell'}\}},
	\]
	using 
	\[
	\begin{split}
	d_{u,n}(h)&=\frac{1}{n}\sum_{k=1+h}^{n+h}\mu_{u,k}\mu_{u,k-h}-\frac{1}{n}\sum_{k=n+1}^{n+h}\mu_{u,k}\mu_{u,k-h}\\
	&=\sum_{\ell=1}^p c_{\ell}^2\frac{\cos(h\theta_{j_{\ell}})}{2\pi f(u,\theta_{j_{\ell}})}-\frac{1}{n}\sum_{k=n+1}^{n+h}\mu_{u,k}\mu_{u,k-h},
	\end{split}
	\] 
	we get that uniformly over $J$, $c$ and $u$, there exists $K_0$ such that
	\[
	\left|d_{u,n}(h)-\sum_{\ell=1}^p c_{\ell}^2\frac{\cos(h\theta_{j_{\ell}})}{2\pi f(u,\theta_{j_{\ell}})}\right|\le K_0\min\left\{\frac{h}{n},1\right\}.
	\]
	Next, we can write $\|T_{u,n}\|^2/n$ as
	\[
	\begin{split}
	&\frac{1}{n}\EE\left(\sum_{k=1}^N\mu_{u,k}\tau\left(\frac{k-\lfloor uN\rfloor}{n}\right)X_k\right)^2\\
	&= d_{u,n}(0)r(u,0)\left[\frac{1}{n}\sum_k \tau\left(\frac{k-\lfloor uN\rfloor}{n}\right)^2\right]\\
	&+2\sum_{h=1}^{\infty}d_{u,n}(h)r(u,h)\left[\frac{1}{n}\sum_{k}\tau\left(\frac{k-\lfloor uN\rfloor}{n}\right)\tau\left(\frac{k+h-\lfloor uN\rfloor}{n}\right)\right]
	+o(1).
	\end{split}
	\]
	Furthermore,  defining
	\[
	f_n(u,\theta):=\frac{1}{2\pi}\sum_{h=0}^{\infty} r(u,h)\cos(h\theta) \left[\frac{1}{n}\sum_k\tau\left(\frac{k-\lfloor uN\rfloor}{n}\right)\tau\left(\frac{k+h-\lfloor uN\rfloor}{n}\right)\right],
	\]
	we have that
	\[
	\begin{split}
	&\sum_h\left\{ \left[\frac{1}{n}\sum_k\tau\left(\frac{k-\lfloor uN\rfloor}{n}\right)\tau\left(\frac{k+h-\lfloor uN\rfloor}{n}\right)\right]r(u,h)\sum_{\ell=1}^p c_{\ell}^2\frac{\cos(h\theta_{j_{\ell}})}{2\pi f(u,\theta_{j_{\ell}})}\right\}\\
	&=\sum_{\ell=1}^p\frac{c_{\ell}^2}{2\pi f(u,\theta_{j_{\ell}})}\sum_h\left\{r(u,h)\cos(h_{\theta_{j_{\ell}}})\left[\frac{1}{n}\sum_k\tau\left(\frac{k-\lfloor uN\rfloor}{n}\right)\tau\left(\frac{k+h-\lfloor uN\rfloor}{n}\right)\right]\right\}\\
	&=\sum_{\ell=1}^p c_{\ell}^2\frac{f_n(u,\theta_{j_l})}{f(u,\theta_{j_l})}.\\
	\end{split}
	\]
	By the assumptions that $\tau\in \mathcal{C}^1([-1/2,1/2])$, $\int \tau^2(x)\dee x=1$, together with $\sup_u|r(u,h)|=o(h^{-2})$, and $\sum_{h=1}^{\infty}|r(u,h)|<\infty$, we have 
	$f_n(u,\theta)=f(u,\theta)+o(1)$, uniformly over $u$ and $\theta$.
	This implies that 
	\[
	\sum_{\ell=1}^p c_{\ell}^2\frac{f_n(u,\theta_{j_l})}{f(u,\theta_{j_l})}=\sum_{\ell=1}^p c_{\ell}^2+o(1)=1+o(1).
	\]
	Therefore, uniformly over $J$ and $c$, we have that
	\[\label{tmp}
	\begin{split}
	&\left|\frac{\|T_{u,n}\|^2}{n}-1\right|-o(1)\\
	&\le 2\sum_{h=0}^{\infty}	\left|d_{u,n}(h)-\sum_{\ell=1}^p c_{\ell}^2\frac{\cos(h\theta_{j_{\ell}})}{2\pi f_n(u,\theta_{j_{\ell}})}\right|r(u,h)\left[\frac{1}{n}\sum_k\tau\left(\frac{k-\lfloor uN\rfloor}{n}\right)\tau\left(\frac{k+h-\lfloor uN\rfloor}{n}\right)\right]\\
	&\le 2\sum_{h=0}^{\infty}K_0\min\left\{\frac{h}{n},1\right\}r(u,h)\left[\frac{1}{n}\sum_k\tau\left(\frac{k-\lfloor uN\rfloor}{n}\right)\tau\left(\frac{k+h-\lfloor uN\rfloor}{n}\right)\right].
	\end{split}
	\]
	Finally, since $\sup_u\sum_h |r(u,h)|<\infty$, we have $\sup_u\sum_{h>n} |r(u,h)|\to 0$. Also, as $n\to \infty$, 
	\[
	\begin{split}
	\sup_u\sum_{h<n} (h/n) r(u,h)&\le\sup_u \sum_{h<\sqrt{n}}(h/n)r(u,h)+\sup_u\sum_{\sqrt{n}\le h<n}(h/n)r(u,h)\\
	&\le \sup_u \sum_{h<\sqrt{n}} r(u,h)/\sqrt{n}+\sup_u \sum_{h>\sqrt{n}} r(u,h)	\\
	&\to 0.
	\end{split}
	\]
	Therefore, $\left|\frac{\|T_{u,n}\|^2}{n}-1\right|\to 0$.

\subsection{Proof of \cref{lemma1b}}\label{proof_lemma1b}
	{Throughout the proof, we write $\tilde{X}_k^{[\ell]}$ as $\tilde{X}_k$ for short.}
	Note that $\int \tau(x)\tau(x+h)\dee x\le \frac{1}{2}\int \left[\tau(x)^2+\tau(x+h)^2\right]\dee x=1$. For simplicity of the proof, we can assume that there exists some finite $\tau_*$ such that
	\[
	\frac{1}{n}\sum_k\tau\left(\frac{k-\lfloor uN\rfloor}{n}\right)\tau\left(\frac{k+h-\lfloor uN\rfloor}{n}\right)\le \tau_*^2.
	\]
	Then we have that
	\[
	\begin{split}
	&\frac{\|T_{u,n}-\tilde{T}_{u,n}\|}{\sqrt{n}}=\left[\frac{1}{n}\sum_{j=-\infty}^{\floor{uN+n/2}}\|\mathcal{P}_j(T_{u,n}-\tilde{T}_{u,n})\|^2\right]^{1/2}\\
	&\le \mu_*\tau_*\left[\frac{1}{n}\sum_{k=1}^n\sum_{j=-\infty}^{\floor{uN+n/2}}\|\mathcal{P}_j(X_{\floor{uN}+k-\floor{n/2}}-\tilde{X}_{\floor{uN}+k-\floor{n/2}})\|^2\right]^{1/2}\\
	&\le\mu_*\tau_*\max_{k\in\{1,\dots,n\}}\sum_{j=-\infty}^{\floor{uN+n/2}}\|\mathcal{P}_j(X_{\floor{uN}+k-\floor{n/2}}-\tilde{X}_{\floor{uN}+k-\floor{n/2}})\|\\
	&\le \mu_*\tau_*\max_{k\in\{1,\dots,n\}}\sum_{j=-\infty}^{\floor{uN+n/2}}\min\left\{2\|\mathcal{P}_j(X_{\floor{uN}+k-\floor{n/2}})\|,\|X_{\floor{uN}+k-\floor{n/2}}-\tilde{X}_{\floor{uN}+k-\floor{n/2}}\|\right\}\\
	&\le \mu_*\tau_*\sup_{k}\sum_{j=-\infty}^{k+n}\min\{2\|\mathcal{P}_j(X_k)\|,\|X_k-\tilde{X}_k\|\}\to 0, \quad \textrm{as } \ell\to \infty.
	\end{split}
	\]
	Since the upper bound does not depend on $u$, the convergence holds uniformly over $u$.

\subsection{Proof of \cref{lemma1c}}\label{proof_lemma1c}
	In this proof, we omit subscript $u$ for simplicity {and write $\tilde{X}_k^{[\ell]}$ as $\tilde{X}_k$ for short}.
		Since $\sup_k\EE(X_k^2)<\infty$, we have that
		\[
		\lim_{t\to \infty}\sup_k\EE[X_k^2\Ind(|X_k|>t)]=0.
		\]
	By the property of conditional expectation, we have $\EE(\tilde{X}_k^2)<\EE(X_k^2)$.
	Therefore, defining
	\[
	g_n(r)=r^2\sup_k \EE[\tilde{X}_k^2\Ind(|\tilde{X}_k|\ge \sqrt{n}/r)],
	\]
	we can get $\lim_{n\to \infty} g_n(r)=0$ for all given $r>0$. Also $g_n$ is non-decreasing with $r$. Then there exists a sequence $\{r_n\}$ such that $r_n\upto \infty$ and $g_n(r_n)\to 0$. Note that $r_n$ does not depend on $u$.
	
	For simplicity, we use $\tilde{X}_{u,k}$ to denote $\tilde{X}_{\floor{uN}+k-\floor{n/2}}$. Let $Y_{u,k}=\tilde{X}_{u,k}\Ind(|\tilde{X}_{u,k}|\le\sqrt{n}/r_n)$ and $T_{u,n,Y}=\sum_{k=1}^n\mu_{u,k}Y_{u,k}$. Since $\EE[\tilde{X}_k^2\Ind(|\tilde{X}_k|\ge \sqrt{n}/r)]=o(1/r_n^2)$ by the definition of $r_n$, we have $\|Y_{u,k}-\tilde{X}_{u,k}\|=o(1/r_n)$. Now since $Y_{u,k}-\tilde{X}_{u,k}$ is $\ell$-dependent, we divide each of $\{Y_{u,k}\}$ and $\{\tilde{X}_{u,k}\}$ into $\ell$ sub-sequences that each sub-sequences has  $\floor{n/\ell}$ independent elements. Then by the triangle inequality we can get
	\[
	\|T_{u,n,Y}-\tilde{T}_{u,n}\|\le \sum_{a=1}^{\ell}\left\|\sum_{b=a,a+\ell,\dots}^n\mu_{u,b}(Y_{u,b}-\tilde{X}_{u,b})\right\|=o(\sqrt{n}/r_n).
	\]
	
	Next, divide the sequence of $\{Y_{u,k}\}$ into pieces of length $p_n+\ell$ where $p_n=\floor{r_n^{1/4}}$.
	\[
	U_{u,t}=\sum_{a\in B_t} \mu_{u,a}Y_{u,a}
	\]
	where $B_t=\{a\in \mathbb{N}: 1+(t-1)(p_n+\ell)\le a\le p_n+(t-1)(p_n+\ell)\}$. Note that for given $u$, $\{U_{u,t}\}$ are independent (but not identically distributed) for different $t$.
	
	Define $V_{u,t}=\sum_{t=1}^{t_n} U_{u,t}$, then the difference between $V_{u,t}$ and $T_{u,n,Y}$ is the sum of those dropped $\ell$ terms in each piece. Since $\ell$ is fixed and there are $t_n$ blocks, we have $\|T_{u,n,Y}-V_{u,t}\|=\bigO(\sqrt{t_n})$. 
	
	Furthermore, since
	\[
	(\sqrt{n}/r_n)^2 \Pr(|\tilde{X}_k|\ge \sqrt{n}/r))\ge \EE[\tilde{X}_k^2\Ind(|\tilde{X}_k|\ge \sqrt{n}/r)]=o(1/r_n^2)
	\]
	we have $P(|\tilde{X}_k|\ge \sqrt{n}/r))=o(1/n)$. Then, using
	\[
	\begin{split}
	[\EE(Y_k)]^2&=[\EE(\tilde{X}_k)-\EE(Y_k)]^2=[\EE\tilde{X}_k\Ind(|\tilde{X}_k|\ge \sqrt{n}/r)]^2\\
	&\le \EE(\tilde{X}_k^2\Ind(|\tilde{X}_k|\ge \sqrt{n}/r))\Pr(|\tilde{X}_k|\ge \sqrt{n}/r))=o(1/r_n^2)o(1/n)
	\end{split}
	\]
	we have $\EE(Y_k)=o(\frac{1}{\sqrt{n}r_n})$, which implies $|\EE(V_n)|=O(n)|\EE(Y_k)|=o(\sqrt{n}/r_n)$.
	
	Next, defining $W=(V_n-\EE(V_n))/\sqrt{n}$ and $\Delta=\tilde{T}_n/\sqrt{n}-W$, we get
	\[
	\begin{split}
	\sqrt{n}\|\Delta\|=\|\tilde{T}_n-V_n+\EE(V_n)\|&\le |\EE(V_n)|+\|V_n-\tilde{T}_n\|\\
	&\le|\EE(V_n)|+\|V_n-T_{n,Y}\|+\|T_{n,Y}-\tilde{T}_n\|\\
	&=o(\sqrt{n}/r_n)+\bigO(\sqrt{t_n}+\sqrt{n}/r_n)=\bigO(\sqrt{t_n}).
	\end{split}
	\]
	
	Next, we apply \cref{lemma_berry-esseen} to $\{U_t-\EE(U_t), t=1,\dots, t_n\}$. Recall that $V_n=\sum_{t=1}^{t_n} U_t$ and $W=(V_n-\EE(V_n))/\sqrt{n}$, then
	\[
	\begin{split}
	&\sup_x\left|\Pr(V_n-\EE(V_n)<x \|V_n-\EE(V_n)\|)-\Phi(x)\right|\\
	&=\sup_x\left|\Pr(W<x \|W\|)-\Phi(x)\right|\\
	&\le C\sum_{t=1}^{t_n}\EE|U_t-\EE(U_t)|^3\|V_n-\EE(V_n)\|^{-3}\\
	&\le C\sum_{t=1}^{t_n}\EE|U_t|^3\|V_n-\EE(V_n)\|^{-3}.
	\end{split}
	\]
	Next, we get upper bounds of $\EE|U_t|^3$ and $\|V_n-\EE(V_n)\|^{-3}$. First, by H\"older's inequality $\sum_{a\in B_t} |Y_a|\le (\sum_{a\in B_t} |Y_a|^3)^{1/3}(\sum_{a\in B_t} 1)^{2/3}$, we have that
	\[ 
	\EE|U_t|^3\le \mu_*^3\EE\left|\sum_{a\in B_t} Y_a\right|^3\le \mu_*^3 p_n^2\sum_{a\in B_t} \EE|Y_a|^3\le \mu_*^3 p_n^2\sum_{a\in B_t} \EE(\frac{\sqrt{n}}{r_n}|Y_a|^2)=\bigO\left(\mu_*^3p_n^3\frac{\sqrt{n}}{r_n}\right).
	\]
	
	For sequences $a_n$ and $b_n$, we define $a_n=\Theta(b_n)$ if both $a_n=\bigO(b_n)$ and $b_n=\bigO(a_n)$. Then, using the definition of $\Theta(\cdot)$, the variance of $\sum_{a\in B_t} \mu_a Y_a$ has the order of $\Theta(p_n)$ because $Y_a$ is $\ell$-dependent. Then the variance of $V_n$ has the order of $\Theta(t_n p_n)=\Theta(n)$. Thus, 
	$\|V_n-\EE(V_n)\|^{-3}$ has an order of $\Theta(n^{-3/2})$. Overall, we have that
	\[
	\sup_x\left|\Pr(W<x \|W\|)-\Phi(x)\right|\le \bigO(\mu_*^3p_n^3(\sqrt{n}/r_n))\Theta(n)=\bigO(p_n^{-2}).
	\]
	
	To complete the proof, we first replace $V_n=\sum_t\sum_{a\in B_t}\mu_a Y_a$ by $\tilde{T}_n=\sum_k \mu_k \tilde{X}_k$ then by $T_n=\sum_k X_k$. Since 
	\[
	\{W\le x-\delta, |\Delta|<\delta\}\subseteq\{W+\Delta\le x\}\subseteq\{W\le x+\delta\}\cup\{|\Delta|\ge\delta\},
	\]
	we have that
	\[
	\Pr(W\le x-\delta)-\Pr(|\Delta|\ge \delta)\le \Pr(W+\Delta\le x)\le \Pr(W\le x+\delta)+\Pr(|\Delta|\ge \delta).
	\]
	Furthermore, one can get
	\[
	\begin{split}
	&\sup_x\left|\Pr(W<x \|W\|)-\Phi(x)\right|\\
	&=\sup_x\left|\Pr(W<x)-\Phi(x/\|W\|)\right|\\
	&=\sup_x\left|\Pr(\tilde{T}_n/\sqrt{n}-\Delta<x)-\Phi(x/\|W\|)\right|.
	\end{split}
	\]
	Using
	\[
	\Pr(W<x-\delta)-\Pr(|\Delta|\ge \delta)\le \Pr(\tilde{T}_n/\sqrt{n}<x)\le \Pr(W<x+\delta)+\Pr(|\Delta|\ge \delta),
	\]
	we get
	\[
	\sup_x \left|\Pr(\tilde{T}_n/\sqrt{n}<x)-\Pr(W<x)\right|\le \Pr(|\Delta|\ge \delta)=\bigO(\|\Delta\|^2/\delta^2)=\bigO(p_n^{-1}/\delta^2).
	\]
	Also
	\[
	\begin{split}
	&\sup_x|\Phi(x/\|W\|)-\phi(x/\|W+\Delta\|)|\\
	&= \bigO(\|W+\Delta\|/\|W\|-1)=\bigO(\|\Delta\|)=\bigO(\sqrt{t_n/n})=\bigO(p_n^{-1/2}).
	\end{split}
	\]
	Letting $\delta=p_n^{-1/4}$
	we have that
	\[
	\sup_x \left|\Pr(\tilde{T}_n/\sqrt{n}<x)-\Phi(x/\|W+\Delta\|)\right|=\bigO(p_n^{-2})+\bigO(p_n^{-1/2})+\bigO(p_n^{-1/2}).
	\]
	Finally, use the above technique again with $\Delta_1=(T_n-\tilde{T}_n)/\sqrt{n}$ and $\delta=\|\Delta_1\|^{1/2}$, we get
	\[
	\sup_x\left|\Pr(T_n/\sqrt{n}<x)-\Phi(\sqrt{n}x/\|T_n\|)\right|=\bigO(\Pr(|\Delta_1|\ge \|\Delta_1\|^{1/2})+p_n^{-1/2}+\|\Delta_1\|).
	\]
\begin{lemma}\label{lemma_bernstein}(Bernstein's inequality)
	Let $X_1,\dots,X_n$ be independent zero-mean random variables. Suppose $|X_i|\le M$ a.s., for all $i$. Then for all positive $t$,
	\[
	\Pr\left(\sum_i X_i>t\right)\le \exp\left(\frac{-\frac{1}{2}t^2}{\sum\EE(X_i^2)+\frac{1}{3}Mt}\right).
	\]
\end{lemma}
		
\begin{definition}\label{def_cum}
	Let $(U_1,\dots,U_k)$ be a random vector. Then the joint cumulant is defined as
	\[
	\textrm{cum}(U_1,\dots,U_k)=\sum (-1)^p (p-1)! \EE\left(\prod_{j\in V_1} U_j\right)\dots \EE\left(\prod_{j\in V_p} U_j\right),
	\]
	where $V_1,\dots,V_p$ is a partition of the set $\{1,2,\dots,k\}$ and the sum is taken over all such partitions.
\end{definition}
\begin{lemma}\label{lemma3a}
	Assume $\textrm{GMC}(\alpha)$ with $\alpha=k$ for some $k\in\mathbb{N}$, and $\sup_t \EE(|X_t|^k)<\infty$ Then there exists a constant $C>0$ such that for all $u$ and $0\le m_1\le\dots\le m_{k-1}$,
	\[
	|\textrm{cum}(X_{u,0},X_{u,m_1},\dots,X_{u,m_{k-1}})|\le C\rho^{m_{k-1}/[k(k-1)]},
	\]
	where $X_{u,i}:=\tau\left(\frac{i-\floor n/2\rfloor}{n}\right)X_{\floor{uN}+i-\floor{n/2}}$.
\end{lemma}
\begin{proof} Since $\tau(\cdot)$ is bounded, we have $\sup_{u}\sup_i \EE(|X_{u,i}|^k)<\infty$. We extend  \cite[Proposition 2]{Wu2004} to the cases of locally stationary time series.
	
	Given $1\le l\le k-1$, by multi-linearity of joint cumulants, we replace $X_{u,m_i}$ by independent $X_{u,m_i}'$ for all $i\ge l$ as follows
	\[
	\begin{split}
	&J:=\textrm{cum}(X_{u,0},X_{u,m_1},\dots,X_{u,m_{k-1}})\\
	&\quad=\textrm{cum}(X_{u,0},X_{u,m_1},\dots,X_{u,m_{l-1}},X_{u,m_l}',\dots,X_{u,m_{k-1}}')\\
	&\qquad+\textrm{cum}(X_{u,0},X_{u,m_1},\dots,X_{u,m_{l-1}},X_{u,m_l}-X_{u,m_l}',\dots,X_{u,m_{k-1}})\\
	&\qquad\dots\\
	&\qquad+\textrm{cum}(X_{u,0},X_{u,m_1},\dots,X_{u,m_{l-1}},X_{u,m_l}',\dots,X_{u,m_{k-1}}-X_{u,m_{k-1}}')\\
	&\qquad 	=:B+\sum_{i=l}^{k-1}A_i.
	\end{split}
	\]
	Note that  $(X_{u,0},X_{u,m_1},\dots,X_{u,m_{l-1}})$ is independent with $(X_{u,m_l}',\dots,X_{u,m_{k-1}}')$. By \cite[pp.35]{Rosenblatt1985}, we have $B=0$.
	Suppose we have that
	\[\label{tmp10}
	|A_i|\le \frac{C}{k}\rho^{(m_{i}-m_{l-1})/k}\le \frac{C}{k}\rho^{(m_{l}-m_{l-1})/k}
	\] 
	for $l\le i\le k-1$ and some constant $C$ that does not depend on $l$. Then $|J|\le C\rho^{(m_l-m_{l-1})/k}$ for any $1\le l\le k-1$. Then we get
	\[
	|J|\le C \min_{l} \rho^{(m_l-m_{l-1})/k}=C \rho^{\max_l\frac{m_l-m_{l-1}}{k}}\le C \rho^{m_{k-1}/k(k-1)}.
	\]
	Next, we show \cref{tmp10}. In particular, we show the case $i=l$ and the other cases can be proven similarly. Note that $\EE(|X_{u,i}|^k)$ is uniformly bounded, by the definition of joint cumulants in \cref{def_cum}, we only need to show that for $V\subset \{0,\dots,k-1\}$ such that $l\notin V$, we have that
	\[
	\EE\left((X_{u,m_l}-X_{u,m_l}')\prod_{j\in V} X_{u,m_j}\right)\le C\rho^{(m_l-m_{l-1})/k}.
	\]
	Letting $|V|$ be the cardinality of the set $V$, then $|V|\le k-1$, and we have
	\[
	\begin{split}
	\left|\EE\left[\left(\prod_{j\in V} X_{u,m_j}\right)^{\frac{1+|V|}{|V|}}\right]\right|&\le \EE\left[\left(\frac{1}{|V|}\sum_{j\in V}|X_{u,m_j}|^{|V|}\right)^{\frac{1+|V|}{|V|}}\right]\\
	&\le \EE\left(\frac{1}{|V|}\sum_{j\in V}|X_{u,m_j}|^{1+|V|}\right)\\
	&\le\max_{j\in V}\EE\left(|X_{u,m_j}|^{1+|V|}\right)\le M.
	\end{split}
	\]
	By H\"older's inequality and Jensen's inequality
	\[
	\begin{split}
	&\left|\EE\left((X_{u,m_l}-X_{u,m_l}')\prod_{j\in V} X_{u,m_j}\right)\right|\\
	&\quad\le\left\|X_{u,m_l}-X_{u,m_l}'\right\|_{1+|V|}\left\|\prod_{j\in V} X_{u,m_j}\right\|_{\frac{1+|V|}{|V|}}\\
	&\quad\le\left\|X_{u,m_l}-X_{u,m_l}'\right\|_{k}M^{\frac{|V|}{1+|V|}}\le (C'\rho^{m_l-m_{l-1}})^{1/k} M'\le C\rho^{(m_l-m_{l-1})/k}.
	\end{split}
	\]
\end{proof}

\subsection{Proof of \cref{lemma3c}}\label{proof_lemma3c}
	{Throughout this proof, we write $\tilde{X}_k^{[\ell]}$ as $\tilde{X}_k$ and $\tilde{X}_{u,i}^{[\ell]}$ as $\tilde{X}_{u,i}$ for short.}
	First, letting $\alpha_k=a(k/B_n)\cos(k\theta)$, we have that
	\[
	\begin{split}
	h_n(u,\theta)&=\frac{1}{2\pi\sqrt{nB_n}}\left(\sum_{k=0}^{B_n}\sum_{j=n-k+1}^nX_{u,j}X_{u,j+k}\alpha_k+\sum_{k=-B_n}^{-1}\sum_{j=n+k+1}^nX_{u,j}X_{u,j+k}\alpha_k\right).
	\end{split}
	\]
	By the summability of cumulants of orders $2$ and $4$ \cite[page 185]{Rosenblatt1985}, one can get
	{\[
	\sup_{\theta}\sup_u \textrm{var}\left(\sum_{k=0}^{B_n}\sum_{j=n-k+1}^nX_{u,j}X_{u,j+k}\alpha_k\right)=\bigO(B_n^2).
	\]
	}
	Therefore, we have {$\sup_{\theta}	\sup_u\|h_n(u,\theta)\|=(nB_n)^{-1/2} \bigO(B_n)$}.
	
	Next, note that by the assumption of GMC($4$) defined in \cref{new_GMC}, 
	we have that
	\[
	\sup_u\sup_i \EE(|X_{u,i}-\tilde{X}_{u,i}|^{4})\le C \rho^{\ell_n}.
	\]
	Then we have
	{
	\[\label{tmp8}
	\begin{split}
	\sup_{\theta} \sup_u\sup_i\|Y_{u,i}-\tilde{Y}_{u,i}\|&\le\sup_u\sup_i\frac{1}{2\pi}\sum_{k=-B_n}^{B_n}
	\|X_{u,i}X_{u,i+k}-\tilde{X}_{u,i}\tilde{X}_{u,i+k}\|(\sup_{\theta} |\alpha_k|)\\
	&\le C\sup_u\sup_i\sum_{k=-B_n}^{B_n}
	\|(X_{u,i}-\tilde{X}_{u,i})X_{u,i+k}+\tilde{X}_{u,i}(X_{u,i+k}-\tilde{X}_{u,i+k})\|\\
	&=\bigO(B_n)\sup_u\sup_i\|X_{u,i}-\tilde{X}_{u,i}\|\\
	&=\bigO(B_n)\sup_u\sup_i(\EE(|X_{u,i}-\tilde{X}_{u,i}|)^4)^{1/4}\\
	&=\bigO(B_n\rho^{\ell_n/4}).
	\end{split}
	\]
	}
	Finally
	{ 
	\[
	\sup_{\theta}
	\sup_u\|g_n(u,\theta)-\tilde{g}_n(u,\theta)\|=\bigO\left(\sup_{\theta}\sup_u\sum_{i=1}^n\|Y_{u,i}-\tilde{Y}_{u,i}\|\right)=\bigO(nB_n\rho^{\ell_n/4})=o(1).
	\]
}

\subsection{Proof of \cref{lemma3d}}\label{proof_lemma3d}
	{We write $\tilde{X}_k^{[\ell]}$ as $\tilde{X}_k$ and $\tilde{X}_{u,i}^{[\ell]}$ as $\tilde{X}_{u,i}$ for short.}
	To show \cref{tmp3},
	since $\alpha_k$ is bounded, letting $z_n=k_n(p_n+q_n)+1-q_n$, we have that
	\[\label{tmp9}
	\sup_u\EE(\max_{\theta}|V_{u,k_n}(\theta)|)\le C \sum_{j=-B_n}^{B_n} \sup_u\EE|\sum_{i=z_n}^{n}\tilde{X}_{u,i}\tilde{X}_{u,i+j}|.
	\]
	Since $\tilde{X}_{u,i}\tilde{X}_{u,i+j}$ is $2\ell_n$-dependent, if $|j|<\ell_n$, we have
	\[
	\sup_u\left\|\sum_{i=z_n}^{n}\tilde{X}_{u,i}\tilde{X}_{u,i+j}\right\|=\bigO(2\ell_n \sqrt{(n-z_n)/2\ell_n})=\bigO(\sqrt{q_n\ell_n})=\bigO(\sqrt{p_n\ell_n}).
	\]
	If $|j|\le\ell_n$, since $\EE(\tilde{X}_{u,i}\tilde{X}_{u,i+j}\tilde{X}_{u,i'}\tilde{X}_{u,i'+j})=0$ if $|i-i'|>\ell_n$, we have that
	\[
	\begin{split}
	\sup_u\left\|\sum_{i=z_n}^{n}\tilde{X}_{u,i}\tilde{X}_{u,i+j}\right\|^2&=\sup_u\sum_{i,i'=z_n}^n \EE(\tilde{X}_{u,i}\tilde{X}_{u,i+j}\tilde{X}_{u,i'}\tilde{X}_{u,i'+j})\\
	&=\sup_u\sum_{i'=i-\ell_n}^{i+\ell_n}\sum_{i=z_n}^n \EE(\tilde{X}_{u,i}\tilde{X}_{u,i+j}\tilde{X}_{u,i'}\tilde{X}_{u,i'+j})\\
	&=\bigO(q_n\ell_n)=\bigO(p_n\ell_n),
	\end{split}
	\]
	where we have used the assumption $\sup_i\EE(|X_i|^{4+\delta})<M$. Therefore, we get \cref{tmp3}.
	
	To show \cref{tmp4}, we first define $\tilde{h}_n(u,\theta)$ by replacing $X_i$ by $\tilde{X}_i$. Then we can prove similarly to \cref{tmp8} that
	\[
	\sup_u \EE(\max_{\theta}|h_n(u,\theta)-\tilde{h}_n(u,\theta)|)=o(1).
	\]
	Therefore, it suffices to show $\sup_u \EE(\max_{\theta}|\tilde{h}_n(u,\theta)|)=o(1)$. Using similar technique to \cref{tmp9} we can show that
	\[
	\sup_u \EE(\max_{\theta}|\tilde{h}_n(u,\theta)|)=\frac{1}{\sqrt{nB_n}}\bigO(\sqrt{B_n\ell_n} B_n)=\bigO(\sqrt{\ell_n}B_n/\sqrt{n})=o(1),
	\]
	where we have used $\eta<\frac{1}{2}$ and $\sqrt{\ell_n}B_n/\sqrt{n}=\bigO((\log n)^{1/2} n^{\eta-1/2})=o(1)$.
	
	To show \cref{tmp5}, we note that GMC($4$) implies the absolute summability of cumulants up to the fourth order. Also, for zero-mean random variables $X,Y,Z,W$, the joint cumulants
	\[\label{tmp11}
	\textrm{cum}(X,Y,Z,W)=\EE(XYZW)-\EE(XY)\EE(ZW)-\EE(XZ)\EE(YW)-\EE(XW)\EE(YZ).
	\]
	Therefore, letting $\mathcal{L}_r$ be the set of the indices $i$'s such that $Y_{u,i}$ belongs to the block corresponding to $U_{u,r}$, we have that
	\[\label{tmp12}
	\begin{split}
	&\textrm{var}(U_{u,r}(\theta))=\left\|\sum_{i\in \mathcal{L}_r}\sum_{k=-B_n}^{B_n}[X_{u,i}X_{u,i+k}-\EE(X_{u,i}X_{u,i+k})]\alpha_k\right\|^2\\
	&=\sum_{i,j\in \mathcal{L}_r}\sum_{k,l=-B_n}^{B_n}\EE\{[X_{u,i}X_{u,i+k}-\EE(X_{u,i}X_{u,i+k})][X_{u,j}X_{u,j+l}-\EE(X_{u,j}X_{u,j+l})]\alpha_k\alpha_{l}\}\\
	&=\sum_{i,j\in \mathcal{L}_r}\sum_{k,l=-B_n}^{B_n} \textrm{cum}(X_{u,i},X_{u,i+k},X_{u,j},X_{u,j+l})\alpha_k\alpha_{l}\\
	&\qquad+\sum_{i,j\in \mathcal{L}_r}\sum_{k,l=-B_n}^{B_n}\EE(X_{u,i}X_{u,j})\EE(X_{u,i+k}X_{u,j+l})\alpha_k\alpha_{l}\\
	&\qquad+\sum_{i,j\in \mathcal{L}_r}\sum_{k,l=-B_n}^{B_n}\EE(X_{u,i}X_{u,j+l})\EE(X_{u,i+k}X_{u,j})\alpha_k\alpha_{l},
	\end{split}
	\]
	where the first term is finite since the fourth cumulants are summable. For the second term (the last term can also be shown similarly), we use the condition \cref{local_autocorrelation_condition}, so that 
	\[\label{tmp13}
	\EE(X_{u,i}X_{u,j})\EE(X_{u,i+k}X_{u,j+l})=[r(u,i-j)+o(1/n)][r(u,i-j+k-j)+o(1/n)].
	\]
	Then using $p_n=o(n)$, $B_n=o(p_n)$ and $\sup_u\sum_{k=-\infty}^{\infty} |r(u,k)|<\infty$, one can get
	\[
	\begin{split}
	&\sup_u \max_r\max_{\theta}\sum_{i,j\in \mathcal{L}_r}\sum_{k,l=-B_n}^{B_n}[r(u,i-j)+o(1/n)][r(u,i-j+k-j)+o(1/n)]\\
	&=\sup_u \max_r\max_{\theta}\sum_{i,j\in \mathcal{L}_r}r(u,i-j)\left[\sum_{k,l=-B_n}^{B_n}r(u,i-j+k-j)+ o(B_n/n)\right]\\
	&\le (2p_n+1)(2B_n+1)(\sup_u\sum_{k=-\infty}^{\infty}|r(u,k)|^2)+o(p_nB_n/n)=\bigO(p_nB_n).
	\end{split}
	\]
	
	To show \cref{tmp6}, we note that
	\[
	\textrm{var}(U_{u,r}')=\textrm{var}(U_{u,r})\left[1+\frac{2\EE(U_{u,r}')\EE(U_{u,r}-U_{u,r}')-2\textrm{var}(U_{u,r}-U_{u,r}')}{\textrm{var}(U_{u,r})}\right].
	\]
	From \cref{lemma3b}, we know that $\textrm{var}(U_{u,r}(\theta))\sim p_nB_n\sigma^2_{u}(\theta)$ and $\sigma^2_{u}(\theta)=[1+\eta(2\theta)]f^2(u,\theta)\int_{-1}^{1}a^2(t)\dee t\ge f_*^2\int_{-1}^{1}a^2(t)\dee t>0$. Thus, it suffices to show that
	\[
	\sup_u\sup_r\sup_{\theta} \EE(U_{u,r}')\EE(U_{u,r}-U_{u,r}')=o(p_nB_n),\quad\sup_u\sup_r\sup_{\theta} \textrm{var}(U_{u,r}-U_{u,r}')=o(p_nB_n).
	\]
	By \cref{lemma3e}, applying similar inequalities as \cref{tmp_inequalities}, we have that
	\[
	\begin{split}
	&\sup_u\sup_i\sup_{\theta}\textrm{var}(U_{u,r}-U_{u,r}')\\
	&\le\sup_u\sup_i\sup_{\theta}\frac{\|U_{u,r}\|^{2+\delta/2}_{2+\delta/2}}{d_n^{\delta/2}}\\
	&=\bigO((\ell_n\sqrt{p_nB_n})^{2+\delta/2}(\sqrt{nB_n}(\log n)^{-1/2})^{-\delta/2})\\
	&=\bigO(p_nB_n)\bigO((\log n)^{2+3\delta/4}(\sqrt{p_nB_n})^{\delta/2}(\sqrt{nB_n})^{-\delta/2})\\
	&=o(p_nB_n).
	\end{split}
	\]
	Finally, since $\EE(U_{u,r}')\le\EE(|U_{u,r}|)\le [\EE(|U_{u,r}|^{2+\delta/2})]^{\frac{1}{2+\delta/2}}$, using again similar inequalities as \cref{tmp_inequalities}, we have that
	\[
	\begin{split}
	&\sup_u\sup_r\sup_{\theta} \EE(U_{u,r}')\EE(U_{u,r}-U_{u,r}')\\
	&\le
	\sup_u\sup_r\sup_{\theta} \|U_{u,r}\|_{2+\delta/2}\frac{\|U_{u,r}\|_{2+\delta/2}^{2+\delta/2}}{d_n^{1+\delta/2}}\\
	&=\bigO(\|U_{u,r}\|_{2+\delta/2}/d_n)o(p_nB_n)\\
	&=\bigO(\sqrt{p_n/n}(\log n)^{3/2})o(p_nB_n)=o(p_nB_n).
	\end{split}
	\]

\subsection{Proof of \cref{lemma3e}}\label{proof_lemma3e}
{In this proof, we write $\tilde{X}_k^{[\ell]}$ as $\tilde{X}_k$ and $\tilde{X}_{u,i}^{[\ell]}$ as $\tilde{X}_{u,i}$ for short.} For simplicity, we first consider that $u$ and $i$ are fixed. Without loss of generality, we consider the first block sum ($i=1$) so
\[
U_{u,1}(\theta)=\sum_{j=1}^{p_n}\tilde{Y}_{u,j}(\theta).
\]
We will first show that
\[
\left\|\sum_{j=1}^{p_n}\sum_{k=-B_n}^{B_n} \tilde{X}_{u,}\tilde{X}_{u,j+k}\alpha_k\right\|_{2+\delta/2}=\bigO(\ell_n\sqrt{p_nB_n}),
\]
where $\alpha_k=a(k/B_n)\cos(k\theta)$. Then we conclude that $\bigO(\ell_n\sqrt{p_nB_n})$ is also uniformly over $u$ and $i$ since the assumption $\sup_u\sup_i \EE(|X_{u,i}|^{4+\delta})<M$. We first write by the triangle inequality
\[
\begin{split}
&\left\|\sum_{j=1}^{p_n}\sum_{k=-B_n}^{B_n} \tilde{X}_{u,j}\tilde{X}_{u,j+k}\alpha_k\right\|_{2+\delta/2}\\
&\le\left\|\sum_{j=1}^{p_n}\sum_{k=-B_n}^{0} \tilde{X}_{u,j}\tilde{X}_{u,j+k}\alpha_k\right\|_{2+\delta/2}+\left\|\sum_{j=1}^{p_n}\sum_{k=0}^{B_n} \tilde{X}_{u,j}\tilde{X}_{u,j+k}\alpha_k\right\|_{2+\delta/2}.
\end{split}
\]
Now consider two cases (i) $\ell_n=o(B_n)$, then
\[\label{eq_tmp_two_terms}
\begin{split}
\sum_{j=1}^{p_n}\sum_{k=-B_n}^{0} \tilde{X}_{u,j}\tilde{X}_{u,j+k}\alpha_k&=\sum_{j=1}^{p_n}\left(\tilde{X}_{u,j}\sum_{k=-B_n}^{-\ell_n}\tilde{X}_{u,j+k}\alpha_k\right)+\sum_{j=1}^{p_n}\sum_{k=1-\ell_n}^{0}\tilde{X}_{u,j}\tilde{X}_{u,j+k}\alpha_k,
\end{split}
\]
where the first term of the right hand side of \cref{eq_tmp_two_terms} satisfies
\[
\begin{split}
&\left\|\sum_{j=1}^{p_n}\left(\tilde{X}_{u,j}\sum_{k=-B_n}^{-\ell_n}\tilde{X}_{u,j+k}\alpha_k\right)\right\|_{2+\delta/2}\\
&\le\sum_{h=1}^{\ell_n}\left\|\sum_{j=1}^{\floor{(p_n-h)/\ell_n}}\tilde{X}_{u,h+(j-1)\ell_n}\sum_{k=-B_n}^{-\ell_n}\tilde{X}_{u,h+(j-1)\ell_n+k}\alpha_k\right\|_{2+\delta/2}.
\end{split}
\]
Continuing to divide the sum of $\sum_{k=-B_n}^{-\ell_n}\tilde{X}_{u,h+(j-1)\ell_n+k}\alpha_k$ into $\ell_n$ parts, then by $\sup_{u,i}\EE(|X_{u,i}|^{4+\delta})<M$, we have that
\[
\begin{split}
\left\|\sum_{j=1}^{p_n}\left(\tilde{X}_{u,j}\sum_{k=-B_n}^{-\ell_n}\tilde{X}_{u,j+k}\alpha_k\right)\right\|_{2+\delta/2}&=\bigO(\ell_n)\bigO(\sqrt{p_n/\ell_n})\bigO(\ell_n)\bigO(\sqrt{B_n/\ell_n})\\
&=\bigO(\ell_n\sqrt{p_nB_n}),
\end{split}
\]
which holds uniformly over $u$ and $i$. Similarly, for the second term of the right hand side of \cref{eq_tmp_two_terms}
\[
\begin{split}
&\left\|\sum_{j=1}^{p_n}\sum_{k=1-\ell_n}^{0}\tilde{X}_{u,j}\tilde{X}_{u,j+k}\alpha_k\right\|_{2+\delta/2}
\le \sum_{k=1-\ell_n}^{0}\left\|\sum_{j=1}^{p_n}\tilde{X}_{u,j}\tilde{X}_{u,j+k}\alpha_k\right\|_{2+\delta/2}\\
&=\sum_{k=1-\ell_n}^{0}\sum_{h=1}^{3\ell_n}\left\|\sum_{j=1}^{\floor{(p_n-h)/3\ell_n}}\tilde{X}_{u,h+3j\ell_n}\tilde{X}_{u,h+3j\ell_n+k}\alpha_k\right\|_{2+\delta/2}=\bigO(\ell_n^2\sqrt{p_n/\ell_n}).
\end{split}
\]
Note that the order  $\bigO(\ell_n^2\sqrt{p_n/\ell_n})$ also holds uniformly over $u$ and $i$. This is because $\|\tilde{X}_{u,h+3j\ell_n}\tilde{X}_{u,h+3j\ell_n+k}\|_{2+\delta/2}$ is uniformly bounded, which can be shown using Cauchy--Schwarz's inequality and $\sup_u\sup_i \EE(|X_{u,i}|^{4+\delta})<M$. Therefore, we have proven that, for case (i), we have
$\sup_u\sup_i\sup_{\theta}\|U_{u,i}(\theta)\|_{2+\delta/2}=\bigO(\ell_n\sqrt{p_nB_n})$.

For the second case (ii) $B_n=\bigO(\ell_n)$, we have that
\[
\begin{split}
&\left\|\sum_{j=1}^{p_n}\sum_{k=-B_n}^{0}\tilde{X}_{u,j}\tilde{X}_{u,j+k}\alpha_k\right\|_{2+\delta/2}
\le \sum_{k=-B_n}^{0}\left\|\sum_{j=1}^{p_n}\tilde{X}_{u,j}\tilde{X}_{u,j+k}\alpha_k\right\|_{2+\delta/2}\\
&=\sum_{k=-B_n}^{0}\sum_{h=1}^{3\ell_n}\left\|\sum_{j=1}^{\floor{(p_n-h)/3\ell_n}}\tilde{X}_{u,h+3j\ell_n}\tilde{X}_{u,h+3j\ell_n+k}\alpha_k\right\|_{2+\delta/2}\\
&=\bigO(B_n\ell_n\sqrt{p_n/\ell_n})=\bigO(\ell_n\sqrt{p_nB_n}),
\end{split}
\]
which is also uniform over $u$ and $i$.

\subsection{Proof of \cref{lemma3b}}\label{proof_lemma3b}		
Using the property of cumulants in \cref{tmp11}, similarly to \cref{tmp12,tmp13}, one can get that
\[\label{tmp14}
\begin{split}
&\left\|\sum_{i=-s_n/2}^{s_n/2} \{Y_{u,i}(\theta)-\EE(Y_{u,i}(\theta))\}\right\|^2\\
&=\sum_{i,j=-s_n}^{s_n}\sum_{k,l=-B_n}^{B_n} \textrm{cum}(X_{u,i},X_{u,i+k},X_{u,j},X_{u,j+l})\alpha_k\alpha_{l}\\
&\qquad+\sum_{i,j=-s_n}^{s_n}\sum_{k,l=-B_n}^{B_n}r(u,i-j)r(u,i+k-j-l)\alpha_k\alpha_{l}+o(s_nB_n/n)\\
&\qquad+\sum_{i,j=-s_n}^{s_n}\sum_{k,l=-B_n}^{B_n}r(u,i-j-l)r(u,i+k-j)\alpha_k\alpha_{l}+o(s_nB_n/n).
\end{split}
\]
By \cref{lemma3a}, we have that
\[
\begin{split}
&\sum_{m_1,m_2,m_3\in\mathbb{Z}}\textrm{cum}(X_{u,0},X_{u,m_1},X_{u,m_2},X_{u,m_3})<C\sum_{s=0}^{\infty}\rho^{s/[4(4-1)]}<\infty,
\end{split}
\]
which implies that the first term of the right hand side of \cref{tmp14} is finite.

Finally, according to \cite[Theorem 2, Eqs. (3.9)--(3.12)]{Rosenblatt1984}, one can show that
\[
\begin{split}
&\sum_{i,j=-s_n}^{s_n}\sum_{k,l=-B_n}^{B_n}r(u,i-j)r(u,i+k-j-l)\alpha_k\alpha_{l}\\
&+\sum_{i,j=-s_n}^{s_n}\sum_{k,l=-B_n}^{B_n}r(u,i-j-l)r(u,i+k-j)\alpha_k\alpha_{l}
\sim s_nB_n\sigma^2_{u}(\theta).
\end{split}
\]	

\begin{lemma}\label{lemma_max_dev_pre1}
	Let $\{X_k\}$ be $\ell$-dependent with $\EE X_k=0$ and $X_k\in \mathcal{L}^p$ with $p\ge 2$. Let $W_n=\sum_{k=1}^n X_k$. Then for any $Q>0$, there exists $C_1,C_2>0$ only depending on $Q$ such that
	\[
	\Pr(|W_n|\ge x)\le C_1\left(\frac{\ell}{x^2}\EE W_n^2\right)^Q+C_1\min\left[\frac{{\ell}^{p-1}}{x^p}\sum_{k=1}^n\|X_k\|_p^p,\sum_{k=1}^n \Pr\left(|X_k|\ge C_2\frac{x}{\ell}\right)\right].
	\]	
\end{lemma}	
\begin{proof}
	See \cite[Lemma 2]{Liu2010}.
\end{proof}

\begin{lemma}\label{lemma_max_dev_pre2} Let $\{X_t\}$ be $\ell$-dependent with $\EE X_t=0$, $|X_t|\le M$ a.s., $\ell\le n$, and $M\ge 1$. Let $S_{k,l}=\sum_{t=l+1}^{l+k}X_t\sum_{s=1}^{t-1} {\alpha}_{n,t-s} X_s$, where $l\ge 0$, $l+k\le n$ and assume that $\max_{1\le t\le n} |\alpha_{n,t}|\le K_0$, $\max_{1\le t\le n}\EE X_t^2\le K_0$, $\max_{1\le t\le n}\EE X_t^4\le K_0$ for some $K_0>0$. Then for any $x\ge 1$, $y\ge 1$, and $Q>0$,
	\[
	\begin{split}
	\Pr(|S_{k,l}-\EE S_{k,l}|\ge x)&\le 2 e^{-y/4}+C_1 n^3 M^2\left(x^{-2}y^2\ell^3(M^2+k)\sum_{s=1}^n \alpha_{n,s}^2\right)^Q\\
	&\quad +C_1n^3M^2\sum_{i=1}^n\Pr\left(|X_i|\ge \frac{C_2 x}{y\ell^2(M+k^{1/2})}\right),
	\end{split}
	\]
	where $C_1,C_2>0$ are constants depending only on $Q$ and $K_0$.
\end{lemma}
\begin{proof}
	See \cite[Proposition 3]{Liu2010}.
\end{proof}

\begin{lemma}\label{lemma_max_dev_pre3}
	Assume that $X_k\in \mathcal{L}^p$, with $p>1$, and $\EE X_k=0$. Let $C_p=18p^{3/2}(p-1)^{-1/2}$ and $p'=\min (2,p)$. Let $\alpha_1,\dots,\in \mathbb{C}$. Then under GMC, we have that
	\[
	\left\|\sum_{k=1}^n \alpha_k(X_k-{\tilde{X}_k^{[\ell]}})\right\|_p\le C_p \left(\sum_{k=1}^n|\alpha_k|^{p'}\right)^{1/p'} o(\rho^{\ell}),
	\]
	and 
	\[
	\|\sum_{k=1}^n \alpha_k X_k\|_p\le C\left(\sum_{k=1}^n|\alpha_k|^{p'}\right)^{1/p'},\quad 	\|\sum_{k=1}^n \alpha_k {\tilde{X}_k^{[\ell]}}\|_p\le C\left(\sum_{k=1}^n|\alpha_k|^{p'}\right)^{1/p'},
	\]
	for some constant $C$.
\end{lemma}
\begin{proof}
	This lemma follows from \cite[Lemma 1]{Liu2010} with $\Theta_{\ell+1,p}=o(\sum_{j=\ell+1}^{\infty}\rho^j)=o(\rho^{\ell})$.
\end{proof}	

\begin{lemma}\label{lemma_max_dev_pre4}
	Assume $\EE X_{u,k}=0$, $\sup_u\EE|X_{u,k}|^{2p}<\infty$, $p\ge 2$. Let 
	\[
	L_{n,u}=\sum_{1\le j\le j'\le n}\alpha_{j'-j} X_{u,j}X_{u,j'},\quad \tilde{L}_{n,u}=\sum_{1\le j\le j'\le n}\alpha_{j'-j} {\tilde{X}_{u,j}^{[\ell]}\tilde{X}_{u,j'}^{[\ell]}},
	\]
	where $\alpha_1,\dots,\in\mathbb{C}$. Then under GMC, we have that
	\[
	\frac{\sup_u\|L_{n,u}-\EE L_{n,u}-(\tilde{L}_{n,u}-\EE\tilde{L}_{n,u})\|_p}{n^{1/2}(\sum_{s=1}^{n-1}|\alpha_s|^2)^{1/2}}=o(\ell \rho^{\ell}).
	\]
\end{lemma}
\begin{proof}
	For fixed $u$, if $\EE|X_{u,k}|^{2p}<\infty$, the result follow from \cite[Proposition 1]{Liu2010} with $\Theta_{0,2p}=o(1)$ and $d_{\ell,2p}=\sum_{t=0}^{\infty}\min\{o(\rho^t),o(\rho^{\ell})\}=o(\ell\rho^{\ell})$. Since we have $\sup_u\EE|X_{u,k}|^{2p}<\infty$ the proof of \cite[Proposition 1]{Liu2010} also holds uniformly over $u$.
\end{proof}	

\begin{lemma}\label{lemma_max_dev_pre9}
	Assume that $\EE X_{u,k}=0$, $\sup_u\EE X_{u,k}^4<\infty$ and $\GMC(2)$.
	Let $\alpha_j=\beta_j\exp(ij\theta)$, where $i=\sqrt{-1}$, $\theta\in\Reals$, $\beta_j\in\Reals$, $1-n\le j\le -1$, $m\in\Naturals$ and $\tilde{L}_{n,u}=\sum_{1\le j<t\le n}\alpha_{j-t}{\tilde{X}_{u,j}^{[\ell]}\tilde{X}_{u,t}^{[\ell]}}$. Define
	\[
	D_{k}(u,\theta)=A_{u,k}-\EE(A_{u,k}\,|\,\mathcal{F}_{u,k-1}),\quad M_{n}(u,\theta)=\sum_{t=1}^{n}D_{t}(u,\theta)^*\sum_{j=1}^{t-1}\alpha_{j-t}D_{j}(u,\theta),
	\]
	where $(\cdot)^*$ denotes the complex conjugate, $A_{u,k}=\sum_{t=0}^{\infty}\EE({\tilde{X}_{u,t+k}^{[\ell]}}\,|\,\mathcal{F}_{u,k})\exp(ij\theta)$ where $\mathcal{F}_{u,k-1}:=\mathcal{F}_{\floor{uN-n/2}+k-1}$. Then
	\[
	\sup_u\frac{\|\tilde{L}_{n,u}-\EE\tilde{L}_{n,u}-M_{n}(u,\theta)\|}{m^{3/2}n^{1/2}\sup_k\|X_{u,k}\|^2_4}\le C V_m^{1/2}(\beta),
	\]
	where
	\[
	V_m(\beta)=\max_{1-n\le i\le -1}\beta_i^2+m\sum_{j=-1}^{-n-1}|\beta_j-\beta_{j-1}|^2.
	\]
\end{lemma}
\begin{proof}
	For fixed $u$, the result comes from \cite[Proposition 2]{Liu2010}. Since here we have assumed $\sup_u \EE X_{u,k}^4<\infty$, following the proof of  \cite[Proposition 2]{Liu2010}, the upper bound also holds uniformly over $u$.
\end{proof}	

\begin{lemma}\label{lemma_max_dev_pre5}
	Suppose that $\EE X_k=0$, $\sup_u \EE X_k^4<\infty$, and $\GMC(2)$ holds, then
	\begin{enumerate}
		\item We have that
		\[
		\left|\frac{\EE[(g_n(u_1,\theta_1)-\EE g_n(u_1,\theta_1))(g_n(u_2,\theta_2)-\EE g_n(u_2,\theta_2))]}{nB_n}\right|=\bigO(1/ (\log B_n)^2),
		\]
		uniformly on $(u_1,u_2,\theta_1,\theta_2)$ such that either $(u_1,u_2)\in \mathcal{U}^2$
		or $(\theta_1,\theta_2)\in \Theta^2$ where
		$\mathcal{U}^2=\{(u_1,u_2): \frac{n}{2N}\le u_1\le u_2\le 1-\frac{n}{2N}, |u_1-u_2|\ge \frac{n}{N}(1-1/(\log B_n)^2)  \}$ and $\Theta^2=\{(\theta_1,\theta_2): 0\le\theta_1<\theta_2\le \pi-B_n^{-1}(\log B_n)^2, |\theta_1-\theta_2|\ge B_n^{-1}(\log B_n)^2\}$.
		
		\item For $\alpha_n>0$ with $\limsup \alpha_n<1$, we have that
		\[
		\left|\frac{\EE[(g_n(u_1,\theta_1)-\EE g_n(u_1,\theta_1))(g_n(u_2,\theta_2)-\EE g_n(u_2,\theta_2))]}{4\pi^2nB_n f(u_1,\theta_1)f(u_2,\theta_2)\int_{t=-1}^{1} a^2(t)\dee t}\right|\le \alpha_n,
		\]		
		uniformly on $(u_1,u_2,\theta_1,\theta_2)$ such that either $(u_1,u_2)\in \mathcal{U}^2$
		or $(\theta_1,\theta_2)\in \bar{\Theta}^2$ where  $\mathcal{U}^2=\{(u_1,u_2): \frac{n}{2N}\le u_1\le u_2\le 1-\frac{n}{2N}, |u_1-u_2|\ge \frac{n}{N}(1-1/(\log B_n))  \}$ and  $\bar{\Theta}^2=\{(\theta_1,\theta_2): B_n^{-1}(\log B_n)^2\le\theta_1<\theta_2\le \pi-B_n^{-1}(\log B_n)^2, |\theta_1-\theta_2|\ge B_n^{-1}\}$.
		
		\item We have that
		\[
		\left|\frac{\EE[g_n(u,\theta)-\EE g_n(u,\theta)]^2}{4\pi^2n B_n f^2(u,\theta)\int_{t=-1}^{1}a^2(t)\dee t}-1\right|=\bigO(1/ (\log B_n)^2),
		\]
		uniformly on $\{(u,\theta): B_n^{-1}(\log B_n)^2\le \theta\le \pi-B_n^2(\log B_n)^2, \frac{n}{2N}< u< 1-\frac{n}{2N} \}$.		 
	\end{enumerate}
\end{lemma}
\begin{proof}
	{Throughout this proof, we write $\tilde{X}_k^{[\ell]}$ as $\tilde{X}_k$ and $\tilde{X}_{u,i}^{[\ell]}$ as $\tilde{X}_{u,i}$ for simplicity.}
	\begin{enumerate}
		\item By \cref{lemma_max_dev_pre4} we approximate $g_n-\EE g_n$ first by $\tilde{g}_n-\EE\tilde{g}_n$. Then by \cref{lemma_max_dev_pre9}, we approximate $\tilde{g}_n-\EE\tilde{g}_n$ by $M_n(u,\theta)$, where $M_n(u,\theta)=\sum_{t=1}^n D_{t}(u,\theta)^*\sum_{j=1}^{t-1}\alpha_{n,j-t}D_{j}(u,\theta)$. Then it is suffices to show that $|\EE [M_n(u_1,\theta_1)-M_n^*(u_1,\theta_1)][M_n(u_2,\theta_2)-M_n^*(u_2,\theta_2)]|\le C\frac{n B_n}{(\log B_n)^2}$ and $|\EE [M_n(u_1,\theta_1)+M_n^*(u_1,\theta_1)][M_n(u_2,\theta_2)+M_n^*(u_2,\theta_2)]|\le C\frac{n B_n}{(\log B_n)^2}$. We only prove the first inequality here, since the other inequality can be proved similarly. Define
		\[
		r_{n}(u_1,\theta_1,u_2,\theta_2):=|\EE [M_n(u_1,\theta_1)+M_n^*(u_1,\theta_1)][M_n(u_2,\theta_2)+M_n^*(u_2,\theta_2)]|.
		\]
		Since the martingale differences $\{D_{t}(u,\theta)\}$ are uncorrelated but not independent, we further define $N_n(u,\theta)=\sum_{t=1}^n D_{t}(u,\theta)^*\sum_{j=1}^{t-\ell-1}\alpha_{n,j-t}D_{j}(u,\theta)$, then $\|M_n(u,\theta)-N_n(u,\theta)\|=\bigO(\sqrt{n\ell})$ and $|r_{n}(u_1,\theta_1,u_2,\theta_2)|\le |\tilde{r}_{n}(u_1,\theta_1,u_2,\theta_2)|+\bigO(\sqrt{n\ell(n B_n)}+\sqrt{n\ell (B_n^2)})$, where
		\[
		\tilde{r}_{n}(u_1,\theta_1, u_2, \theta_2):=|\EE [N_n(u_1,\theta_1)+N_n^*(u_1,\theta_1)][N_n(u_2,\theta_2)+N_n^*(u_2,\theta_2)]|.
		\]
		Since $\ell=\floor{n^{\gamma}}$ where $\gamma$ is small enough, it suffices to show that $\tilde{r}_{n}(u_1,\theta_1,u_2,\theta_2)=\bigO(n B_n/(\log B_n)^2)$. Now we substitute $N_n(u,\theta)=\sum_{t=1}^n D_{t}(u,\theta)^*\sum_{j=1}^{t-\ell-1}\alpha_{n,j-t}D_{j}(u,\theta)$ to $\tilde{r}_{n}(u_1,\theta_1,u_2,\theta_2)$.
		
		If $\theta_1\neq \theta_2$ and $u_1=u_2$, we have that
		\[
		\sum_{t=1}^n\sum_{j=1}^{t-\ell-1}2\EE |D_{t}(u,\theta) D_{j}(u,\theta)|^2 a^2\left(\frac{t-j}{B_n}\right)[\cos((t-j)(\theta_1+\theta_2))+\cos((t-j)(\theta_1-\theta_2))].
		\]
		Now it suffices to show that $$\sum_{t=1}^n\sum_{j=1}^{t-\ell-1}a^2\left(\frac{t-j}{B_n}\right)\cos((t-j)(\theta_1\pm\theta_2))=\bigO(nB_n/(\log B_n)^2).$$
		Since $|\theta_1-\theta_2|\ge B_n^{-1}(\log B_n)^2$, using $1+2\sum_{k=1}^n \cos(k\theta)=\sin((n+1)\theta/2)/\sin(\theta/2)\le 1/\sin(\theta/2)$, $\sin(x)=\Theta(x)$ when $x\to 0$, and denoting $j=t-s$, we have that
		\[
		\sum_{t=1}^n \left|\sum_{j=1}^{B_n} a^2(j/B_n) \cos[j(\theta_1\pm\theta_2)] \right|
		\le C n/(B_n^{-1}(\log B_n)^2)=\bigO(n B_n/(\log B_n)^2).
		\]
		If $\theta_1= \theta_2$ but $u_1\neq u_2$, using \cref{tmp15} and $n-N|u_1-u_2|\le n/(\log B_n)^2$, we have that $$\tilde{r}_n(u_1,\theta,u_2,\theta)\le C \tilde{r}_{n-N|u_1-u_2|}(u,\theta,u,\theta)=\bigO((n-N|u_1-u_2|)B_n)
		=\bigO(nB_n/(\log B_n)^2).
		$$
		\item When $\theta_1\neq \theta_2$,  using \cite[Lemma 3.2(ii)]{Woodroofe1967} with the assumption on the continuity of $a(\cdot)$ in \cref{thm_max_dev}, we have that
		\[
		\lim_{n} \sup\, 2(n B_n)^{-1} \sum_{t=1}^n \sum_{j=1}^{t-\ell-1} a^2\left(\frac{t-j}{B_n}\right)\cos((t-j)(\theta_1-\theta_2))<\int a^2(t)\dee t.
		\]
		If $\theta_1=\theta_2$ and $u_1\neq u_2$ then
		\[
		\begin{split}
		&\lim_{n} \sup\, 2(n B_n)^{-1} \sum_{t=1}^{n-N|u_1-u_2|} \sum_{j=1}^{t-\ell-1} a^2\left(\frac{t-j}{B_n}\right)\\
		&\le\lim_{n} \sup\, 2(n B_n)^{-1}(n-N|u_1-u_2|)\sum_{j=-B_n}^{B_n}a^2\left(\frac{t-j}{B_n}\right)\\
		&\le \lim_{n} \sup\, 2(n B_n)^{-1} [nB_n/(\log B_n)^2]\int a^2(t)\dee t< \int a^2(t)\dee t.
		\end{split}
		\]
		\item Since $\|D_{t}(u,\theta)\|^2=\sum_{j=\ell}^{\ell} \EE (\tilde{X}_{u,t} \tilde{X}_{u,t+j}) \exp(ij\theta)$, we have that
		\[\label{tmp15}
		\begin{split}
		\tilde{r}_n(u,\theta,u,\theta)&=\bigO(nB_n/(\log B_n)^2)+\sum_{t=1}^n\|D_{t}(u,\theta)\|^2\,\sum_{s=-B_n}^{B_n} a^2(s/B_n)\\
		&=\bigO(n B_n/(\log B_n)^2) +4\pi^2f^2(u,\theta)nB_n\int a^2(t)\dee t.
		\end{split}
		\]
	\end{enumerate}
\end{proof}	

\begin{lemma}\label{lemma_max_dev_pre6}
	Let $X_1,\dots,X_m$ be independent mean zero $d$-dimensional random vectors such that $|X_i|\le M$. If the underlying probability space is rich enough, one can define independent normally distributed mean zero random vectors $V_1,\dots, V_m$ such that the covariance matrices of $V_i$ and $X_i$ are equal, for all $1\le i\le m$; furthermore
	\[
	\Pr\left(\left|\sum_{i=1}^m (X_i-V_i)\right|\ge \delta\right)\le c_1 \exp(-c_2\delta/M).
	\]
\end{lemma}	
\begin{proof}
	See \cite[Fact 2.2]{Einmahl1997}.
\end{proof}


\begin{lemma}\label{lemma_max_dev_pre8}
	If $X$ and $Y$ have a bi-variate normally distributed distribution with expectations $0$, unit variances, and correlation coefficient $r$, then
	\[
	\lim_{c\to\infty}\frac{\Pr(\{X>c\}\cap \{Y>c\})}{[2\pi (1-r)^{\frac{1}{2}} c^2]^{-1}\exp\left(-\frac{c^2}{1+r}\right)(1+r)^{\frac{3}{2}}}=1,
	\]
	uniformly for all $r$ such that $|r|\le \delta$, for all $0<\delta<1$.
\end{lemma}
\begin{proof}
	See \cite[Lemma 2]{Berman1962}.
\end{proof}


\subsection{Proof of \cref{lemma_max_dev_1}}\label{proof_lemma_max_dev_1}
	By Markov's inequality, we have that
	\[
	\begin{split}
	&\Pr\left(\max_{u\in\mathcal{U}}\max_{0\le i\le B_n} \frac{|g_n(u,\theta_i)-\tilde{g}_n(u,\theta_i)|}{\sqrt{nB_n}}
	\ge 1/\log D_n\right)\\
	&\le \sum_{u\in\mathcal{U}}\sum_{0\le i\le B_n} \Pr\left( \frac{|g_n(u,\theta_i)-\tilde{g}_n(u,\theta_i)|}{\sqrt{nB_n}}\ge 1/\log D_n\right)\\
	&\le C B_n C_n   \frac{\EE[ \frac{|g_n(u,\theta_i)-\tilde{g}_n(u,\theta_i)|}{\sqrt{nB_n}}]^{p/2} }{(1/\log D_n)^{p/2}}.
	\end{split}
	\]
	By \cref{lemma_max_dev_pre4},
	$\EE |g_n(u,\theta_i)-\tilde{g}_n(u,\theta_i)|=o(n^{1+\gamma}\rho^{\floor{n^{\gamma}}})$ uniformly on $u$ and $\theta_i$.
	Since $D_n=B_nC_n$ is polynomial of $n$, the GMC assumption guarantees
	\[
	\Pr\left(\max_{u\in\mathcal{U}}\max_{0\le i\le B_n} \frac{|g_n(u,\theta_i)-\tilde{g}_n(u,\theta_i)|}{\sqrt{nB_n}}
	\ge 1/\log D_n\right)=o(1).
	\]

\subsection{Proof of \cref{lemma_max_dev_2}}\label{proof_lemma_max_dev_2}	
\begin{lemma}\label{lemma_Aven85}
	Let $X_i, i=1,\dots,n$ be an arbitrary sequence of real-valued random variables with finite mean and variance. Then
	\[
	\EE(\max_{1\le i\le n} X_i)\le \max_{1\le i\le n} \EE X_i +\sqrt{\frac{n-1}{n}\sum_{i=1}^{n}\var(X_i)}.
	\]
\end{lemma}	
\begin{proof}
	See \cite[Theorem 2.1]{Aven1985}.
\end{proof}
	{In this proof, we write $\tilde{X}_k^{[\ell]}$ and $\tilde{X}_{u,i}^{[\ell]}$ as $\tilde{X}_k$ and $\tilde{X}_{u,i}$ for simplicity.}
	First of all, since $a(\cdot)$ has bounded support $[-1,1]$, we only need to consider the case that $|s-k|\le B_n$. Furthermore, let $\alpha_*$ be an upper bound of $\alpha_{n,i}$ uniformly over $u$. Then, we have that
	\[\label{eq_tmp2_two_terms}
	\begin{split}
	&\EE\left(\max_{u\in\mathcal{U}}\max_{\theta}|\tilde{g}_n(u,\theta)-\bar{g}_n(u,\theta)|\right)\\
	&\le \alpha_*\, \EE\left[\max_{u\in\mathcal{U}} \sum_{2\le k\le n, \max(1,k-B_n)\le s\le k-1}\left|\tilde{X}_{k,u}\tilde{X}_{s,u}-\EE(\tilde{X}_{k,u}\tilde{X}_{s,u})-\bar{X}_{k,u}\bar{X}_{s,u}+\EE(\bar{X}_{k,u}\bar{X}_{s,u})\right|\right]\\
	&\le 2\alpha_*\,\EE\left[\max_{u\in\mathcal{U}}\sum_{2\le k\le n, \max(1,k-B_n)\le s\le k-1}\left|\tilde{X}_{k,u}\tilde{X}_{s,u}-\bar{X}_{k,u}\bar{X}_{s,u}\right|\right]\\
	&= 2\alpha_*\,\EE\left[\max_{u\in\mathcal{U}}\sum_{2\le k\le n, \max(1,k-B_n)\le s\le k-1}\left|\tilde{X}_{k,u}\tilde{X}_{s,u}-\bar{X}_{k,u}\bar{X}_{s,u}-\tilde{X}_{k,u}\bar{X}_{s,u}+\tilde{X}_{k,u}\bar{X}_{s,u}\right|\right]\\
	&\le 2\alpha_*\, \EE\left[\max_{u\in\mathcal{U}}\left(\sum_{k=2}^n |\tilde{X}_{k,u}|\sum_{s=\max\{1,k-B_n\}}^{k-1} |\tilde{X}_{s,u}-\bar{X}_{s,u}|\right)\right]\\
	&\qquad+ 2\alpha_*\,\EE\left[\max_{u\in\mathcal{U}}\left(\sum_{k=2}^n|\tilde{X}_{k,u}-\bar{X}_{k,u}|\sum_{s=\max\{1,k-B_n\}}^{k-1}|\tilde{X}_{s,u}|\right)\right].
	\end{split}
	\]
	Next, we show that the first term of the right hand side of \cref{eq_tmp2_two_terms} satisfies
	\[
	\EE\left(\max_{u\in\mathcal{U}}\frac{\sum_{k=2}^n |\tilde{X}_{k,u}|\sum_{s=\max\{1,k-B_n\}}^{k-1} |\tilde{X}_{s,u}-\bar{X}_{s,u}|}{\sqrt{nB_n}}\right)=o(1).
	\]
	Similar arguments yield the same result for the second term of the right hand side of \cref{eq_tmp2_two_terms}. Note that
	\[
	\begin{split}
	&\EE\left(\max_{u\in\mathcal{U}}\frac{\sum_{k=2}^n |\tilde{X}_{u,k}|\sum_{s=\max\{1,k-B_n\}}^{k-1} |\tilde{X}_{u,s}-\bar{X}_{u,s}|}{\sqrt{nB_n}}\right)\\
	&\le\EE\left(\max_{u\in\mathcal{U}}\frac{\sum_{k=2}^n |\tilde{X}_{u,k}|\sum_{s=\max\{1,k-B_n\}}^{k-\ell} |\tilde{X}_{u,s}-\bar{X}_{u,s}|}{\sqrt{nB_n}}\right)\\
	&\qquad+
	\EE\left(\max_{u\in\mathcal{U}}\frac{\sum_{k=2}^n |\tilde{X}_{u,k}|\sum_{s=\max\{1,k-\ell+1\}}^{k-1} |\tilde{X}_{u,s}-\bar{X}_{u,s}|}{\sqrt{nB_n}}\right).
	\end{split}
	\]
	Applying \cref{lemma_Aven85} and using $\ell$-independence and H\"older's inequality, we have that uniformly on $u$
	\[
	\begin{split}
	&\EE\left(\sum_{k=2}^n |\tilde{X}_{u,k}|\sum_{s=\max\{1,k-B_n\}}^{k-\ell} |\tilde{X}_{u,s}-\bar{X}_{u,s}|\right)\\
	&=\EE\left(\sum_{k=2}^n |\tilde{X}_{u,k}|\right)\EE\left(\sum_{s=\max\{1,k-B_n\}}^{k-\ell} |\tilde{X}_{u,s}-\bar{X}_{u,s}|\right)
	\\
	&=\bigO(n)\bigO(B_n)\EE\left|\tilde{X}_{u,k}\Ind_{|\tilde{X}_{u,k}|>(nB_n)^{\alpha}}-\EE \tilde{X}_{u,k}\Ind_{|\tilde{X}_{u,k}|>(nB_n)^{\alpha}}\right|\\
	&\le\bigO(nB_n)(\EE\tilde{X}_{u,k}^p)^{1/p}(\Pr(|\tilde{X}_{u,k}|^p>(nB_n)^{\alpha p}))^{1-1/p}\\
	&=\bigO(nB_n)\bigO((nB_n)^{-\alpha p})^{1-1/p}
	=\bigO((nB_n)^{1-\alpha(p-1)}).
	\end{split}
	\]
	Furthermore, uniformly on $u$, we also have that
	\[
	\begin{split}
	&\sqrt{\EE\left(\sum_{k=2}^n |\tilde{X}_{u,k}|\sum_{s=\max\{1,k-B_n\}}^{k-\ell} |\tilde{X}_{u,s}-\bar{X}_{u,s}|\right)^2}\\
	&=\sqrt{\EE\left(\sum_{k=2}^n |\tilde{X}_{u,k}|\right)^2\EE\left(\sum_{s=\max\{1,k-B_n\}}^{k-\ell} |\tilde{X}_{u,s}-\bar{X}_{u,s}|\right)^2}\\
	&=\sqrt{\bigO(n^2)}\sqrt{\bigO(B_n^2)\EE\left|\tilde{X}_{u,k}\Ind_{|\tilde{X}_{u,k}|>(nB_n)^{\alpha}}-\EE \tilde{X}_{u,k}\Ind_{|\tilde{X}_{u,k}|>(nB_n)^{\alpha}}\right|^2}\\
	&\le\bigO(nB_n)(\EE\tilde{X}_{u,k}^p)^{1/p}(\Pr(|\tilde{X}_{u,k}|^p>(nB_n)^{\alpha p}))^{1-1/p}\\
	&=\bigO(nB_n)\bigO((nB_n)^{-\alpha p})^{1-1/p}
	=\bigO((nB_n)^{1-\alpha(p-1)}).
	\end{split}
	\]
	By the assumptions $p>4$ and  $(p-1)\alpha>3/4$, we have that
	\[
	\EE\left(\max_{u\in\mathcal{U}}\frac{\sum_{k=2}^n |\tilde{X}_{u,k}|\sum_{s=\max\{1,k-B_n\}}^{k-\ell} |\tilde{X}_{u,s}-\bar{X}_{u,s}|}{\sqrt{nB_n}}\right)=\bigO\left(\frac{C_n^{1/2}(nB_n)^{1-\alpha(p-1)}}{(nB_n)^{1/2}}\right)=o(1),
	\]
	since we have assumed $C_n^{1/2}=o[(nB_n)^{\alpha(p-1)-\frac{1}{2}}]$.
	Next, uniformly on $u$, the second term of the right hand side of \cref{eq_tmp2_two_terms} satisfies that
	\[
	\begin{split}
	&\EE\left(\sum_{k=2}^n |\tilde{X}_{u,k}|\sum_{s=\max\{1,k-\ell+1\}}^{k-1} |\tilde{X}_{u,s}-\bar{X}_{u,s}|\right)\\
	&=\bigO(n\ell)\EE\left|\tilde{X}_{u,k}^2\Ind_{\tilde{X}_{u,k}^2>(nB_n)^{2\alpha}}-\EE \tilde{X}_{u,k}^2\Ind_{\tilde{X}_{u,k}^2>(nB_n)^{2\alpha}}\right|\\
	&\le\bigO(n\ell)\left(\EE|\tilde{X}_{u,k}|^p\right)^{2/p}\left(\Pr(\tilde{X}_{u,k}^p<(nB_n)^{p\alpha})\right)^{1-2/p}\\
	&=\bigO(n\ell)\bigO((nB_n)^{-\alpha p})^{1-2/p}\\
	&=\bigO(n\ell)\bigO(nB_n)^{-\alpha(p-2)}.
	\end{split}
	\]
	Furthermore, uniformly on $u$, we have that
	\[
	\begin{split}
	&\sqrt{\EE\left(\sum_{k=2}^n |\tilde{X}_{u,k}|\sum_{s=\max\{1,k-\ell+1\}}^{k-1} |\tilde{X}_{u,s}-\bar{X}_{u,s}|\right)^2}\\
	&=\sqrt{\bigO(n^2\ell^2)\EE\left|\tilde{X}_{u,k}^2\Ind_{\tilde{X}_{u,k}^2>(nB_n)^{2\alpha}}-\EE \tilde{X}_{u,k}^2\Ind_{\tilde{X}_{u,k}^2>(nB_n)^{2\alpha}}\right|^2}\\
	&\le\bigO(n\ell)\left(\EE|\tilde{X}_{u,k}|^p\right)^{2/p}\left(\Pr(\tilde{X}_{u,k}^p<(nB_n)^{p\alpha})\right)^{1-2/p}\\
	&=\bigO(n\ell)\bigO((nB_n)^{-\alpha p})^{1-2/p}\\
	&=\bigO(n\ell)\bigO(nB_n)^{-\alpha(p-2)}.
	\end{split}
	\]
	Overall, we have that
	\[
	\EE\left(\max_{u\in\mathcal{U}}\frac{\sum_{k=2}^n |\tilde{X}_{u,k}|\sum_{s=\max\{1,k-\ell+1\}}^{k-1} |\tilde{X}_{u,s}-\bar{X}_{u,s}|}{\sqrt{nB_n}}\right)=\bigO\left(\frac{C_n^{1/2}n\ell (nB_n)^{-\alpha(p-2)}}{\sqrt{nB_n}}\right)=o(1),
	\]
	since we have assumed $C_n=o(B_n^{1+2\alpha(p-2)}n^{-2-2\gamma})$.

\subsection{Proof of \cref{lemma_max_dev_3}}\label{proof_lemma_max_dev_3}
	We prove this lemma by first showing that 
	\[\label{proof_lemma_max_dev_3_step1}
	\max_{u\in\mathcal{U}}\max_{0\le i\le B_n} \left|\frac{\sum_{j=1}^{k_n+1}V_j(u,\theta_i)}{\sqrt{nB_n}}\right|=o_{\Pr}(1),
	\]
	and then showing
	\[\label{proof_lemma_max_dev_3_step2}
	\max_{u\in\mathcal{U}}\max_{0\le i\le B_n}\left|\frac{\sum_{j=1}^{k_n}U_j(u,\theta_i)-\sum_{j=1}^{k_n}\bar{U}_j(u,\theta_i)}{\sqrt{nB_n}}\right|=o_{\Pr}(1).
	\]
	
	To show \cref{proof_lemma_max_dev_3_step1}, we note that $\{V_j\}$ are independent. Applying \cref{lemma_max_dev_pre1}, we have that
	\[
	\Pr\left(\left|\frac{\sum_{j=1}^{k_n+1}V_j}{\sqrt{nB_n}}\right|\ge \frac{1}{\log B_n}\right)\le C_1\left(\frac{\sum_{j=1}^{k_n+1}\EE V_j^2}{nB_n(\log B_n)^{-2}}\right)^Q+C_1\sum_{j=1}^{k_n+1}\Pr\left(\frac{|V_j|}{\sqrt{nB_n}}\ge \frac{C_2}{\log B_n}\right).
	\]
	Similar to the proof of \cref{lemma_max_dev_4}, one can show $\sum_{j=1}^{k_n+1}\EE V_j^2=\bigO(n^{1+\gamma}B_n^{1-\beta})$. Therefore, by choosing $\gamma$ close to zero and $Q$ large enough, we have that
	\[
	\left(\frac{\sum_{j=1}^{k_n+1}\EE V_j^2}{nB_n(\log B_n)^{-2}}\right)^Q=\bigO(n^{-c}),
	\]
	for any $c>0$. For the other term 
	\[
	\sum_{j=1}^{k_n+1}\Pr\left(\frac{|V_j|}{\sqrt{nB_n}}\ge \frac{C_2}{\log B_n}\right),
	\]
	we apply \cref{lemma_max_dev_pre2} with $M=(nB_n)^{\alpha}$, $k=B_n+\ell$, $\ell=\floor{n^{\gamma}}$ and $y=(\log B_n)^2$, which yields
	\[
	\begin{split}
	&\Pr\left(\frac{|V_j|}{\sqrt{nB_n}}\ge \frac{C_2}{\log B_n}\right)\\
	&\le 2\exp\left(-\frac{(\log B_n)^2}{4}\right)+\bigO\left(n^3(nB_n)^{2\alpha}\left(\frac{(
		\log B_n)^2}{nB_n}(\log B_n)^4 \floor{n^{3\gamma}}((nB_n)^{2\alpha}+B_n)\right)^Q\right)\\
	&\quad+\bigO\left(n^3(nB_n)^{2\alpha}\sum_{i=1}^n\Pr\left(|\bar{X}_{i,\ell}|\ge \frac{C_2 \frac{\sqrt{nB_n}}{\log B_n}}{(\log B_n)^2 \floor{n^{2\gamma}}((nB_n)^{\alpha}+(B_n+\floor{n^{\gamma}})^{1/2})}\right)\right),
	\end{split}
	\]
	where the second term of the right hand side is $\bigO(n^{-c})$ by choosing $Q$ large enough.  Since $\alpha<1/4$ and $|\bar{X}_{i,\ell}|<(nB_n)^{\alpha}$ almost surely, the last term of the right hand side converges to zero almost surely if 
	\[
	(nB_n)^{\alpha}=o\left(\frac{ \frac{\sqrt{nB_n}}{\log B_n}}{(\log B_n)^2 \floor{n^{2\gamma}}((nB_n)^{\alpha}+(B_n+\floor{n^{\gamma}})^{1/2})}\right),
	\]
	which can be satisfied by choosing $\gamma$ close enough to zero. Therefore, by choosing $Q$ large enough so that $\bigO(C_nB_n n^{-c})=o(1)$ (Note that this only requires $C_n=o(n^c B_n^{-1})$ for some $c$, which is always satisfied when $C_n$ is polynomial of $n$), we have that
	\[
	\Pr\left(\max_{u\in\mathcal{U}}\max_{0\le i\le B_n} \left|\frac{\sum_{j=1}^{k_n+1}V_j}{\sqrt{nB_n}}\right|\ge \frac{1}{\log B_n}\right)\le \bigO(C_nB_n) \Pr\left(\left|\frac{\sum_{j=1}^{k_n+1}V_j}{\sqrt{nB_n}}\right|\ge \frac{1}{\log B_n}\right)=o(1),
	\]
	which implies \cref{proof_lemma_max_dev_3_step1}.
	
	To prove \cref{proof_lemma_max_dev_3_step2}, note that 
	\[
	U_j-\bar{U}_j(u,\theta)=U_j(u,\theta)\Ind\left(\frac{|U_j(u,\theta)|}{\sqrt{nB_n}}>\frac{1}{(\log B_n)^4}\right)-\EE U_j(u,\theta)\Ind\left(\frac{|U_j(u,\theta)|}{\sqrt{nB_n}}>\frac{1}{(\log B_n)^4}\right).
	\]
	Therefore, other than using $p_n=B_n^{1+\beta}$ instead of $q_n=B_n+\ell$, the proof of \cref{proof_lemma_max_dev_3_step2} is essentially the same as the proof of \cref{proof_lemma_max_dev_3_step1}.

\subsection{Proof of \cref{lemma_max_dev_4}}\label{proof_lemma_max_dev_4}
	By \cref{lemma_bernstein},
	we have that
	\[
	\begin{split}
	&\Pr\left(\max_{u\in\mathcal{U}}\max_{i\notin [(\log B_n)^2, B_n-(\log B_n)^2]}\left|\frac{\sum_{j=1}^{k_n}\bar{U}_j}{\sqrt{nB_n}}\right|\ge x\sqrt{\log (B_nC_n)}\right)\\
	&= \bigO(C_n)\sum_{i\notin [(\log B_n)^2, B_n-(\log B_n)^2]}\Pr\left(\left|\frac{\sum_{j=1}^{k_n}\bar{U}_j}{\sqrt{nB_n}}\right|\ge x\sqrt{\log (B_nC_n)}\right)\\
	&=\bigO(C_nB_n+C_n(\log B_n)^2)
	\Pr\left(\frac{\sum_{j=1}^{k_n}|\bar{U}_j|}{\sqrt{nB_n}}\ge x\sqrt{\log (B_nC_n)}\right)\\
	&= \bigO(C_nB_n) \exp\left(\frac{-\frac{1}{2}x^2 nB_n (\log B_n+\log C_n)}{\sum_{j=1}^{k_n}\EE \bar{U}_j^2+\frac{1}{3}\frac{\sqrt{nB_n}}{(\log B_n)^4}x\sqrt{nB_n (\log B_n+\log C_n)}}\right).
	\end{split}
	\]
	Note that 
	$U_j=\sum_{k\in H_j}(\bar{Y}_{k,\ell}-\EE\bar{Y}_{k,\ell})$, 
	we first divide $\sum_{k\in H_j}(\bar{Y}_{k,\ell}-\EE\bar{Y}_{k,\ell})$ into $\ell$ sums of sub-sequences. Note that $\bar{Y}_{k,\ell}=\bar{X}_{k,\ell}\sum_{s=1}^{k-1}\alpha_{n,k-s}\bar{X}_{s,\ell}=\bar{X}_{k,\ell}\sum_{s=\max(1,k-B_n)}^{k-1}\alpha_{n,k-s}\bar{X}_{s,\ell}$. Thus, one can get $\|\bar{U}_j\|^2=\bigO(\ell B_n^2)$. Then using $\ell=\bigO(n^{\gamma})$ and $k_n=\floor{n/(p_n+q_n)}=\bigO(n/B_n^{1+\beta})$, one can get $\sum_{j=1}^{k_n}\EE\bar{U}_j^2=\bigO(n^{1+\gamma}B_n^{1-\beta})=o(nB_n)$ by choosing $\gamma$ and $\beta$ such that $n^\gamma B_n^{-\beta}=o(1)$. 
	
	Finally, we have that
	\[
	\begin{split}
	&\bigO(C_nB_n)\exp\left(\frac{-\frac{1}{2}x^2 nB_n (\log B_n+\log C_n)}{\sum_{j=1}^{k_n}\EE \bar{U}_j^2+\frac{1}{3}\frac{\sqrt{nB_n}}{(\log B_n)^4}x\sqrt{nB_n (\log B_n+\log C_n)}}\right)\\
	&=o\left[C_nB_n\exp\left(\frac{-\frac{1}{2}x^2\log(B_nC_n)}{o(nB_n)/(nB_n)+\frac{1}{3}x\frac{\log(B_nC_n)}{(\log B_n)^4}}\right)\right]\\
	&\to o\left[C_nB_n\exp\left(-\frac{3}{2}x(\log B_n)^4\right)\right]=o(1),
	\end{split}
	\]
	since $\log C_n+\log B_n=o(\log B_n)^4$ when $C_n$ and $B_n$ are polynomials of $n$.

\subsection{Proof of \cref{lemma_max_dev_5}}\label{proof_lemma_max_dev_5}
	\begin{enumerate}
		\item We first show that for $|i_1-i_2|\ge (\log B_n)^2/B_n$ or $|u_1-u_2|\ge \frac{n}{N}(1-1/(\log B_n)^2)$, we have that
		\[\label{proof_lemma_max_dev_eq1}
		\left|\frac{\EE \sum_{j=1}^{k_n}\bar{U}_j(u_{1},\theta_{i_1})\sum_{j=1}^{k_n}\bar{U}_j(u_{2},\theta_{i_2})}{nB_n}\right|=\bigO(1/(\log B_n)^2).
		\]
		Note that $\sum_{j}\bar{U}_j$ can be approximated by $\bar{g}_{n}$. This is because according to the proof of \cref{lemma_max_dev_3}, we have that
		\[
		\max_{u}\max_{i}\frac{\EE|\sum_{j=1}^{k_n}\bar{U}_j(u,\theta_i)-\bar{g}_{n}(u,\theta_i)|^2}{nB_n}=\bigO( B_n^{-\epsilon/2}).
		\]
		Next, we can approximate $\bar{g}_{n}$ by $\tilde{g}_{n}$. This is because by \cref{lemma_max_dev_2} we have that
		\[
		\max_{u}\max_{\theta} \frac{\EE|\tilde{g}_{n}(u,\theta)-\bar{g}_{n}(u,\theta)|^2}{nB_n}=\bigO(1/(\log B_n)^2).
		\]
		Finally, we only need to show that
		\[
		\frac{|\Cov(\tilde{g}_n(u_1,\theta_{i_1}),\tilde{g}_n(u_2,\theta_{i_2}))|}{nB_n}=\bigO(1/(\log B_n)^2),
		\]
		which has been proved in \cref{lemma_max_dev_pre5}{(i)}. 
		\item For convenience, we assume $\int a^2(t)\dee t=1$. Select $d$ distinct tuples $(\theta_{i_1},u_i), i=1,\dots,d$ that $(\log B_n)^2\le i_1\le \dots\le i_d\le B_n-(\log B_n)^2$ and $u_i\in\mathcal{U}, i=1,\dots,d$.
		Let $\bW_n=\sum_{j=1}^{k_n}W_j$ where
		\[
		W_j=\left(\frac{\bar{U}_j(u_1,\theta_{i_1})}{f(u_1,\theta_{i_1})},\dots,\frac{\bar{U}_j(u_d,\theta_{i_d})}{f(u_d,\theta_{i_d})}\right),
		\quad 1\le j\le k_n.
		\]
		Note that by \cref{lemma_max_dev_pre5}(iii), we have that
		\[
		\left|\EE\frac{\left(\sum_{j=1}^{k_n}\bar{U}_j(u,\theta)\right)^2}{nB_n}-4\pi^2f^2(u,\theta)\right|=\bigO(1/(\log B_n)^2).
		\]
		Together with \cref{proof_lemma_max_dev_eq1}, we have that
		\[
		\left|\frac{\Cov(\bW_n)}{nB_n}-4\pi^2\bI_d\right|=\bigO(1/(\log B_n)^2).
		\]
		Then we approximate $\bW_n$ by $\bW'_n=\sum_{j=1}^{k_n} W'_j$ using \cref{lemma_max_dev_pre6}, where $\{W'_j\}$ are independent centered normally distributed random vectors. 
		Then by \cref{lemma_max_dev_pre6}, we have $\Cov(
		W_j)=\Cov(W'_j)$, for $1\le j\le k_n$, and
		\[\label{proof_lemma_max_dev_eq2}
		\Pr\left(\frac{|\bW_n-\bW'_n|}{\sqrt{nB_n}}\ge 1/\log B_n\right)=\bigO(e^{-(\log B_n)^3}).
		\]
		
		Therefore, we have that
			\[\label{proof_lemma_max_dev_eq3}
			\left|\frac{\Cov(\bW'_n)}{nB_n}-4\pi^2\bI_d\right|=\bigO(1/(\log B_n)^2).
			\]
		\item Next, for $z=(z_1,\dots,z_d)$, we define the minimum of $\{z_i\}$ by $|z|_d:=\min_{1\le i\le d}\{z_i\}$. Then we show that
			\[\label{proof_lemma_max_dev_eq4}
			\Pr\left(\frac{|\bW_n|_d}{\sqrt{nB_n}}\ge y_n\right)=(1+o(1))\left(\sqrt{8\pi}y_n^{-1}\exp\left(-\frac{y_n^2}{8\pi^2}\right)\right)^d,
			\]
			uniformly on distinct tuples of $\{(u_j,\theta_{i_j}), j=1,\dots,d: (\log B_n)^2\le j_1\le\cdots\le j_d\le B_n-(\log B_n)^2, \frac{n}{2N}< u_j<1-\frac{n}{2N}\}$ such that for any two tuples $(u_{j_1},\theta_{i_{j_1}})$ and $(u_{j_2},\theta_{i_{j_2}})$, if $u_{j_1}=u_{j_2}$ then $|\theta_{i_{j_1}}-\theta_{i_{j_2}}|\ge (\log B_n)^2/B_n$; if $\theta_{i_{j_1}}=\theta_{i_{j_2}}$ then $|u_{j_1}-u_{j_2}|\ge \frac{n}{N}(1-1/(\log B_n)^2)$.
		
		According to \cref{proof_lemma_max_dev_eq2}, we have that
		\[
		\begin{split}
		&\Pr\left(\frac{|\bW'_n|_d}{\sqrt{nB_n}}\ge y_n-\frac{1}{\log B_n}\right)-\bigO(e^{-(\log B_n)^3})\\
		&\le \Pr\left(\frac{|\bW_n|_d}{\sqrt{nB_n}}\ge y_n\right)\le\Pr\left(\frac{|\bW'_n|_d}{\sqrt{nB_n}}\ge y_n-\frac{1}{\log B_n}\right)
		+\bigO(e^{-(\log B_n)^3}).
		\end{split}
		\]
		From \cref{proof_lemma_max_dev_eq3}, we have that
		\[
		\left|\frac{\Cov^{1/2}(\bW'_n)}{\sqrt{nB_n}}-2\pi\bI_d\right|=\bigO(1/(\log B_n)^2),
		\]
		so that for a standard normally distributed $R^d$-valued random vector, $\tilde{W}$, the tail probability of $\frac{\Cov^{1/2}(\bW'_n)}{\sqrt{nB_n}}\tilde{W}-2\pi\bI_d\tilde{W}$ satisfies that
		\[
		\begin{split}
		&\Pr\left(\left|\left(\frac{\Cov^{1/2}(\bW'_n)}{\sqrt{nB_n}}-2\pi\bI_d\right)\tilde{W}\right|\ge 1/\log B_n\right)\\
		&\quad=\bigO(e^{-(\log B_n)^2/4}).
		\end{split}
		\]
		Putting together the above results we can use $2\pi|\tilde{W}|_d$ (recall that we defined the minimum of $\{z_i\}$ by $|z|_d:=\min_{1\le i\le d}\{z_i\}$) instead of $\frac{|\bW'_n|_d}{\sqrt{nB_n}}$ to bound the tail probability of $\frac{|\bW_n|_d}{\sqrt{nB_n}}$:
		\[
		\begin{split}
		&\Pr(2\pi|\tilde{W}|_d\ge y_n-2/\log B_n)-\bigO(e^{-(\log B_n)^2/4})\\
		&\le\Pr\left(\frac{|\bW_n|_d}{\sqrt{nB_n}}\ge y_n \right)\\
		&\le\Pr(2\pi|\tilde{W}|_d\ge y_n-2/\log B_n)
		+\bigO(e^{-(\log B_n)^2/4}).
		\end{split}
		\]
		Using the following approximation of tail probability of a standard normally distributed random variable $Z$,
		\[
		\Pr(Z>z)=1-\Phi(z)\le \frac{1}{z\sqrt{2\pi}}\exp\left(-\frac{z^2}{2}\right),
		\]
		we can get that
		\[
		\Pr\left(|Z|>\frac{y_n}{2\pi}\right)=2\,\Pr\left(Z>\frac{y_n}{2\pi}\right) \le \sqrt{8\pi}y_n^{-1}\exp\left(-\frac{y_n^2}{8\pi^2}\right).
		\]
		Then we have shown that
		\[
		\Pr\left(\frac{|\bW_n|_d}{\sqrt{nB_n}}\ge y_n \right)=(1+o(1))\left(\sqrt{8\pi}y_n^{-1}\exp\left(-\frac{y_n^2}{8\pi^2}\right)\right)^d.
		\]
		
		Similarly, using \cref{lemma_max_dev_pre8} and \cref{lemma_max_dev_pre5}{(ii)}, we can also have that
			\[\label{proof_lemma_max_dev_eq5}
			\begin{split}
			&\Pr\left(\left|\frac{\sum_{j=1}^{k_n}\bar{U}_j(u_k,\theta_{i_k})}{\sqrt{nB_n}f(u_k,\theta_{i_k})}\right|\ge y_n , k=1,\dots,d\right)\\
			&\quad \le C \left(\sqrt{8\pi}y_n^{-1}\exp\left(-\frac{y_n^2}{8\pi^2}\right)\right)^{d-2}y_n^{-2}\exp\left(-\frac{y_n^2}{8\pi^2}(1+\delta)\right), 
			\end{split}
			\]
			for some $\delta>0$, 	uniformly on distinct tuples of $\{(u_j,\theta_{i_j}), j=1,\dots,d: (\log B_n)^2\le j_1\le\cdots\le j_d\le B_n-(\log B_n)^2, \frac{n}{2N}< u_j<1-\frac{n}{2N}\}$  such that for any two tuples $(u_{j_1},\theta_{i_{j_1}})$ and $(u_{j_2},\theta_{i_{j_2}})$with $j_1\le j_2$, if $u_{j_1}=u_{j_2}$ then if $\theta_{i_{j_1}}=\min_j \theta_{i_{j}}$ then $|\theta_{i_{j_1}}-\theta_{i_{j_2}}|\ge B_n^{-1}$; otherwise $|\theta_{i_{j_1}}-\theta_{i_{j_2}}|\ge (\log B_n)^2/B_n$; if $\theta_{i_{j_1}}=\theta_{i_{j_2}}$ then $|u_{j_1}-u_{j_2}|\ge \frac{n}{N}(1-1/(\log B_n)^2)$.
		\item Finally, we define
		\[
		A_{u,i}=\left\{\frac{|\sum_{j=1}^{k_n}\bar{U}_j(u,\theta_i)|^2}{4\pi^2n B_n f^2(u,\theta_i)}\ge 2\log B_n+2\log C_n-\log(\pi\log B_n+\pi\log C_n)+x\right\}
		\]
		and we show
		\[
		\Pr\left(\bigcup_{(\log B_n)^2\le i\le B_n-(\log B_n)^2, u \in\mathcal{U}} A_{u,i}\right)\to 1-e^{-e^{-x/2}}.
		\]
		To this end, we define
		\[
		\tilde{A}_u=\bigcup_{(\log B_n)^2\le i\le B_n-(\log B_n)^2} A_{u,i}
		\]
		and
		\[
		P_{t,u}:=\sum_{(\log B_n)^2\le i_1<\dots<i_t\le B_n-(\log B_n)^2} \Pr(A_{u,i_1}\cap\dots\cap A_{u,i_t}).
		\]
		Then by Bonferroni's inequality, we have for every fixed $k$ and $u$ that
		\[
		\sum_{t=1}^{2k}(-1)^{t-1}P_{t,u}\le \Pr(\tilde{A}_u)\le \sum_{t=1}^{2k-1}(-1)^{t-1}P_{t,u}.
		\]
		Next following the proof of  \cite[Theorem]{Watson1954} 
		and \cite[Theorem 3.3]{Woodroofe1967} based on \cref{proof_lemma_max_dev_eq4} and \cref{proof_lemma_max_dev_eq5}, we can show that
		\[
		P_{t,u}\to [B_n\Pr(A_{u,i})]^t/t!
		\]
		as $n\to \infty$. 
		As shown in \cite[pp.799]{Watson1954}, with  \cref{proof_lemma_max_dev_eq4} and \cref{proof_lemma_max_dev_eq5}, when $n\to\infty$, we have that
		\[
		P_{t,u}\to[(B_n-2(\log B_n)^2)^t/t!+\bigO(B_n-2(\log B_n)^2)^{t-1}]\Pr(A_{u,i})^t.
		\]
		Therefore, we have shown that
		\[
		\Pr(\tilde{A}_u)\to 1-e^{-[B_n \Pr(A_{u,i})]}.
		\]
		Finally, we use the above techniques again to show
		\[
		\Pr\left(\bigcup_{u\in\mathcal{U}}\tilde{A}_u\right)\to 1-e^{-e^{-x/2}},
		\]
		which means we only need to show
		\[
		C_n\Pr(\tilde{A}_u)\to \exp(-x/2).
		\]
		Letting $y_n^2/4\pi^2=2\log B_n+2\log C_n-\log(\pi\log B_n+\pi\log C_n)+x$, as in \cref{proof_lemma_max_dev_eq4}, we have that
		\[
		\begin{split}
		&C_n\Pr(\tilde{A}_u)\to C_nB_n \Pr(A_{u,i})\to C_n B_n\Pr\left(|N|>\frac{y_n}{2\pi}\right)\\
		&\to \frac{C_nB_n}{y_n}\sqrt{8\pi}\exp\left(-\frac{y_n^2}{8\pi^2}\right)\\
		&\to C_n B_n\frac{\sqrt{8\pi}}{\sqrt{8\pi^2}\sqrt{\log B_n+\log C_n}}\exp\left(-\frac{x}{2}\right)\frac{\sqrt{\pi \log B_n+\pi \log C_n}}{B_nC_n}\\
		&\to \exp\left(
		-\frac{x}{2}\right).
		\end{split}
		\] 
	\end{enumerate}

\subsection{Proof of \cref{remark_expectation_consistency}}\label{proof_remark_expectation_consistency}
First of all, by the assumption $\GMC(2)$
\[
\begin{split}
\EE\hat{f}_n(u,\theta)-f(u,\theta)&=\frac{1}{2\pi}\left[\sum_{k=-B_n}^{B_n}\EE\hat{r}(u,k)a(k/B_n)-\sum_{k\in\mathbb{Z}}r(u,k)\right]\exp(\sqrt{-1}k\theta)\\
&=\frac{1}{2\pi}\sum_{k=-B_n}^{B_n}\left[\EE\hat{r}(u,k)a(k/B_n)-r(u,k)\right]\exp(\sqrt{-1}k\theta)+\bigO(\rho^{B_n}).
\end{split}
\]
Next, by the SLC condition, we know $r(u,k)$ is Lipschitz. Together with the Lipschitz condition of $\tau(\cdot)$, we have that
\[
\EE\hat{r}(u,k)&=\frac{1}{n}\sum_{i=1}^N \tau\left(\frac{i-\lfloor uN \rfloor}{n}\right)\tau\left(\frac{i+k-\lfloor uN \rfloor}{n}\right)\EE(X_i X_{i+k})\\
&=\frac{1}{n}\sum_{i=\lfloor uN \rfloor-\frac{n}{2}}^{\lfloor uN \rfloor+\frac{n}{2}} \tau\left(\frac{i-\lfloor uN \rfloor}{n}\right)\tau\left(\frac{i+k-\lfloor uN \rfloor}{n}\right)\left[r(i/N,k)+\bigO(k/N)\right]\\
&=\frac{1}{n}\sum_{i=\lfloor uN \rfloor-\frac{n}{2}}^{\lfloor uN \rfloor+\frac{n}{2}} \left[\tau\left(\frac{i-\lfloor uN \rfloor}{n}\right)^2+o(k/n)\right]r(i/N,k)+\bigO(k/N).
\]
Since $r(u,k)$ is twice continuously differentiable with respect to $u$, we have that
\[
\EE\hat{r}(u,k)&=\frac{1}{n}\sum_{i=\lfloor uN \rfloor-\frac{n}{2}}^{\lfloor uN \rfloor+\frac{n}{2}} \tau\left(\frac{i-\lfloor uN \rfloor}{n}\right)^2 \left[r(u,k)+\left(\frac{i-\lfloor uN\rfloor }{N}\right)\frac{\partial r(u,k)}{\partial u}+\bigO(n^2/N^2)\right]\\
&\qquad +o(k/n)r(i/N,k)+\bigO(k/N).
\]
Furthermore, since $\tau(\cdot)$ is an even function
\[
\frac{1}{n}\sum_{i=\lfloor uN \rfloor-\frac{n}{2}}^{\lfloor uN \rfloor+\frac{n}{2}} \tau\left(\frac{i-\lfloor uN \rfloor}{n}\right)^2\left(\frac{i-\lfloor uN\rfloor }{N}\right)\frac{\partial r(u,k)}{\partial u}=0.
\]
Therefore, we have that
\[
\EE\hat{r}(u,k)&=\left[\int\tau^2(x)\dee x+o(1/n)\right]r(u,k)+\bigO(n^2/N^2)+o(k/n)r(u,k)+\bigO(k/N)\\
&=r(u,k)+o(k/n+1/n)r(u,k)+\bigO(k/N+n^2/N^2).
\]
Therefore, by the locally quadratic property of $a(\cdot)$ at $0$, we have that
\[
\begin{split}
&\EE\hat{r}(u,k)a(k/B_n)-r(u,k)\\
&=\EE\hat{r}(u,k)\left[a(0)+a'(0)k/B_n+\frac{1}{2}a''(0)k^2/B_n^2+o(k^2/B_n^2)\right]-r(u,k)\\
&=-C\left(\frac{k^2}{B_n^2}+o(k/n)\right)r(u,k)+\bigO(k/N+n^2/N^2).
\end{split}
\]
Then, using the fact that if $\theta\notin \{0,\pi\}$, we know that
\[
\sum_{k=0}^{B_n}\cos(k\theta)=\frac{1}{2}+\frac{\sin(\frac{2B_n+1}{2}\theta)}{2\sin(\theta/2)},\quad \sum_{k=1}^{B_n}\sin(k\theta)=\frac{\sin\frac{B_n\theta}{2}\sin\frac{(B_n+1)\theta}{2}}{\sin (\theta/2)}.
\]
Then, for fixed $\theta\notin\{0,\pi\}$, we have that
\[
\sum_{k=0}^{B_n}\cos(k\theta)=\bigO(1),\quad \sum_{k=0}^{B_n}k\cos(k\theta)=\bigO(B_n).
\]
If $\sup_u\sum_{k\in \mathbb{Z}}|r(u,k)|k^2<\infty$ and $B_n=o(n)$, then
\[
\EE\hat{f}_n(u,\theta)-f(u,\theta)+  \frac{C}{2\pi}\sum_{k\in \mathbb{Z}}\frac{k^2r(u,k)\exp(\sqrt{-1}k\theta)}{B_n^2}=\bigO(B_n/N+n^2/N^2).
\]
Finally, $B_n=o(N^{1/3})$ implies $\bigO(B_n/N)=o(1/B_n^2)$. Also, $n=o(N^{2/3})$ and $B_n=o(N^{1/3})$ implies $\bigO(B_n^2n^2/N^2)=o(1)$.

\subsection{Proof of \cref{thm_near_optimality}}\label{proof_thm_near_optimality}
For simplicity, we denote $\delta_{u,n}$ as  $\delta_u$ and $\delta_{\theta,n}$ as $\delta_{\theta}$. First, we write
\[
\begin{split}
\hat{f}_n(u,\theta)-\hat{f}_n(u_i,\theta_j)&=\hat{f}_n(u,\theta)-\hat{f}_n(u_i,\theta_j)\\
&-\EE[\hat{f}_n(u,\theta)-\hat{f}_n(u_i,\theta_j)]+\EE[\hat{f}_n(u,\theta)-\hat{f}_n(u_i,\theta_j)].
\end{split}
\]
Then by continuity we have that
\[
\max_{\{u_i,\theta_j\}}\sup_{\{u: |u-u_i|\le \delta_u, \theta: |\theta-\theta_j|\le \delta_{\theta}\}}|\EE\hat{f}_n(u,\theta)-\EE\hat{f}_n(u_i,\theta_j)|=o_{\Pr}(\sqrt{\log n}).
\]

Letting $\hat{g}_n(u,u_i,\theta,\theta_j):=\hat{f}_n(u,\theta)-\hat{f}_n(u_i,\theta_j)$, it suffices to show that
\[
\max_{\{u_i,\theta_j\}}\sup_{\{u: |u-u_i|\le \delta_u, \theta: |\theta-\theta_j|\le \delta_{\theta}\}}|\hat{g}_n(u,u_i,\theta,\theta_j)-\EE\hat{g}_n(u,u_i,\theta,\theta_j)|=o_{\Pr}(\sqrt{\log n}).
\]
Note that
\[
\begin{split}
\hat{g}_n(u,u_i,\theta,\theta_j)=&[\hat{f}_n(u,\theta)-\EE\hat{f}_n(u,\theta)]\left[1-\frac{\hat{f}_n(u_i,\theta_j)}{\hat{f}_n(u,\theta)} \right]\\
&+\EE\hat{f}_n(u,\theta)\left[1-\frac{\hat{f}_n(u_i,\theta_j)}{\hat{f}_n(u,\theta)} \right].
\end{split}
\]
Then we can write
\[
\begin{split}
&\sup_{\{u,\theta\}} \hat{g}_n(u,u_i,\theta,\theta_j)\\
&\le\sup_{\{u,\theta\}} \left[\hat{f}_n(u,\theta)-\EE\hat{f}_n(u,\theta)\right] \sup_{\{u,\theta\}}\left|\frac{\hat{f}_n(u_i,\theta_j)}{\hat{f}_n(u,\theta)}-1\right|\\
&+ \sup_{\{u,\theta\}} \EE\hat{f}_n(u,\theta)\sup_{\{u,\theta\}}\left|\frac{\hat{f}_n(u_i,\theta_j)}{\hat{f}_n(u,\theta)}-1\right|.
\end{split}
\]
Since by \cref{thm_consistency}, we have that
\[
\sup_{\{u,\theta\}} [\hat{f}_n(u,\theta)-\EE\hat{f}_n(u,\theta)]=\bigO_{\Pr}(\sqrt{\log n}).
\]
Therefore, the following result completes the proof.
\begin{lemma}\label{lemma_near_optimality} If $\delta_u=\bigO(\frac{n}{N(\log n)^{\alpha}})$ and $\delta_{\theta}=\bigO(\frac{1}{B_n(\log n)^{\alpha}})$ for some $\alpha>0$, then
	\[
	\max_{\{u_i,\theta_j\}}\sup_{\{u: |u-u_i|\le \delta_u, \theta: |\theta-\theta_j|\le \delta_{\theta}\}}\left|\frac{\hat{f}_n(u_i,\theta_j)}{\hat{f}_n(u,\theta)}-1\right|=o_{\Pr}(1).
	\]
\end{lemma}
\begin{proof}
	See  \cref{proof_lemma_near_optimality}.
\end{proof}

\subsection{Proof of \cref{lemma_near_optimality}}\label{proof_lemma_near_optimality}
First, we pick any $(u_0,\theta_0)$ such that $|u_0-u|\le \delta_u$ and $|\theta_0-\theta|\le \delta_{\theta}$. Then 
\[
\hat{f}_n(u_0,\theta_0)-\hat{f}_n(u,\theta)=\frac{1}{2\pi} \sum_{k=-B_n}^{B_n}a(k/B_n)[\hat{r}(u_0,k)\exp(\sqrt{-1}k\theta_0)-\hat{r}(u,k)\exp(\sqrt{-1}k\theta)].
\]
Using
$\tau\left(\frac{i-\lfloor u_0N\rfloor}{n}\right)=\tau\left(\frac{i-\lfloor uN\rfloor}{n}\right)+\bigO\left(\frac{\delta_uN}{n}\right)$,
we have that
\[
&\hat{r}(u_0,k)\exp(\sqrt{-1}k\theta_0)
=\frac{1}{n}\sum_{i=1}^N \tau\left(\frac{i-\lfloor u_0N \rfloor}{n}\right)\tau\left(\frac{i+k-\lfloor u_0N \rfloor}{n}\right)(X_i X_{i+k})\exp(\sqrt{-1}k\theta_0)\\
&=\frac{1}{n}\sum_{i=\lfloor uN\rfloor-\frac{n}{2}}^{\lfloor uN\rfloor+\frac{n}{2}} \left[\tau\left(\frac{i-\lfloor uN \rfloor}{n}\right)\tau\left(\frac{i+k-\lfloor uN \rfloor}{n}\right)+\bigO\left(\frac{\delta_uN}{n}\right)\right](X_i X_{i+k})\exp(\sqrt{-1}k\theta_0).
\]
Note that
 $\exp(\sqrt{-1}k\theta_0)=\exp(\sqrt{-1}k\theta)[\exp(\sqrt{-1}k(\theta_0-\theta))]$ and  $\cos(k\theta_0)=\cos(k\theta)\cos(k(\theta_0-\theta))-\sin(k\theta)\sin(k(\theta_0-\theta))$.  Therefore, we have that
\[
&\hat{f}_n(u_0,\theta_0)=\frac{1}{2\pi} \sum_{k=-B_n}^{B_n}a(k/B_n)\hat{r}(u_0,k)\exp(\sqrt{-1}k\theta)\exp(\sqrt{-1}k(\theta_0-\theta))\\
&=\frac{1}{2\pi} \sum_{k=-B_n}^{B_n}a(k/B_n)\hat{r}(u,k)\cos(k\theta)\left[1+ \bigO\left(\frac{\delta_u N}{n}\right)\right]\left[1+\bigO(k\delta_{\theta})\right]\\
&\quad -\frac{1}{2\pi} \sum_{k=-B_n}^{B_n}a(k/B_n)\hat{r}(u_0,k)\sin(k\theta)\bigO(k\delta_{\theta})\\
&=\hat{f}_n(u,\theta)\left[1+ \bigO\left(\frac{\delta_u N}{n}\right)\right]\left[1+\bigO(B_n\delta_{\theta})\right]+\bigO_{\Pr}(B_n\delta_{\theta}),
\]
where we have used the fact that the GMC condition implies
$\sum_{k=0}^{B_n}k r(u,k)=\bigO(\sum_{k=0}^{B_n}k\rho^k)=\bigO(B_n)$. Note that we have assumed that $f(u,\theta)>f_*>0$ uniformly over $u$ and $\theta$, so we can write $\bigO_{\Pr}(B_n\delta_{\theta})=(B_n\delta_{\theta})\bigO_{\Pr}(\hat{f}_n(u,\theta))$. Therefore, we have that
\[
\hat{f}_n(u_0,\theta_0)-\hat{f}_n(u,\theta)=\bigO(\delta_u N/n+B_n\delta_{\theta})\bigO_{\Pr}(\hat{f}_n(u,\theta)),
\]
which implies that
\[
\left|\frac{\hat{f}_n(u_0,\theta_0)}{\hat{f}_n(u,\theta)}-1\right|=\bigO_{\Pr}(\delta_uN/n+B_n\delta_{\theta}).
\]
In order to make it equal to $o_{\Pr}(1)$, we only need $\delta_u=o(n/N)$ and $\delta_{\theta}=o(1/B_n)$. Therefore, choosing $\alpha>0$, $\delta_u=\bigO\left(\frac{n}{N(\log n)^{\alpha}}\right)$ and $\delta_{\theta}=\bigO\left(\frac{1}{B_n(\log n)^{\alpha}}\right)$ is sufficient.

\subsection{Proof of \cref{remark_SLC}}\label{proof_remark_SLC}
By the triangle inequality and H\"older's inequality, we have that
\[
&|r(u,k)-r(s,k)|\\
&=\left|\EE\left[G(u,\mathcal{F}_{i})G(u,\mathcal{F}_{i+k})-G(s,\mathcal{F}_{i})G(s,\mathcal{F}_{i+k})
\right]\right|\\
&\le \|\left[G(u,\mathcal{F}_{i})-G(s,\mathcal{F}_{i})\right]G(u,\mathcal{F}_{i+k})\|_1+\|\left[G(u,\mathcal{F}_{i+k})-G(s,\mathcal{F}_{i+k})\right]G(s,\mathcal{F}_{i})\|_1\\
&\le\|G(u,\mathcal{F}_{i})-G(s,\mathcal{F}_{i})\|_q \|G(u,\mathcal{F}_{i+k})\|_p +\|G(u,\mathcal{F}_{i+k})-G(s,\mathcal{F}_{i+k})\|_q\|G(s,\mathcal{F}_{i})\|_p\\
&\le C|u-s|.
\]


\end{document}